\documentclass[11pt]{amsart}

\usepackage{DKstyle}
\usepackage{bm}
\usepackage{float}
\usepackage[toc,page]{appendix}

\newcommand{\col}{\: : \:}
\newcommand{\Mod}[1]{\ (\mathrm{mod}\ #1)}
\newcommand{\vertiii}[1]{{\left\vert\kern-0.25ex\left\vert\kern-0.25ex\left\vert #1
    \right\vert\kern-0.25ex\right\vert\kern-0.25ex\right\vert}}
 \makeatletter
\newcommand*{\rom}[1]{\expandafter\@slowromancap\romannumeral #1@}
\makeatother
\theoremstyle{plain}

\usepackage{caption}
\usepackage[labelfont=rm]{subcaption}

\usepackage[pagebackref=false, colorlinks]{hyperref}

\hypersetup{pdffitwindow=true,linkcolor=blue,citecolor=blue,urlcolor=blue,menucolor=blue}

\usepackage{comment}

\def\d{\text{d}}

\subjclass{11F12, 11F55, 11D45}%

\date{\today}%
\dedicatory{}%
\commby{}%

\title{Fourier expansion of light-cone Eisenstein series}
\author{Dubi Kelmer}
\address{Department of Mathematics, Boston College, Chestnut Hill MA 02467-3806, USA}
\email{kelmer@bc.edu}

\author{Shucheng Yu}
\address{Department of Mathematics, Uppsala University, Box 480, SE-75106, Uppsala, SWEDEN}
\email{shucheng.yu@math.uu.se}
\address{Current address:  School of Mathematical Sciences, University of Science and Technology of China (USTC), 230026, Hefei, China}
\email{yusc@ustc.edu.cn}

\thanks{D.K. and S.Y. were partially supported by NSF CAREER grant DMS-1651563. S.Y. was supported by the Knut and Alice Wallenberg Foundation}
\begin{document}
\begin{abstract}
In this work we give an explicit formula for the Fourier coefficients of Eisenstein series corresponding to certain arithmetic lattices acting on hyperbolic $n+1$-space. 
As a consequence we obtain results on location of all poles of these Eisenstein series as well as their supremum norms. We use this information to get new results on counting rational points on spheres.
\end{abstract}
\maketitle

\tableofcontents

\section{Introduction}
This work is concerned with the Fourier expansion of Eisenstein series corresponding to a non-uniform arithmetic lattice, $\G$, acting on the hyperbolic space $\bH^{n+1}$. 
When $n=1$ and $\Gamma=\SL_2(\Z)$ or a congruence subgroup, the analytic continuation of the Eisenstein series can be obtained directly via its explicit Fourier expansion.
Similar results were also obtained for $\G=\SL_2(\cO)$ acting on $\bH^3$ with $\cO$ an order in a quadratic (imaginary) number field \cite{ElstrodtGrunewaldMennicke1985}, as well as for certain lattices corresponding to orders in quaternion algebras acting on $\bH^4$ \cite{Gritsenko1987,Krieg1988,GritsenkoSchulzePillot1994} and $\bH^5$ \cite{KrafftOsenberg1990}. 
In this paper we describe a novel approach for obtaining an explicit Fourier expansion for arithmetic lattices acting on $\bH^{n+1}$ which allows us to treat all dimensions at once.
Using our expansion we obtain precise information regarding the location of all poles of these Eisenstein series, estimates on their sup-norms as well as some interesting applications for counting rational points on spheres.

\begin{rem}
An alternative approach to obtain the analytic properties of Eisenstein series for arithmetic lattices, uses automorphic representations following the ideas of Langlands. Using this approach it was shown in \cite{Reznikov1993}, that the constant term of Eisenstein series corresponding to congruence lattices acting on $\bH^{n+1}$ is always a ratio of entire functions of order one, however, pinning down {precisely what those functions are} and using this to extract information on the location of the poles is not obvious. 
\end{rem}

We now describe our main result in more detail. 
Given $n\in \N$ and  $Q$ a rational quadratic form of signature $(n+1,1)$, let $G=\SO_Q^+(\R)$ denote the connected component of the orthogonal group of $Q$. There is a natural (transitive) action of $G$ on the upper half space
$$\bH^{n+1}=\left\{z=(\bm{x},y)\in \R^{n}\times \R^+\right\}$$ 
with the stabilizer of $(\bm{0},1)$ a maximal compact group $K$, allowing us to identify the quotient $G/K$ with the upper half space $\bH^{n+1}$ (see 
 \rmkref{rem:Coordinates} below).
 The group of integer points $\Gamma=\SO^+_Q(\Z)$ is a lattice and we identify $\G\bk \bH^{n+1}$ with $\G\bk G/K$.  To simplify the discussion we will assume for the rest of the introduction that our quadratic form is 
\begin{equation}\label{e:Qn}
Q_n(\bm{v})=\sum_{j=1}^{n+1}v_j^2-v_{n+2}^2.\end{equation}

Next we define a series we call the {\em light-cone Eisenstein series} defining a function on $\Gamma\bk G/K\cong \G\bk \bH^{n+1}$. For this consider the light cone  
$$\cV_{Q_n}=\left\{\bm{v}\in \R^{n+2}\setminus\{\bm{0}\}: Q_n(\bm{v})=0\right\},$$ 
 and let $\cV_{Q_n}^+=\{\bm{v}\in \cV_{Q_n}:  v_{n+2}>0\}$ denote the positive half branch. 
The (completed) light-cone Eisenstein series is defined for any $g\in G$ and  $\Re(s)>n$ by the series
$$E^*_{Q_n}(s,g)=\sum_{\bm{v}\in \cV_{Q_n}^+(\Z)} \|\bm{v} g\|^{-s},$$
where $\cV_{Q_n}^+(\Z)=\cV_{Q_n}^+\cap \Z^{n+2}$ is the set of integer lattice points in $\cV_{Q_n}^+$ and $\|\cdot\|$ is the Euclidean norm on $\R^{n+2}$.
We then define the light-cone Eisenstein series, $E_{Q_n}(s,g)$, by the equation $E^*_{Q_n}(s,g)=2^{{-}s/2}\zeta(s)E_{Q_n}(s,g)$ (which can also be expressed as a sum over primitive integer points on the light cone). 

Note that for  $Q_n$, the maximal compact subgroup $K$ preserves $\|\cdot\|$; see section \ref{sec:explicit}.
It is then clear that both  $E^*_{Q_n}(s,g)$ and $E_{Q_n}(s,g)$ are right $K$- and left $\G$-invariant and hence are functions on 
$\G\bk G/K$. Using the identification $\G\bk G/K\cong \G\bk \bH^{n+1}$ we can also think of $E^*_{Q_n}(s,g)$ and $E_{Q_n}(s,g)$ as functions on $\G\bk \bH^{n+1}$,
and with a slight abuse of notation we write $E_{Q_n}^*(s,z)$ and $E_{Q_n}(s,z)$ for the corresponding function on the upper half space.

As we will show below, when $\G$ has one cusp we have that,  $E_{Q_n}(s,z)$ equals the Eisenstein series corresponding to this cusp (with the appropriate scaling). When there are several cusps, it is given by the following weighted sum over the Eisenstein series of all cusps,
\begin{equation}\label{e:EQ2E}E_{Q_n}(s,z)=\sum_{j=1}^\kappa v_{P_j}^{\frac{s}{n}} E_{\G,j}(s,z),\end{equation}
where $E_{\G,j}(s,z)$ is the Eisenstein series and $P_j$ is the parabolic subgroup corresponding to the $j$th cusp, and $v_{P_j}$ is the volume of that cusp (see section \ref{sec:cusps} for more details and the definition of these cusp volumes). 
 
 \subsection{Fourier expansion}
 The advantage of working with the series $E_Q(s,z)$,  is that it is possible to get an explicit formula for Fourier coefficients  in terms of products of the zeta functions and Dirichlet $L$-functions. 
Explicitly, for any integer $D\not\equiv 3\Mod{4}$ let $\left(\frac{D}{\cdot }\right)$ denote the Kronecker symbol, let $\chi_D$ denote the unique primitive Dirichlet character satisfying that $\chi_D(j)=\left(\frac{D}{j}\right)$ for all $j$ co-prime to $D$, and let $L(s,\chi_D)$ denote the corresponding $L$-function. We also denote by 
$$
\xi(s)=\pi^{-s/2}\G\left(\tfrac{s}{2}\right)\zeta(s)\quad \textrm{and }\quad L^*(s,\chi_{-4})=\left(\tfrac{\pi}{4}\right)^{-\frac{s+1}{2}} \G\left(\tfrac{s+1}{2}\right)L(s,\chi_{-4}),
$$
the corresponding completed zeta and $L$-functions.  
 \begin{Thm}\label{thm:mainthmintr}
For $Q=Q_n$ as in \eqref{e:Qn}, the function $E_{Q}(s,z)$ is invariant under $\bm{x}\mapsto \bm{x}+\bm{v}$ for all $\bm{v}\in \Lambda$ with 
 $$\Lambda:=\{\bm{v}\in \Z^n: \|\bm{v}\|^2\equiv 0\Mod{2}\}.$$ 
 The normalized function $E_{Q}(s,z)=\frac{2^{s/2}}{\zeta(s)}E_{Q}^*(s,z)$  has a Fourier expansion 
 $$E_{Q}(s,z)=y^s+\Phi_n(s)y^{n-s}+ \frac{2^{s-n+1}\pi^s}{\G(s)\zeta(s)}\sum_{\bm\lambda\in \Lambda^*\setminus\{\bm{0}\}}\|\bm{\lambda}\|^{s-\frac{n}{2}}\Phi_{n}(s;\bm{\lambda})y^{\frac{n}{2}}K_{s-n/2}(2\pi\|\bm{\lambda}\|y))e(\bm{\lambda}\cdot \bm{x}),$$
 where $\Lambda^*$ is the dual lattice of $\Lambda$, $\G(s)$ is the gamma function and $K_{s}(y)$ is the Bessel function of the second kind.
 Furthermore
 $$\Phi_{n}(s):=\e^{(2)}_{n}(s)\left\lbrace\begin{array}{ll}
		\frac{\Gamma(\frac{2s-n+2}{4})^2}{\G(\frac{s+1}{2})\G(\frac{s-n+1}{2})} \frac{\xi(s-n+1)\xi(s-\frac{n}{2})}{\xi(s)\xi(s-\frac{n}{2}+1)} & n\equiv 0\Mod{4},\\
		\frac{\G(\frac{2s-n}{4})\G(\frac{2s-n+4}{4})}{\G(\frac{s+1}{2})\G(\frac{s-n+1}{2}) }\frac{ \xi(s-n+1)L^*(s-\frac{n}{2},\chi_{-4})}{\xi(s) L^*(s-\frac{n}{2}+1,\chi_{-4})} & n\equiv 2\Mod{4},\\
		\frac{\G(\frac{2s-n+1}{4})\G(\frac{2s-n+3}{4})}{\G(\frac{s+1}{2}) \Gamma(\frac{s-n+1}{2})}\frac{\xi(s-n+1)\xi(2s-n)}{\xi(s)\xi(2s-n+1)} &n\equiv 1\Mod{2},
	\end{array}\right.
	$$
	and for any nonzero $\bm{\lambda}\in \Lambda^*$,
	\begin{align*}
		\Phi_{n}(s;\bm{\lambda})&:=\e_n(s; \bm{\lambda})\left\lbrace\begin{array}{ll}
			\frac{1}{\zeta(s-\frac{n}{2}+1)} & n\equiv 0\Mod{4},\\
			\frac{1}{L(s-\frac{n}{2}+1,\chi_{-4})} & n\equiv 2\Mod{4},\\
			\frac{L\left(s-\frac{n-1}{2},\chi_{D}\right)}{\zeta(2s-n+1)} & n\equiv 1\Mod{2},
		\end{array}\right.
	\end{align*}
with $D:=(-1)^{\frac{n-1}{2}}\|2\bm{\lambda}\|^2$ when $n$ is odd. Here
	\begin{align}\label{equ:e2factor}
		\e_{n}^{(2)}(s):= \left\lbrace\begin{array}{ll}
			\frac{2^{s-\frac{n}{2}-1}-\left(1-(-1)^{\frac{n}{4}}\right)2^{-1}}{1-2^{\frac{n}{2}-1-s}} & n\equiv 0\Mod{4},\\
			2^{s-\frac{n}{2}} & n\equiv 2\Mod{4},\\
			\frac{2^{s-\frac{n+1}{2}}+(-1)^{\frac{n^2-1}{8}}}{1+(-1)^{\frac{n^2-1}{8}}2^{\frac{n-1}{2}-s}}  & n\equiv 1\Mod{2},
		\end{array}\right.
	\end{align}
	and $\e_n(s; \bm{\lambda})$ is some (explicitly computable) function which is holomorphic off the line $\Re(s)=\left \lfloor{\frac{n-1}{2}}\right \rfloor$. 
 \end{Thm}
 
 \begin{remark}\label{rmk:morepreefac}
 Our analysis will give us more precise information on these $\e$-factors, $\e_n(s;\bm{\lambda})$, corresponding to nonzero $\bm{\lambda}\in \Lambda^*$. Indeed, we will see that $\e_n(s;\bm{\lambda})=\prod_{p\mid 2\|2\bm{\lambda}\|^2}\e_n^{(p)}(s; \bm{\lambda}) $ with each $\e_n^{(p)}(s;\bm{\lambda})$ a polynomial in $p^{-s}$ except when $p=2$ and $n\not\equiv 2\Mod{4}$ in which case 
\begin{align*}
 \left\lbrace\begin{array}{ll}
			(1-2^{\frac{n}{2}-1-s})\e_n^{(2)}(s;\bm{\lambda}) & n\equiv 0\Mod{4},\\
			(1-2^{n-1-2s})\e_n^{(2)}(s;\bm{\lambda})  & n\equiv 1\Mod{2},
				\end{array}\right.
\end{align*} 
is a polynomial in $2^{-s}$; see \rmkref{rmk:polyargu}. In particular, this immediately implies that $\e_n(s;\bm{\lambda})$ is holomorphic off the line $\Re(s)=\left \lfloor{\frac{n-1}{2}}\right \rfloor$. Moreover, the potential poles contributed from these $\e$-factors are periodic points on the line $\Re(s)=\left \lfloor{\frac{n-1}{2}}\right \rfloor$ (solutions to $1-2^{2\left \lfloor{\frac{n-1}{2}}\right \rfloor-2s}=0$).
We also note that these factors are closely related to twisted divisor functions, see \rmkref{rmk:efacaritfun}. 
 \end{remark}
 Next we describe some interesting consequences of these formulas.

  \subsection{Poles of $E_{Q_n}(s,z)$}
   From the general theory of Eisenstein series corresponding to a lattice $\G$ acting on hyperbolic space, 
 the Eisenstein series corresponding to each cusp has an analytic continuation to the complex plane, where in the half plane $\Re(s)\geq \frac{n}{2}$ it has a simple pole at $s=n$ with constant residue $\omega_\G = \frac{1}{\vol(\G\bk \bH^{n+1})}$, and possibly finitely many simple poles in the interval $(\tfrac{n}{2},n)$ that we call \textit{exceptional poles}. The residues of these poles (if they exist) are residual forms, which are square integrable eigenfunctions of the Laplacian, with eigenvalue $\lambda_s=s(n-s)$. There are also infinitely many poles in the half plane $\Re(s)<\tfrac{n}{2}$ (all located in some vertical strip) whose residues are not square integrable.
 
 When $n=1$ and $\Gamma\leq \SL_2(\Z)$ is a congruence subgroup it is well known that there are no exceptional poles \cite[Theorem 11.3]{Iwaniec2002}, and was conjectured by Selberg that in fact all eigenfunctions of the Laplacian in $L^2(\G\bk \bH^{n+1})$ satisfy $\lambda\geq \frac{1}{4}$. The situation is similar for congruence subgroups $\G$ acting on $\bH^3$, where it is known that there are no exceptional poles  \cite{ElstrodtGrunewaldMennick1998} and conjectured that all eigenvalues satisfy $\lambda\geq 1$.  In his 1962 ICM lecture \cite[p. 183]{Selberg1962}, Selberg raised the question of the existence of exceptional poles for arithmetically defined lattices acting on higher dimensional spaces, and we revisit this question now since 
having Eisenstein series with no exceptional poles is very useful in some recent applications to Diophantine approximations (see \cite{KelmerYu2022b}). 
We note that for $n\geq 3$, it was shown in \cite{CogdellLiPiatetskiShapiroSarnak}, that there are congruence groups for which $\lambda_s=s(n-s)$ is a Laplace eigenvalue for any $s\in (\tfrac{n}{2},n)\cap \Z$. However, it is not clear whether these eigenvalues come from residues of Eisenstein series and if they can appear at full level.

Using {the Fourier expansion in \thmref{thm:mainthmintr}} we can read off the precise location of all exceptional poles. In particular we have the following result which was our main motivation in computing the Fourier coefficients.
 \begin{Thm}\label{c:poles1}
For $n=1$ or any $n\not\equiv 1\Mod{8}$ the light-cone Eisenstein series $E_{Q_n}(s,z)$ has no exceptional poles in the interval $(\tfrac{n}{2},n)$ while for $n\geq 9$ with $n\equiv 1\Mod{8}$ there is one such pole at $s=\tfrac{n+1}{2}$.
 In particular, for $1\leq n\leq 7$ where $\G=\SO^+_{Q_n}(\Z)$ has one cusp, the Eisenstein series $E_{\Gamma}(s,z)$ has no exceptional poles. 
 \end{Thm}

 In addition to the exceptional poles (or lack thereof) in the half plane $\Re(s)>\tfrac{n}{2}$, using our formula we can also read off the location of all poles in the half plane $\Re(s)<\tfrac{n}{2}$ in terms of the non-trivial zeros of the Riemann zeta function and a certain Dirichlet $L$-function.
 \begin{Thm}\label{c:poles2}
For $n\geq 2$, except for (possibly) finitely many real poles  and a set of periodic poles on the line $\Re(s)=\left \lfloor{\frac{n-1}{2}}\right \rfloor$, the set of poles of $E_{Q_n}(s, z)$ in the half plane $\Re(s)<\tfrac{n}{2}$ is given by 
 \begin{align*}
 \left\lbrace\begin{array}{ll}
	\left\{\rho, \tfrac{n}{2}-1+\rho: \xi(\rho)=0 \right\}	 & n\equiv 0\Mod{4},\\
	\left\{\rho, \tfrac{n}{2}-1+\rho': \xi(\rho)=0,\ L^*(\rho',\chi_{-4})=0\right\} & n\equiv 2\Mod{4},\\
	\left\{\rho, \tfrac{n-1+\rho}{2}: \xi(\rho)=0\right\}	 & n\equiv 1\Mod{2}.
	\end{array}\right.
 \end{align*}
For $n=1$ the poles of $E_{Q_1}(s, z)$ in the half plane $\Re(s)<\tfrac{1}{2}$ are the zeroes of $(1+2^{-s})\xi(2s)$.
 \end{Thm}

 \begin{rem}
 When there is a single cusp, all poles come from the constant term, however, when there is more than one cusp there could potentially be poles where the residue of the constant term vanishes, so to see all poles one needs to look at all coefficients.  
 Moreover, when there is more than one cusp, the Eisenstein series corresponding to the different cusps could have additional exceptional poles that get cancelled out when summing over all cusps to form $E_{Q_n}(s,z)$.
 While we could not prove, or exclude, that such cancellation occurs, we note that for the poles in the half plane $\Re(s)<\tfrac{n}{2}$ there is a lot of cancellation.
To see this we note that by \cite[Theorem 1.3]{Kelmer2015}, for a lattice with $\kappa$ cusps, if $\rho=\beta+i\gamma$ are all the poles of the Eisenstein series (for all cusps) then for any $T>2$,
$$\sum_{|\gamma|\leq T}\left(\frac{n}{2}-\beta\right)=\frac{\kappa n}{2\pi}T\log(T)+O(T).$$
On the other hand, when summing over the poles of the light-cone Eisenstein series described above, we get  $\frac{n}{2\pi}T\log(T)+O(T)$.
 Comparing the two asymptotics we see that when there is more than one cusp many poles do indeed get canceled in the sum \eqref{e:EQ2E}.
\end{rem}

 \subsection{Volume and cusp volume}
 Let $\omega_{Q_n}=\Res_{s=n}E_{Q_n}(s,z)$  which, by our formula is given by 
  \begin{equation}\label{eq:omegaQ}
  \omega_{Q_n}=\left\lbrace\begin{array}{ll}
		\frac{2^{\frac{n}{2}-1}-\left(1-(-1)^{\frac{n}{4}}\right)2^{-1}}{1-2^{-1-\frac{n}{2}}}\frac{\Gamma(\frac{n+2}{4})^2}{\G(\frac{n+1}{2})\G(\frac{1}{2})} \frac{\xi(\frac{n}{2})}{\xi(n)\xi(\frac{n}{2}+1)} & n\equiv 0\Mod{4},\\
		2^{\frac{n}{2}} \frac{\G(\frac{n}{4})\G(\frac{n}{4}+1)}{\G(\frac{n+1}{2})\G(\frac{1}{2}) }\frac{ L^*(\frac{n}{2},\chi_{-4})}{\xi(n) L^*(\frac{n}{2}+1,\chi_{-4})} & n\equiv 2\Mod{4},\\
		\frac{2^{\frac{n-1}{2}}+(-1)^{\frac{n^2-1}{8}}}{1+(-1)^{\frac{n^2-1}{8}}2^{-\frac{n+1}{2}}} \frac{\G(\frac{n+1}{4})\G(\frac{n+3}{4})}{\G(\frac{n+1}{2}) \Gamma(\frac{1}{2})}\frac{1}{\xi(n+1)} & n\equiv 1\Mod{2}\ \text{and}\ n>1,\\
		\frac{4}{\pi} & n=1.
	\end{array}\right.\end{equation}
 
For $\G=\SO_{Q_n}^+(\Z)$ using the relation \eqref{e:EQ2E}, between the light-cone Eisenstein series and the Eisenstein series corresponding to all cusps, together with the fact that $\Res_{s=n}E_{\G, i}(s,z)=\frac{1}{\vol(\G\bk \bH^{n+1})}$  we see that 
 \begin{equation}\label{eq:omegaQvol}
 \sum_{j=1}^{\kappa} v_{P_j}=\vol(\G\bk \bH^{n+1})\omega_{Q_n}.
 \end{equation} 
For the cusp at infinity $P=P_1$ we can calculate the cusp volume explicitly and show that $v_P=2^{2-n}(n!)^{-1}$ (see \rmkref{rmk:cuspvol1}). The volume $\vol(\G\bk \bH^{n+1})$ was already calculated in \cite[Section 3]{RatcliffeTschantz1997}.
For $n\leq 7$, there is one cusp and our formula \eqref{eq:omegaQ} together with \eqref{eq:omegaQvol}  recovers the volume formula in \cite{RatcliffeTschantz1997}. For $n\geq 8$ there is more than one cusp, and our result gives new information on the volume of the remaining cusps. To illustrate this, in the following table we compare  $\vol(\G\bk \bH^{n+1})\omega_{Q_n}$ to $v_{P_1}$ for $1\leq n\leq 12$.  Below we abbreviate $L(s,\chi_{-4})$ simply by $L(s)$.
\begin{displaymath}
\renewcommand{\arraystretch}{1.5}
\begin{tabular}{| c | c | c | c | c|}
  \hline			
  $n$ & $\vol(\G\bk \bH^{n+1})$ & $\omega_{Q_n}$ &  $\vol(\G\bk \bH^{n+1})\omega_{Q_n}$ & $v_{P_1}$\\
  \hline
  $1$ &  $\frac{\pi}{2}$ & $\frac{4}{\pi}$ & $2$ & $2$\\
  \hline
  $2$ & $\frac{L(2)}{6}$ & $\frac{3}{L(2)}$ & $\frac12$ & $\frac12$\\
  \hline
  $3$ & $\frac{\pi^2}{720}$ & $\frac{60}{\pi^2}$ & $\frac{1}{12}$ & $\frac{1}{12}$\\
  \hline
  $4$ & $\frac{7\zeta(3)}{7680}$ & $\frac{80}{7\zeta(3)}$ & $\frac{1}{96}$ & $\frac{1}{96}$\\
  \hline
  $5$ & $\frac{\pi^3}{388800}$ & $\frac{405}{\pi^3}$ & $\frac{1}{960}$ & $\frac{1}{960}$\\
  \hline
$6$ & $\frac{L(4)}{181440}$ & $\frac{63}{4L(4)}$ &  $\frac{1}{11520}$ & $\frac{1}{11520}$\\
  \hline
$7$ & $\frac{17\pi^4}{4572288000}$ & $\frac{28350}{17\pi^4}$ & $\frac{1}{161280}$ & $\frac{1}{161280}$\\
\hline
$8$ & $\frac{527\zeta(5)}{22295347200}$ & $\frac{512}{31\zeta(5)}$ & $\frac{17}{43545600}$ & $\frac{1}{2580480}$\\
\hline
$9$ & $\frac{\pi^5}{164602368000}$ & $\frac{16065}{4\pi^5}$ & $\frac{17}{696729600}$  & $\frac{1}{46448640}$\\
\hline
$10$ & $\frac{L(6)}{5748019200}$ & $\frac{165}{16 L(6)}$ & $\frac{1}{557383680}$  & $\frac{1}{928972800}$\\
\hline
$11$ & $\frac{691\pi^6}{31070342983680000}$ & $\frac{20945925}{2764 \pi^6}$ & $\frac{31}{183936614400}$ & $\frac{1}{20437401600}$\\
\hline
$12$ & $\frac{87757 \zeta(7)}{24485642108928000}$ & $\frac{515840}{87757 \zeta(7)}$ & $\frac{31}{1471492915200}$ & $\frac{1}{490497638400}$ \\
  \hline
 \end{tabular}
\end{displaymath}
\begin{remark}
We note that the volume listed above is twice the volume given in \cite[p.\ 56]{RatcliffeTschantz1997}. This is because there they considered the lattice of orientation-preserving integral isometries (not necessarily of determinant one) which contains $\G$ as an index $2$ subgroup.
 \end{remark}
 \begin{remark}
 We see from the above table that when $1\leq n\leq 7$, $\vol(\G\bk \bH^{n+1})\omega_{Q_n}$ agrees with $v_{P_1}$, confirming that in these cases $\G$ has only one cusp and giving an alternative calculation for the volume. For $n=8, 9$, it follows from the analysis in \cite{Vinberg1972}  that $\G$ has two cusps. Our computation then implies that in these two cases the second cusp volume $v_{P_2}=\frac{1}{348364800}$. We also note that in all cases we computed we have that $\vol(\G\bk \bH^{n+1})\omega_{Q_n}$ is rational and it seems likely that the cusp volumes are always rational.
 \end{remark}

 \subsection{Functional equation}
 When $n\leq7$ the lattice $\G=\SO^+_{Q_n}(\Z)$ has one cusp and so $E_{Q_n}(s,z)=E_\G(s,z)$ is the Eisenstein series corresponding to this single cusp, and hence,  satisfies the functional equation 
 \begin{equation}\label{e:funeq}
 E_{Q_n}(n-s,z)=\Phi_n(n-s)E_{Q_n}(s,z).\end{equation}
 When $n\geq 8$ there is more than one cusp, and such a functional equation is no longer automatic. In fact, when $n\equiv 0\Mod{8}$ we have that $\Phi_n(s)\Phi_{n}(n-s)\neq 1$ so \eqref{e:funeq} can not hold; see \rmkref{rmk:fueqfail} below.
 Nevertheless, we show the following surprising result.
 \begin{Thm}\label{thm:funequ}
 For any $n\not\equiv 0\Mod{8}$ the light-cone Eisenstein series $E_{Q_n}(s,z)$ satisfies \eqref{e:funeq}. 
 \end{Thm}
 
 \subsection{Evaluation at special points}
 Similar to the Eisenstein series for $\SL_2(\Z)$, which can be evaluated at a special point to give the standard  Epstein zeta function which is a ratio of zeta functions and $L$-functions, it is possible to interpret the light-cone Eisenstein series at the point $z_0=(\bm{0},1)$ in terms of other arithmetic functions.
 Explicitly, let $r_k(m)$ denote  the number of representations of a positive integer $m$ as a sum of $k$ squares and let 
 $$R_k(s)=\sum_{q=1}^{\infty} \frac{r_k(q^2)}{q^s}.$$
 Noting that any $\bm{v}\in \cV^+_{Q_n}(\Z)$ is of the form $\bm{v}=(\bm{u},q)$ with $\bm{u}\in \Z^{n+1}$ 
 and $q\in \N$ and that $\|\bm{v}\|^2=2q^2$ we see that the completed light-cone Eisenstein series satisfies
 $$E^*_{Q_n}(s,z_0)=2^{-s/2}R_{n+1}(s),$$
 (and hence $E_{Q_n}(s,z_0)=\frac{R_{n+1}(s)}{2^s\zeta(s)}$).
 For small values of $k$ there are explicit formulas for $r_k(q^2)$ and it is possible to identify the series $R_k(s)$. 
Explicitly, for $k=3$ we have the following formula of Hurwitz (see e.g. \cite[p. 237]{Duke2003})
 $$R_3(s)=6(1-2^{1-s})\frac{\zeta(s)\zeta(s-1)}{L(s,\chi_{-4})}.$$
 For even values of $2\leq k\leq 8$ there are well known explicit formulas for $r_k(m)$ in terms of divisor functions 
 \cite[p. 186-187]{Iwaniec1997}, and one can use them to express $R_k(s)$ as a ratio of zeta functions and Dirichlet $L$-functions. Explicitly, doing this for $k=2,4,6,8$ we get 
 $$R_2(s)=\frac{4\zeta(s)^2L(s,\chi_{-4})}{(1+2^{-s})\zeta(2s)},$$
 $$R_4(s)=8(1-2^{2-s})\frac{\zeta(s-1)\zeta(s-2)\zeta(s)}{\zeta(2s-2)},$$
 $$ R_6(s)=\frac{12 L(s-2,\chi_{-4})\zeta(s-4)\zeta(s)}{(1-2^{2-s})\zeta(2s-4)},$$
 and 
 $$R_8(s)=16\frac{1 +3\cdot 2^{1-s} +2^{7-2s}}{1+2^{3-s}}\frac{\zeta(s-3)\zeta(s-6)\zeta(s)}{\zeta(2s-6)}.$$

Comparing the value at this point with the Fourier expansion we get new formulas, analogous to the Selberg-Chowla formula \cite{SelbergChowla1967},
 expressing some simple combinations of zeta functions and Dirichlet $L$-functions as an infinite series involving Bessel functions and divisor functions. For larger values of $k$ and for odd $k$ the formula for $r_k(m)$ includes Fourier coefficients of modular forms and is no longer as explicit. Hence, it is not clear whether such explicit formulas still hold in general.

 \subsection{Supremum norms}
 Another application of our results regards the problem of the supremum norm of the light-cone Eisenstein series $E_{Q_n}(\tfrac{n}{2}+it,z)$ for $z$ restricted to some compact set. 
 
We define the \textit{supremum norm exponent} of $E_{Q_n}(s,g)$ to be the smallest real number $\nu_n$ so that for any $\nu>\nu_n$ and any compact $\Omega\subseteq \G\bk \bH^{n+1}$ there is a constant $C=C(\Omega,\nu)$ such that for all $z\in \Omega$ and $t\in\R$
\begin{equation}\label{e:supnorm}|E_{Q_n}(\tfrac{n}{2}+it,z)| \leq C(|t|+1)^\nu.\end{equation}
We also define the $L^2$-exponent, $\tilde\nu_{n}$ by the smallest real number such that for any $\nu>\tilde{\nu}_n$ and for any $T>0$,
\begin{equation}\label{e:supnormave}\int_{-T}^T|E_{Q_n}(\tfrac{n}{2}+it,z)|^2\,\d t \leq C T^{{2\nu+1}}.\end{equation}
We note that from general results on Eisenstein series we have that \eqref{e:supnormave} holds with $\nu=\tfrac{n}{2}$; see \rmkref{r:WboundE} below. 
\begin{rem}
These exponents are related to the problem of bounding the sup-norm of eigenfunctions of the Laplacian in terms of the eigenvalues (recall that $E_{Q_n}(\tfrac{n}{2}+it,z)$ is an eigenfunction of the hyperbolic Laplacian with eigenvalue $\tfrac{n^2}{4}+t^2$). For a compact Riemannian manifold, $X$, one can show that if $\phi\in L^2(X)$ has $\|\phi\|_2=1$ and eigenvalue $\lambda$ then $\|\phi\|_{\infty}\ll \lambda^{\frac{\dim(X)-1}{4} }$. While this bound is sharp in general, when $X$ has negative curvature it is believed that this can be improved and there are some results of this nature for cusp forms on some arithmetic hyperbolic manifolds \cite{IwaniecSarnak1995,BlomerHarcosMilicevic2016}.  
Here the upper bound $\nu_n\leq \tfrac{n}{2}$ is usually referred to as the \textit{convexity bound} and any improvement is called \textit{sub-convex}. While for the $L^2$-exponent the convexity bound $\tilde\nu_n\leq \tfrac{n}{2}$ is known to hold in general, for the supremum norm, even the convexity bound is not known except for some small values of $n$.
Explicitly,  for $n=1$ the first sub-convex bound of $\nu_1\leq \tfrac{3}{8}$ was proved by Young \cite{Young2018} which was improved by Blomer \cite{Blomer2020}  to $\nu_1\leq \frac{1}{3}$, it is conjectured that in fact $\nu_1=0$ and was recently shown in \cite{KelmerKontorovichLutsko2023}  that $\tilde\nu_1=0$.
For $n=2$ the work of Assing \cite{Assing2019}, implies the sub-convex bound $\nu_2\leq \frac{7}{8}$ which was recently improved in  \cite{KelmerKontorovichLutsko2023} who showed that $\nu_2=\tfrac{1}{2}$.  \end{rem}

By using the explicit formulas for the Fourier coefficients as well as the functional equation, we prove the following lower and upper bounds on these exponents for all values of $n$.
\begin{Thm}\label{t:supnorm}
For any {even} $n\geq 2$ we have that $\tfrac{n}{2}-1\leq \tilde\nu_{n}\leq \nu_n\leq \tfrac{n}{2}$ {while for odd $n\geq 3$ we have $\tfrac{n}{2}-1\leq \tilde\nu_{n}\leq \nu_n\leq \tfrac{3n}{4}$}.
\end{Thm}

\begin{rem}
When $n\not\equiv 0\pmod{8}$ the upper bound follows from the functional equation by the standard application of the Phragmén–Lindelöf convexity principle.
It is curious that for odd $n\geq 3$, this method gives the bound $\nu_n\leq \tfrac{3n}{4}$ instead of the expected convexity bound of $\tfrac{n}{2}$ which is obtained when $n$ is even.
When $n\equiv 0\pmod{8}$, we no longer have the functional equation, however, the expected convexity bound $\nu_n\leq \tfrac{n}{2}$ can still be obtained via a more direct application of the Fourier expansion. 
\end{rem}

\begin{rem}
Using the evaluation of $E_{Q_n}(s,z_0)$ in terms of $R_{n+1}(s)$ when $n\in\{1, 2,3,5,7\}$, we can say more about the size of $E_{Q_n}(\frac{n}{2}+it,z_0)$ for $z_0=(\bm{0},1)$. Explicitly, using the known estimates for the mean square of zeta functions \cite{HardyLittlewood1916}, Dirichlet $L$-functions \cite{Rane1980} and Dedekind zeta functions \cite{Muller1989} on the critical line, it follows that for odd $n\leq 7$ 
we have for any $\e>0$ 
\begin{equation}\label{eq:meannoundz0}
T^{n-\epsilon}\ll_\e \int_{-T}^T|E_{Q_n}(\tfrac{n}{2}+it,z_0)|^2\ \d t\ll_\e T^{n+\epsilon}. 
\end{equation}
In particular, in these cases we can get the stronger lower bound $\tilde\nu_{n}\geq \frac{n-1}{2}$. 
For $n=2$ we have $|E_{Q_2}(1+it,z_0)|= 3\left|(1-2^{-it})\frac{\zeta(it)}{L(1+it,\chi_{-4})}\right|$ and from the known upper and lower bounds on $\zeta(s)$ and $L(s,\chi_{-4})$ for $\Re(s)=1$  we get the pointwise upper bound
\begin{align*}
|E_{Q_2}(1+it,z_0)|\ll_\e |t|^{\frac12+\epsilon}.
\end{align*}
We also note that it is not hard to show that \eqref{eq:meannoundz0} still holds for $n=2$ implying that also $\tilde\nu_2\geq \frac{1}{2}$.
\end{rem}

\subsection{Counting rational points on the sphere}
 Next we give a couple of interesting applications of our result to the problem of counting rational points on the sphere of bounded height. 
 Let $S^n=\{\bm{x}\in \R^{n+1}: \|\bm{x}\|=1\}$ be the unit $n$-sphere.
 For a large parameter $T$ we are interested in counting 
 $$N(T)=\#\left\{\tfrac{\bm{p}}{q}\in S^n: (\bm{p},q)\in \Z^{n+2}_{\rm pr},\ 1\leq q\leq T\right\}.$$
Here $\Z^{n+2}_{\rm pr}$ is the set of primitive points in $\Z^{n+2}$.
 We also consider a smoothed counting function given by 
 $$N_h(T)=\sum_{\substack{\frac{\bm{p}}{q}\in S^n\\ (\bm{p},q)\in \Z^{n+2}_{\rm pr}}}h\left(\tfrac{q}{T}\right),$$
where $h\in C^\infty_c(\R^+)$ is a smooth compactly supported function. We have the following smoothed counting estimate.
 \begin{Thm}\label{t:smoothcounting}
 For any $h\in C^\infty_c(\R^+)$ and any $k>0$ we have 
  $$N_h(T)= \omega_{Q_n}\hat{h}(-n)T^n+c_n\hat{h}\left(-\tfrac{n+1}{2}\right)T^{\frac{n+1}{2}}+O_{h,k}(T^{\frac{n}{2}}\log(T)^{-k}),$$
 where $\hat{h}(s)=\int_0^\infty h(y)y^{-(s+1)}\ \d y$ is the Mellin transform of $h$ and $c_n=0$ unless $n>1$ satisfies  $n\equiv 1\Mod{8}$. Moreover, when $n$ is even a similar estimate holds with a remainder of order  $O_h(T^{\frac{n}{2}-\alpha})$ {(and respectively $O_h(T^{\frac{n-\alpha}{2}})$ when $n$ is odd)} for some $\alpha<1/2$ if and only if the Riemann zeta function $\zeta(s)$ (respectively the Dirichlet L-function $L(s,\chi_{-4})$ if $n\equiv 2\Mod{4}$ or both if $n=2$) has no zeros in $\Re(s)\geq 1-\alpha$. 
 \end{Thm}
 
 \begin{rem}
 {As we show in \propref{pr:AC2SC} below, the improved bound on the remainder of order $O_h(T^{\frac{n}{2}-\alpha})$ is equivalent to having an analytic continuation for $E_Q(s,z)$ with no poles in the strip $\tfrac{n}{2}-\alpha<\Re(s)\leq \tfrac{n}{2}$. Given the explicit description of poles in \thmref{c:poles2}, such a pole-free strip is equivalent to a zero-free region for the corresponding zeta and $L$-functions.}
 \end{rem}

 \begin{rem}
 This problem can also be studied via the circle method  \cite{Brown1996,Getz2018,Tran2020} where the results of  \cite{Brown1996} imply that for $n>2$
  $$N_h(T)= c_nT^n+O_{h,\e}(T^{\left \lfloor{\frac{n+2}{2}}\right \rfloor+\epsilon}),$$
 while its refinements in \cite{Getz2018,Tran2020} imply that 
 $$N_h(T)= c_nT^n+c_n'  T^{\left \lfloor{\frac{n+2}{2}}\right \rfloor}+O_{h,\e}(T^{\frac{n}{2}+\epsilon}),$$
 with secondary term $c_n'=0$ when $n\equiv 2\Mod{4}$. Comparing these results with ours we see that in this case the  secondary term $c_n'$ in fact vanishes for all $n\not\equiv 1\Mod{8}$, and that one can improve the remainder 
 $O_{h,\e}(T^{\frac{n}{2}+\epsilon})$ to $O_{h,k}(T^{\frac{n}{2}}\log(T)^{-k})$ (or even to $O_{h,\e}(T^{\frac{n-1}{2}+\e })$ assuming the Riemann Hypothesis).
 \end{rem}

By localizing the smoothing function it is possible to translate the smoothed counting estimate to an estimate with a sharp cutoff. In this case, one cannot expect such strong bounds for the remainder. In particular, noting that when $T$ is a half integer there are no rational points with denominator in $[T, T+1/4]$, the best bound we can hope for is $O(T^{n-1})$. Nevertheless, using bounds on the size of the light-cone Eisenstein series restricted on the line $\Re(s)=\tfrac{n}{2}$ we get the following. 
\begin{Thm}\label{t:supnorm2counting}
Let $\nu>0$ be such that \eqref{e:supnormave} holds for $z_0=(\bm{0},1)$. {If $n\geq 9$ with $n\equiv 1\Mod{8}$ further assume that $\nu\geq \frac{1}{n-1}$}. Then
$$\#\left\{\tfrac{\bm{p}}{q}\in S^n: (\bm{p},q)\in \Z^{n+2}_{\rm pr},\ 1\leq q\leq T\right\}=\frac{\omega_{Q_n}}{n} T^n+
O(T^{n-\frac{n}{2(\nu+1)}}).$$
\end{Thm}
\begin{rem}
The convexity bound $\tilde{\nu}_{n}\leq \frac{n}{2}$ gives the bound of $O(T^{n-\frac{n}{n+2}})$ for the remainder.
We note that one can consider this problem for a general integral quadratic form and the same result holds in this generality (see \propref{p:supnorm2counting} below).
Moreover, any improvement over the convexity bound for the (averaged) sup-norm will result in an improved bound on the remainder.
In particular, the bounds \eqref{eq:meannoundz0} 
imply that for $n\in \{1,2,3,5,7\}$ we have a better estimate of $O_\e(T^{n-\frac{n}{n+1}+\e})$ for the remainder.
\end{rem}


\begin{rem}
In addition to the circle method, such counting problems can be studied by other methods such as spectral theory of automorphic forms \cite{DukeRudnickSarnak1993}, the ergodic theory of group actions \cite{GorodnikNevo2012},  and elementary theory of modular forms \cite{Duke2003}. The results of Duke \cite{Duke2003}, using properties of modular form, not only counts these rational points, but also show that they become equidistributed on the sphere. We note that using an appropriate generalization of the light-cone Eisenstein series (obtained by applying certain raising operators as in \cite{Yu17}), it is also possible to get estimates for counting rational points lying in a prescribed subsets of the sphere. For $n=1$ this was done in the work of \cite{BurrinNevoRuhrWeiss2020}, and for larger $n$ in \cite{KelmerYu2022b}. \end{rem}

\subsection{General quadratic forms}
The light-cone Eisenstein series  $E_Q^*(s,z)$ and $E_Q(s,z)$ can be defined for a general rational form $Q$  (see section \ref{sec:LCES} below),
 and our approach for calculating the Fourier coefficients can be applied more generally. For example we carry this out explicitly for the forms
\begin{equation}\label{e:Qnd}
Q_{n,d}(\bm{v})=\sum_{j=1}^{n+1}v_j^2-d^2v_{n+2}^2,\end{equation}
with $d\in \N$ odd and square-free and a slight modification of these argument also apply to the form 
$\tilde{Q}_{n,d}(\bm{v})=\sum_{j=2}^{n+1}v_j^2+2dv_{1}v_{n+2}$ (for any $d\in\N$).
While some of the results (such as the functional equation) only hold when $d=1$, other results 
 (such as the precise description of all poles) hold for this larger family of forms.

While it is possible to carry out such calculations also for other families of forms, it is not clear to us if it is possible to do it for a general rational form. Nevertheless,  we conclude with another argument that gives very strong limitations on the possible exceptional poles of $E_Q(s,z)$, without using the explicit calculation of the Fourier coefficients. Instead, we use the relation between the analytic properties of $E_Q(s,g)$ and the problem of counting  the number of integer solutions to $Q(v)=0$ with bounded norm.
 Such counting problems were studied (in smooth form) using the circle method in  \cite{Brown1996,Getz2018,Tran2020} and  their results have the following striking implication.
 
 \begin{Thm}\label{t:circle}
Assume $n\geq 3$. For any rational form $Q(\bm{v})=\bm{v}J\bm{v}^t$ of signature $(n+1,1)$, the corresponding light-cone Eisenstein series $E_Q(s,g)$ has at most one exceptional pole at $\sigma=\left \lfloor{\frac{n+2}{2}}\right \rfloor$. Moreover, if $n\equiv 0\Mod{4}$ and  $\sqrt{-\det(J)}$ is irrational, or if $n\equiv 2\Mod{4}$, then there are no exceptional poles. 
\end{Thm} 

\begin{rem}
The results of \cite{Getz2018,Tran2020} apply to indefinite rational quadratic forms of any signature, and even to quadratic forms defined over number fields. Using similar arguments it is likely that one can obtain results on exceptional poles of corresponding Eisenstein series in these cases as well. We leave the treatment of this more general application to future work.  We also note that it is possible for $E_Q(s,z)$ to have no exceptional poles even when the conditions above are not satisfied, as is the case for the form $Q_{n,d}$.
\end{rem}

\begin{rem}
Examining the formula for the constant term for $Q_n$ (and more generally for $Q_{n,d}$) we see when $n\equiv 2\pmod{4}$ the ratio of zeta functions have no poles in $(\tfrac{n}{2},n)$, however, when $n\equiv 0\pmod{4}$ there is a pole at $s=\tfrac{n+2}{2}$ coming from the pole of $\zeta(s-\frac{n}{2})$ which is canceled by a zero of one of the $\Gamma$-factors. This is also the case for the potential pole at $s=\tfrac{n+1}{2}$ coming from the pole of $\zeta(2s-n)$ when $n\equiv 3\pmod{4}$. This seems to suggest that such a cancelation will hold for more general quadratic forms. However, when  $n\equiv 5\pmod{4}$ the pole at $\tfrac{n+1}{2}$ is cancelled by the local factor at $2$, a phenomena that is less likely to hold in general. Indeed, similar calculation for the form $\tilde{Q}_{n,d}$ show that when $d$ is even and $n\equiv 1\pmod{4}$ this pole is no longer canceled and $E_{\tilde{Q}_{n,d}}(s,z)$ has a pole at $s=\tfrac{n+1}{2}$.
\end{rem}
\subsection{Outline of paper}
For the reader's convenience we end this introduction with a brief outline of this paper. After section \ref{sec:prelim} which covers some necessary preliminaries and background on Eisenstein series, in section \ref{sec:LCES} we define the light-cone Eisenstein series corresponding to a general rational quadratic form and show how it is related to the standard Eisenstein series. We then prove some general applications to counting integer and primitive points on the corresponding light cone, and use them to prove \thmref{t:circle}. Next, in section \ref{s:Fourier} we state the formulas for Fourier coefficients for the more general case of the form $Q_{n,d}$ and show how the results regarding the poles, functional equation and sup-norms follow from these {formulas}. Finally in section \ref{sec:FourierCalc} we give the proof of  the formulas for the Fourier coefficients. We first prove a preliminary formula, expressing the Fourier coefficients in terms of gamma functions, Bessel functions and certain Dirichlet series. Next we show that these Dirichlet series have an Euler product, reducing the problem to the calculation of local factors for each prime. For almost all primes these local factors turn out to be {local factors of} ratios of classical zeta and Dirichlet $L$-functions. Finally, we calculate the local factors for the finitely many ramified primes, and use these calculations to prove the functional equation as well as show that some potential poles in fact get canceled by zeroes of these local factors.

\subsection*{Acknowledgements}
We would like to thank Tim Browning, Andre Reznikov, Peter Sarnak and Matthew Stover for their valuable comments regarding various aspects of this work.
We also thank the anonymous referee for his very thoughtful comments. 
\section{Preliminaries}\label{sec:prelim}
\subsection{Notation and conventions}
For any $d\in\N$ the notation $\prod_{p\mid d}$ denotes the product over all the prime divisors of $d$ while we use $\prod_{p\mid d}'$ to denote the product over all the odd prime divisors of $d$. For any lattice $\Lambda$ in $\R^n$ we denote by $\Lambda^*$ its dual lattice {(taken with respect to the standard inner product on $\R^n$)} and $\Lambda_{\rm pr}$ the set of primitive points of $\Lambda$. We use $\|\cdot\|$ exclusively for the Euclidean norm in various Euclidean spaces. We denote by $\R^+$ the set of positive real numbers. For two positive quantities $A$ and $B$, we will use the notation $A\ll B$ {or $A=O(B)$} to mean that there is a constant $c>0$ such that $A\leq cB$, and we will use subscripts to indicate the dependence of the constant on parameters. We will write $A\asymp B$ for $A\ll B\ll A$. \\

Let $n\in\N$ be a positive integer. For $Q$ an isotropic integral quadratic form of signature $(n+1,1)$ let $G=\SO^+_Q(\R)$ denote the identity component of the special orthogonal group preserving this form.
Let $\cV_Q^+$ denote a one-sheeted light cone of $Q$. Explicitly, when $Q=Q_n$ as in \eqref{e:Qn}, we can take
\begin{align*}
\cV^+_{Q}=\left\{ \bm{v}\in\R^{n+2} : Q(\bm{v})=0,\ v_{n+2}>0\right\}
\end{align*}
the positive light cone of $Q$. In general, for $Q(\bm{v})=Q_n(\bm{v}\tau)$ for some $\tau\in \GL_{n+2}(\R)$, we can take $\cV_Q^+=\cV_{Q_n}^+\tau^{-1}$. Fix a base point $\bm{e}_0\in \cV^+_Q$, the stabilizer $P\leq G$ of the line spanned by $\bm{e}_0$ is a parabolic subgroup, and we have that $G=PK$ with $K$ a maximal compact subgroup.
Moreover, the group $P$ has a Langlands decomposition $P=UAM$ with $U\cong\R^n$ the unipotent radical of $P$, $A$ diagonalizable one-parameter group and $M\cong \SO_n(\R)$ the centralizer of $A$ in the maximal compact group $K$.
We note that the groups $M$ and $U$ fix $\bm{e}_0$ while the group $A$ acts on $\bm{e}_0$ by scaling. We denote by $L:=UM$ the group stabilizing $\bm{e}_0$.

\subsection{Coordinates and measures}\label{sec:cormea1}
Let $G=\SO^+_Q(\R)$ and its subgroups as above. It has an Iwasawa decomposition $G=UAK$ which gives natural coordinates on $G$. Explicitly, fixing an isomorphism of $U$ with $\R^n$, for any $\bm{x}\in\R^n$ we denote by $u_{\bm{x}}\in U$ the corresponding element.
 Next, for $y>0$ let $a_y\in A$ satisfy $\bm{e}_0 a_y=y^{-1}\bm{e}_0$. 
 Then any $g\in G$ can be written as $g=u_{\bm{x}}a_yk$ with $k\in K$ and in these coordinates the Haar measure of $G$ is given by (up to scalars)
\begin{align}\label{equ:Haar1}
\d\mu_G(g)=y^{-(n+1)}\,\d\bm{x}\d y\d\mu_K(k),
\end{align}
where $\d\bm{x}$ is the usual Lebesgue measure on $\R^n$ and $\mu_K$ is the probability Haar measure of $K$.

Note that the subgroup $L=UM$ is unimodular and its Haar measure (up to scalars) is given by 
\begin{align}\label{equ:Lhaar}
\d\mu_L(u_{\bm{x}}m)=\d\bm{x}\d\mu_M(m),
\end{align}
where $\mu_M$ is the probability Haar measure of $M$. Let $\cV^+_{Q}$
be the one-sheeted light cone as above. The map from $G$ to $\cV^+_Q$ sending $g$ to $\bm{e}_0g$ then induces an identification between $\cV^+_Q$ and the homogeneous space $L\bk G$. This is true since $\Stab_{G}(\bm{e}_{0})=L$ and the right multiplication action of $G$ on $\cV^+_Q$ is transitive (see e.g. \cite{HoweTan1993}). In particular, since both $L$ and $G$ are unimodular groups, there exists a unique (up to scalars) right $G$-invariant measure on $\cV^+_Q$, which we denote by $m_{\cV^+_Q}$.
Further identifying $L\bk G$ with $A\times M\bk K$ gives 
 natural polar coordinates on $\cV^+_Q$: Every $\bm{v}\in \cV^+_Q$ can be written uniquely as $\bm{v}=\bm{e}_0a_yk$ for some $y>0$ and $k\in M\bk K$. In these coordinates $m_{\cV^+_Q}$ is given by
\begin{align}\label{equ:lebmea}
\d m_{\cV^+_Q}(\bm{e}_0a_yk)=y^{-(n+1)}\,\d y\d\mu_{M\bk K}(k),
\end{align}
where $\mu_{M\bk K}$ is the unique right $K$-invariant probability measure on the homogeneous space $M\bk K$ which is homeomorphic to the unit sphere $S^n$.

\subsubsection{Explicit coordinates}\label{sec:explicit}
For some calculations it is convenient to work with explicit coordinates. 
For any $n, d\in \N$ let $Q=Q_{n,d}$ denote the quadratic form in $n+2$ variables described in \eqref{e:Qnd}.
For this form we can take our base point to be  $\bm{e}_{0}:=(-d,0, \cdots,0,1)\in\R^{n+2}$. With these choices, we have that the compact groups $K$ and $M$ are given by
\begin{equation}\label{eq:KM}
K=\left\{\begin{pmatrix} k & 0\\ 0 &1\end{pmatrix}: k\in \SO_{n+1}(\R)\right\}\,\mbox{ and } M=\left\{\begin{pmatrix} 1 & 0 &0 \\ 0 &k &0\\ 0&0&1\end{pmatrix}: k\in \SO_n(\R)\right\},\end{equation}
and the parametrization for $u_{\bm{x}}\in U$ and $a_y\in A$ are given by 
\begin{align}\label{equ:uxay}
u_{\bm{x}}=\left(\begin{smallmatrix} 1-\frac{d^2\|\bm{x}\|^2}{2} & d\bm{x} & \frac{d\|\bm{x}\|^2}{2}\\ 
-d\bm{x}^t & I_n & \bm{x}^t\\ 
-\frac{d^3\|\bm{x}\|^2}{2} & d^2\bm{x} & 1+\frac{d^2\|\bm{x}\|^2}{2}\end{smallmatrix}\right)\;
\mbox{ and }\;
a_y=\left(\begin{smallmatrix}
\tfrac{y+y^{-1}}{2} &  & \tfrac{y-y^{-1}}{2d}\\
 & I_n & \\
\tfrac{d(y-y^{-1})}{2} &  &\tfrac{y+y^{-1}}{2} \end{smallmatrix}\right).
\end{align}

\begin{rem}\label{rem:Coordinates}
With these explicit coordinates the identification of $G/K$ with the upper half space $\bH^{n+1}$ can be made explicit, as the element $gK$ with $g=u_{\bm{x}}a_yk$ is mapped to $z=(\bm{x},y)\in \bH^{n+1}$.
\end{rem}

\subsection{Cusps}\label{sec:cusps}
Let $G=\SO^+_Q(\R)$ and let $P=UAM$ be the parabolic subgroup as above. Since $G$ is of real rank one, every parabolic subgroup of $G$ is conjugate to $P$. Moreover, since $P$ is self-normalizing, the space of parabolic subgroups in $G$ can be parameterized by the homogeneous space $P\bk G\cong M\bk K$. In particular, every parabolic subgroup $P'$ of $G$ is of the form $P'=kPk^{-1}$ for some $k\in M\bk K$ and it has a Langlands decomposition $P'=U'A'M'$ with $U'=kUk^{-1}$, $A'=kAk^{-1}$ and $M'=kMk^{-1}$.

Let $\G< G$ be a \textit{non-uniform lattice} of $G$, that is, $\G$ is a discrete subgroup and the homogeneous space $\G\bk G$ is non-compact and of finite volume.
The \textit{cusps} of $\Gamma$ are the $\Gamma$-conjugacy classes of parabolic subgroups of $G$ whose unipotent radicals intersecting $\Gamma$ nontrivially. 
Let $P_1,\ldots, P_\kappa$ denote a full set of representatives for these classes and let $k_i\in M\bk K$ be such that $P_i=k_iPk_i^{-1}$. For each cusp $P_i=U_iA_iM_i$, let $\Gamma_{P_i}:=\Gamma\cap P_i$ and $\Gamma_{U_i}:=\Gamma\cap U_i$. By definition, $\G_{U_i}$ is nontrivial. Indeed, since $\G\bk G$ has finite volume, $\G_{U_i}< U_i$ ($\cong \R^n$) is a lattice (so that $\G_{U_i}\bk U_i$ is of finite volume). 
We first prove a simple lemma regarding these two groups $\Gamma_{P_i}$ and $\G_{U_i}$.
\begin{Lem}\label{lem:finindtg}
Let $\G$ be a non-uniform lattice of $G$, and suppose $\G$ has a cusp at $P'$ for some parabolic subgroup $P'=U'M'A'$. Then $\G_{P'}\leq U'M'$ and $[\G_{P'} : \G_{U'}]<\infty$.
\end{Lem}
\begin{proof}
Up to conjugating $\G$ by elements in $G$, we may assume $P'=P$. We first show $\G_{P}< L=UM$. Suppose not, then there exists $\gamma\in \G_P$ such that $\gamma=uam$ for some $u\in U$, $m\in M$ and nontrivial $a\in A$. Take a nontrivial element $u'\in \G_{U}$. Then up to replacing $\gamma$ by $\gamma^{-1}$ if necessary, we have $\gamma^nu'\gamma^{-n}\to \textrm{id}$ as $n\to\infty$. On the other hand, the sequence $\{\gamma^nu'\gamma^{-n}\}_{n\in\N}$ is a sequence of nontrivial elements in $\G$. This contradicts the fact that $\G$ is a discrete subgroup of $G$. 

For the second assertion, consider the map $\mathfrak{c}: L\to \SO(U)$ sending $h\in L$ to the conjugation map $\mathfrak{c}(h)(u):=huh^{-1}$. Explicitly, one easily verifies that if $h=u'm$ for some $u'\in U$ and $m=\diag(1,m',1)\in M$ with $m'\in \SO_n(\R)$, then $\mathfrak{c}(h)(u_{\bm{x}})=u_{\bm{x}m'{^{-1}}}$. Thus this map induces an identification between $U\bk L$ and $\SO(U)$. In particular, since $\G_P\subseteq L$ and $\G_U=\G\cap U=\G_P\cap U$, we have an embedding from $\G_U\bk \G_P$ to $\SO(U)$. On the other hand, since the conjugation map of elements in $\G_P$ preserves $\G_U$, we have $\mathfrak{c}(\G_P)\subseteq \GL(\G_U)$, which again induces an embedding from $\G_U\bk \G_P$ to $\GL(\G_U)$. But $\GL(\G_U)\cap \SO(U)$ is finite (since $\G_U< U$ is a lattice), we also have $\G_U\bk \G_P$ is finite. This finishes the proof.
\end{proof}

Next we give another interpretation of the cusps in {terms} of the action of $\G$ on  $\cV_Q^+$. 
Recall that the homogeneous space $P\bk G$ parameterizes the space of parabolic subgroups of $G$. On the other hand, $P\bk G=LA\bk G\cong A\bk \cV^+_Q$, where the latter parameterizes the space of \textit{rays} on $\cV^+_Q$. Here by a ray on $\cV^+_Q$ we mean a ray which starts from the origin and lies on $\cV^+_Q$. For any $\bm{v}\in \cV^+_Q$ we denote by $[\bm{v}]\in A\bk \cV^+_Q$ the
unique ray on $\cV^+_Q$ which passes through $\bm{v}$. The above identification $A\bk \cV^+_Q\cong P\bk G$ then reads as $[\bm{v}]\in A\bk \cV^+_Q\leftrightarrow P_{[\bm{v}]}\in P\bk G$, where $P_{[\bm{v}]}\in P\bk G$ is the unique parabolic subgroup fixing $[\bm{v}]$. More explicitly, let $P_{[\bm{v}]}=U_{[\bm{v}]}M_{[\bm{v}]}A_{[\bm{v}]}$ be the Langlands decomposition for $P_{[\bm{v}]}$, then $U_{[\bm{v}]}M_{[\bm{v}]}$ fixes $\bm{v}$ and $A_{[\bm{v}]}$ acts on $\bm{v}$ via positive scaling.

Using this identification, we  say \textit{$\G$ has a cusp at $[\bm{v}]$} if $\G$ has a cusp at $P_{[\bm{v}]}$. We say two rays $[\bm{v}_1], [\bm{v}_2]$ are \textit{$\G$-equivalent} if $P_{[\bm{v}_1]}$ and $P_{[\bm{v}_2]}$ are $\G$-conjugate, otherwise we say $[\bm{v}_1]$ and $[\bm{v}_2]$ are \textit{$\G$-inequivalent}. One easily verifies that $[\bm{v}_1]$ and $[\bm{v}_2]$ are $\G$-equivalent if and only if $[\bm{v}_1\gamma]=[\bm{v}_2]$ for some $\gamma\in \G$.
Finally, we say a ray $[\bm{v}]\in A\bk \cV^+_Q$ is \textit{rational} if $[\bm{v}]$ contains a rational vector.
For the special case of $\G=\SO^+_Q(\Z)$ (note that $\SO^+_Q(\Z)$ is a non-uniform lattice since we have assumed $Q$ is isotropic) we have that $\G$ has a cusp at a ray $[\bm{v}]$ if and only if $[\bm{v}]$ is rational (see \cite[Proposition 3.15 (ii)]{FishmanKleinbockMerrillSimmons2022}). Hence the  
 cusps of $\G$ correspond to the rational rays on $\cV_Q^+$, or equivalently, after clearing denominators, to $\G$-orbits of $\cV^+_Q(\Z)_{\rm pr}$, the set of primitive integer points on the light cone $\cV_Q^+$.
We summarize this correspondence in the following.	
\begin{Lem}\label{lem:cusporbits}
The set of primitive integral points on the light cone can be decomposed as 
\begin{align}\label{equ:orbidecp}
		\cV^+_Q(\Z)_{\rm pr}=\bigsqcup_{i=1}^\kappa\bm{v}_i\Gamma
	\end{align}
with $\bm{v}_1,\ldots, \bm{v}_{\kappa}\in \cV_Q^+(\Z)_{\rm pr}$ so that $\{[\bm{v}_1],\ldots, [\bm{v}_\kappa]\}$ is a full set of $\G$-inequivalent cusps.
\end{Lem}
\begin{proof}
Let $[\bm{v}_1],\ldots, [\bm{v}_\kappa]$ be a full set of $\G$-inequivalent cusps. By \cite[Proposition 3.15 (ii)]{FishmanKleinbockMerrillSimmons2022}) they are all rational rays, thus we may take a representative $\bm{v}_i\in \cV^+_Q(\Z)_{\rm pr}$. For each $1\leq i\leq \kappa$, let $\cO_i:=\bm{v}_i\G$. Since $[\bm{v}_1],\ldots, [\bm{v}_\kappa]$ are pairwise $\G$-inequivalent, the orbits $\cO_i$ are pairwise disjoint. 
	On the other hand, for any $\bm{v}\in \cV^+_Q(\Z)_{\rm pr}$, $[\bm{v}]$ is a rational ray so $\G$ has a cusp at $[\bm{v}]$, implying that $[\bm{v}]$ is $\G$-equivalent to some $[\bm{v}_i]$. That is, there exists some $\gamma\in \G$ such that $[\bm{v}]=[\bm{v}_i\gamma]$, or equivalently, $\bm{v}=\lambda \bm{v}_i\gamma$ for some $\lambda>0$. But both $\bm{v}$ and $\bm{v}_i\gamma$ are elements from $\cV_Q^+(\Z)_{\rm pr}$, we must have $\lambda=1$, implying that $\bm{v}=\bm{v}_i\gamma\in \cO_i$. 
\end{proof}

This correspondence between the cusps and the $\G$-orbits of primitive vectors gives a natural notion of a volume of a cusp.		
\begin{Lem}\label{lem:cuspvol}
	For each $1\leq i\leq \kappa$, let $\bm{v}_i\in \cO_i$ and let $P_i$ be the parabolic subgroup fixing the ray $[\bm{v}_i]$. Let $g_i\in G$ such that $\bm{v}_i=\bm{e}_0g_i^{-1}$. Then the cusp volume
	\begin{align}\label{equ:lamcome}
	v_{P_i}:=\mu_L\left(g_i^{-1}\G_{P_i}g_i\bk L\right)
	\end{align}
is well defined and independent of the choices of representative $\bm{v}_i\in \cO_i$ and  $g_i\in G$.
\end{Lem}	
\begin{proof}
It suffices to show that the quantity $\mu_L\left(g_i^{-1}\G_{P_i}g_i\bk L\right)$ is independent of the choices of $\bm{v}_i\in \cO_i$ and $g_i\in G$. Suppose $\bm{w}\in \cO_i=\bm{v}_i\G$ is another element in $\cO_i$, i.e. $\bm{w}=\bm{v}_i\gamma$ for some $\gamma\in \G$. Let $P'$ be the parabolic subgroup fixing $[\bm{w}]$ and let $g_0\in G$ be such that $\bm{w}=\bm{e}_0g_0^{-1}$. We thus would like to show 
	$$
	\mu_L\left(g_i^{-1}\G_{P_i}g_i\bk L\right)=\mu_L\left(g_0^{-1}\G_{P'}g_0\bk L\right).
	$$
	First note that since $\bm{w}=\bm{v}_i\gamma$ and $\bm{w}=\bm{e}_0g_0^{-1}$, we have $P'=\gamma^{-1}P_i\gamma$ and $\bm{v}_i=\bm{e}_0g_0^{-1}\gamma^{-1}$. Combining the latter relation with the relation $\bm{v}_i=\bm{e}_0g_i^{-1}$ we get $\bm{e}_0g_0^{-1}\gamma^{-1}g_i=\bm{e}_0$; thus $g_0^{-1}\gamma^{-1}g_i\in L$, or equivalently, $g_i=\gamma g_0 h$ for some $h\in L$. Then
	$$
	g_i^{-1}\G_{P_i}g_i=h^{-1}g_0^{-1}\gamma^{-1}(\G\cap P_i)\gamma g_0 h=h^{-1}g_0^{-1}(\G\cap P') g_0h=h^{-1}g_0^{-1}\G_{P'} g_0h.
	$$
	Thus 
	$$
	\mu_L\left(g_i^{-1}\G_{P_i}g_i\bk L\right)=\int_{h^{-1}g_0^{-1}\G_{P'}g_0h\bk L}\,\d\mu_L.
	$$
	Now let $\cF\subseteq L$ be a fundamental domain for $g_0^{-1}\G_{P'}g_0\bk L$, then $h^{-1}\cF$ is a fundamental domain for $h^{-1}g_0^{-1}\G_{P'} g_0h\bk L$. Thus by invariance of the Haar measure we get 
	\begin{align*}
		\mu_L\left(g_i^{-1}\G_{P_i}g_i\bk L\right)=\mu_L(h^{-1}\cF)=\mu_L(\cF)=\mu_L\left(g_0^{-1}\G_{P'}g_0\bk L\right),
	\end{align*}
	proving the desired identity. 
\end{proof}

\begin{remark}\label{rmk:cuspvol1}
	{Let $\G=\SO_Q^+(\Z)$ be as above.} Even for the explicit family of forms $Q_{n,d}$ defined in \eqref{e:Qnd}, when there is more than one cusp it is not clear to us how to compute these cusp volumes. Except for the cusp corresponding to the orbit $\bm{e}_0\G$ (or the parabolic group 		$P$)  for which we have  $v_P=2^{2-n}(n!)^{-1}$ if $d$ is odd and $v_P=2^{1-n}(n!)^{-1}$ if $d$ is even.
	To see this,  take $\bm{v}_1=\bm{e}_0$ so  $v_P=\mu_L(\G_P\bk L)$. By definition of $\G$ and the fact that $\G_P< L$ (cf. \lemref{lem:finindtg}), we have $\G_P=L(\Z):=L\cap \SL_{n+2}(\Z)$. It is then not difficult to see that $L(\Z)=U(\Z)M(\Z)$, 
	where $U(\Z):=U\cap \SL_{n+2}(\Z)$ 
	and $M(\Z):=M\cap \SL_{n+2}(\Z)\cong \SO_n(\Z)$ and the latter has order $\#\SO_n(\Z)=2^{n-1}n!$.
	Finally  
	we have that $U(\Z)=\left\{u_{\bm{x}} : \bm{x}\in \Lambda \right\}$, where $\Lambda=\Z^n$ if $d$ is even, while for $d$ odd,
	\begin{align}\label{equ:lambdalatt}
		\Lambda:=\left\{\bm{x}\in \Z^n: \|\bm{x}\|^2\equiv 0\Mod{2}\right\},
	\end{align}
	 is an index $2$ sub-lattice of $\Z^n$. 
\end{remark}

\subsection{Eisenstein series}\label{sec:Eisenstein}
In this subsection we review necessary backgrounds on the spectral theory of Eisenstein series on real hyperbolic manifolds. We refer the reader to \cite[Section 4]{Sodergren2012} and the references therein for more details.

For each $1\leq i\leq \kappa$, let $P_i=k_iPk_i^{-1}$ be as above. Since $A$ normalizes $L$, we have $L_i=k_ia_yLa_y^{-1}k_i^{-1}$ for any $a_y\in A$. In particular, since by \lemref{lem:finindtg}, $\G_{P_i}< L_i$, we have $a_y^{-1}k_i^{-1}\G_{P_i}k_ia_y< L$ for any $a_y\in A$. We fix the \textit{scaling matrix} $\tau_i:=k_ia_{y_i}$ with $y_i>0$ the unique positive number such that $\mu_L(\tau_i^{-1}\G_{P_i}\tau_i\bk L)=1$.
The \textit{(spherical) Eisenstein series} corresponding to the $i$-th cusp is then defined 
for $\Re(s)>n$ by the convergent series
\begin{equation}\label{e:EisensteinSph}
	E_i(s, g):=\sum_{\gamma\in \Gamma_{\!P_i}\bk \Gamma} y(\tau_i^{-1}\gamma g)^s,
\end{equation}
where $y(g)$ is defined by the Iwasawa decomposition $g=u_{\bm{x}}a_{y(g)}k$.
The constant term of $E_i(s,g)$ with respect to the $j$-th cusp is defined by
\begin{equation}\label{e:ConstSph} 
	c_{ij}(s,g)=\frac{1}{\vol(\tau_j^{-1}\G_{U_j}\tau_j\bk U)}\int_{\tau_j^{-1}\G_{U_j}\tau_j\bk U}E_i(s, \tau_ju_{\bm x} g)\,\d\bm{x},
\end{equation}
where $\vol(\cU):=\int_{\cU}\,\d\bm{x}$ for any Borel subset $\cU\subseteq U (\cong \R^n)$. These constant terms are of the form
\begin{equation}\label{e:ConstSph1}
c_{ij}(s,g)=\delta_{ij}y(g)^s+\varphi_{ij}(s)y(g)^{n-s},
\end{equation}
where $\varphi_{ij}(\cdot)$ ($1\leq i,j\leq \kappa$) are some holomorphic functions defined on the half plane $\Re(s)>n$.

The matrix $\Psi(s):=(\varphi_{ij}(s))_{1\leq i,j\leq \kappa}$ is called the \textit{scattering matrix}. 
Let $\mathcal{E}(s,g)$ be the column vector with the $i$-th coordinate given by $E_i(s,g)$. Then $\mathcal{E}(s,g)$ (and hence also $\Psi(s)$) has a meromorphic continuation to the whole
$s$-plane, and satisfies the functional equation
\begin{align}\label{equ:funcequ1}
\mathcal{E}(s,g)=\Psi(s)\mathcal{E}(n-s, g).
\end{align}
As a consequence, the scattering matrix $\Psi(s)$ satisfies the functional equation $\Psi(s)\Psi(n-s)=I_\kappa$, and due to our choice of scaling matrices, it is also symmetric and satisfies that {$\overline{\Psi(\bar{s})}=\Psi(s)$}. In particular, $\Psi(s)$ is unitary on the critical line $\Re(s)=\frac{n}{2}$, implying that $\Psi(s)$ is holomorphic on $\Re(s)=\frac{n}{2}$ and
\begin{align}\label{equ:funcons1}
	|\varphi_{ij}(s)|\leq 1, \qquad\forall\ \Re(s)=\tfrac{n}{2},\ 1\leq i,j\leq \kappa.
	\end{align}
More generally, for each $1\leq i\leq \kappa$, the (meromorphically continued) Eisenstein series $E_i(s,g)$ (and the coefficient function $\varphi_{ii}(s)$) is holomorphic on the half plane $\Re(s)\geq \frac{n}{2}$ except for a simple pole at $s=n$ (called the \textit{trivial pole}) and possibly finitely many simple poles on the interval $(\frac{n}{2}, n)$ (called \textit{exceptional poles}). 

\begin{rem}\label{r:WboundE}
For a general lattice $\Gamma$ not much is known regarding the growth {of} the Eisenstein series $E_i(\tfrac{n}{2}+it,g)$ in terms of $t$.
We record here that the best known general bound (see \cite[Proposition 7.13 and {Theorem} 7.14]{CohenSarnak1980}) states that there is a nonnegative function $W(t)\geq 1$ with $\int_{0}^T W(t)\,\d t\ll T^{n+1}$
such that $|E_i(\sigma+it,g)|^2\ll_g W(t)t^{n+1}$ for all $\sigma\in [\tfrac{n}{2},n]$ and $|t|\geq 1$.
On the critical line, $\sigma=\frac{n}{2}$ one can do better and show that $\int_{-T}^T \left|E_i(\frac{n}{2}+it,g)\right|^2\,\d t\ll {y(g)}^nT+T^{n+1}$ \cite[Corollary 7.7]{CohenSarnak1980}, giving an averaged version of the convexity bound.

\end{rem}

We note that the residue of $E_i(s,g)$ at $s=n$ is an eigenfunction of the Laplacian with eigenvalue zero and hence a constant.  
This constant is the same for Eisenstein series at all cusps and is the reciprocal of the measure of the quotient space $\G\bk G$.
While this is well known, and was shown for example in \cite[Lemma 2.15]{Sarnak1983} for Eisenstein series on hyperbolic $3$-manifolds,  since we could not find a reference proving this in general we include a short proof below.
\begin{Prop}
For each $1\leq i\leq \kappa$ we have that 
$\Res_{s=n}E_i(s,g)=\frac{1}{\mu_G(\G\bk G)}$.
\end{Prop}
\begin{proof}
Let $\rho$ be a smooth  compactly supported function on $\R^+$ and define $f$ on $\cV_Q^+$ by $f(\bm{e}_0a_yk)=\rho(y)$.
The  \textit{incomplete Eisenstein series attached to $f$ at the cusp $P_i$} is the function on $\G\bk G$ defined by the sum 
\begin{align}\label{e:Egf}
E_i(g;f):=\sum_{\gamma\in \G_{P_i}\bk \G}f(\bm{e}_0\tau_i^{-1}\gamma g),\qquad (g\in G),
\end{align}
which is a finite sum since $f$ is compactly supported.
Using the standard unfolding trick and the relation $\mu_L(\tau_i^{-1}\G_{P_i}\tau_i\bk L)=1$ we see that 
\begin{align}\label{equ:firstmoment}
		\int_{\G\bk G}E_i(g;f)\,\d\mu_G(g)=\int_0^\infty \rho(y)\frac{\d y}{y^{n+1}}=\hat\rho(n),
		\end{align}
		where 
\begin{align*}
	\hat{\rho}(s):=\int_0^{\infty}\rho(y)y^{-(s+1)}\,\d y,\qquad (s\in\C),
\end{align*}
is the Mellin transform of $\rho$.
On the other hand using Mellin inversion formula
\begin{align}\label{equ:mellininv}
	\rho(y)=\frac{1}{2\pi i}\int_{(\sigma)}\hat{\rho}(s)y^s\,\d s,\qquad (\sigma\in\R),
\end{align}
we can write for any $\sigma>n$ 
\begin{align*}
	E_i(g;f)=\frac{1}{2\pi i}\int_{(\sigma)}\hat{\rho}(s)E_i(s, g)\,\d s.
\end{align*}
The fast decay of $\hat\rho(s)$ together with the bounds on the Eisenstein series given in \rmkref{r:WboundE} allow us to shift the contour to the line $\Re(s)=\tfrac{n}{2}$.
Collecting the contribution of the {trivial} pole with residue $\Res_{s=n}E_i({s},g)$ and perhaps finitely many exceptional poles with residues 
$\vf_j(g):=\Res_{s=s_j}E_i(s,g)$ we have 
\begin{align*}
	E_i(g;f)=\Res_{s=n}E_i({s},g)\hat\rho(n)+\sum_j \hat\rho(s_j) \vf_j(g)+\frac{1}{2\pi i}\int_{(\tfrac{n}{2})}\hat{\rho}(s)E_i(s, g)\,\d s.
\end{align*}
{Using that the Fourier expansion of the Eisenstein series  $E_i(s,\tau_j g)$ at the $j$-th cusp is given by $c_{i,j}(s,g)$ plus terms that decay exponentially at cusp (see \cite[Section 4.1]{Sodergren2012}), after taking residues it is not hard to see that the residual forms $\vf_j(g)$ are square integrable, and since they are Laplacian eigenfunctions with non zero eigenvalues they are orthogonal to constant functions. Next using the above identity and the fact that $E_i(g;f)$ is square integrable, we see that $\frac{1}{2\pi i}\int_{(\tfrac{n}{2})}\hat{\rho}(s)E_i(s, g)\,\d s$ is also integrable, and since it is given by an integral over Laplacian eigenfunctions with eigenvalues in $[\tfrac{n^2}{4},\infty)$ this term is also orthogonal to constant functions.} 
Thus, integrating this identity over the fundamental domain for $\G\bk G$, 
we get that 
$$\int_{\G\bk G}E_i(g;f)\,\d\mu_G(g)= \Res_{s=n}E_i(s,g)\hat\rho(n)\mu_G(\G\bk G).$$
Comparing this with \eqref{equ:firstmoment} we see that $\Res_{s=n}E_i(s,g)=\tfrac{1}{\mu_G(\G\bk G)}$.
\end{proof}

\subsection{Kronecker symbol and $L$-functions}\label{sec:krsylfun}
We end this section by reviewing necessary backgrounds on Kronecker symbol and its relations with real Dirichlet character and the functional equation satisfied by the $L$-functions attached to these characters. The main reference is \cite[Chapters 9 and 10]{MontgomeryVaughan2007}.

Recall that for any nonzero integer $D\not\equiv 3\Mod{4}$ we denote by $\chi_{D}$ the unique primitive Dirichlet character that induces the real Dirichlet character $k\mapsto \left(\tfrac{D}{k}\right)$, {where} $\left(\tfrac{\cdot}{\cdot}\right)$ is the Kronecker symbol; see \eqref{equ:prchardes} below for a more precise description of $\chi_D$. We note that all the real Dirichlet characters can be realized this way. When $p$ is an odd prime $\left(\tfrac{\cdot}{p}\right)$ agrees with the more classical Legendre symbol, $\left(\tfrac{a}{2}\right)=(-1)^{\frac{a^2-1}{8}}$ for any odd $a$ and $\left(\tfrac{a}{-1}\right)=\left\lbrace\begin{array}{ll}
		-1& a<0,\\
		1 & a\geq 0.
	\end{array}\right.$ 

For $\Re(s)>1$ let $L(s,\chi_D):=\sum_{k=1}^{\infty}\frac{\chi_D(k)}{k^s}$ the $L$-function attached to $\chi_D$. Let $q$ be the modulus of $\chi_D$ and define
\begin{align}\label{equ:comlfun}
L^*(s,\chi_D):=\left(\tfrac{q}{\pi}\right)^{\frac{s+a}{2}}\G(\tfrac{s+a}{2})L(s,\chi_D)
\end{align}
the complected $L$-function. Here $a=0$ if $\chi_D(-1)=1$ (i.e. $D>0$) and $a=1$ if $\chi_D(-1)=-1$ (i.e. $D<0$). Since $\chi_D$ is real (so that $\overline{\chi}_D=\chi_D$), $L^*(s,\chi_D)$ satisfies the following functional equation (see e.g. \cite[p. 333]{MontgomeryVaughan2007})
\begin{align*}
L^*(s,\chi_D)=\tau(\chi_D)L^*(1-s,\chi_D),
\end{align*}
where $\tau(\chi_D):=\frac{\cG(\chi_D)}{i^a\sqrt{q}}$ with $\cG(\chi_D):=\sum_{j=1}^q\chi_D(j)e\left(\tfrac{j}{q}\right)$ the Gauss sum attached to $\chi_D$. Since $\chi_D$ is primitive, we have by \cite[Theorem 9.17]{MontgomeryVaughan2007} that $\cG(\chi_D)=i^a\sqrt{q}$, or equivalently, $\tau(\chi_D)=1$. Thus $L^*(s,\chi_D)$ satisfies the following simple functional equation 
\begin{align}\label{equ:funequlfun}
L^*(s,\chi_D)=L^*(1-s,\chi_D).
\end{align}

\section{The light-cone Eisenstein series}\label{sec:LCES}
Fix a rational isotropic quadratic form, $Q$, of signature $(n+1,1)$ and let $G=\SO_Q^+(\R)$ and $\G=\SO_Q^+(\Z)$. Let $\cV_Q^+$ be a fixed one-sheeted light cone of $Q$ and $\bm{e}_0\in \cV_Q^+$ a fixed base point as before. Let $P=UAM\leq G$ and $K\leq G$ be the corresponding subgroups as before. Let $\|\cdot\|_K$ be a $K$-invariant norm on $\R^{n+2}$. For instance, for $Q(\bm{v})=Q_n(\bm{v}\tau)$ with $Q_n$ as in \eqref{e:Qn} and $\tau\in \GL_{n+2}(\R)$, we may take $\|\bm{v}\|_K:=\|\bm{v}\tau\|$. 
When $Q=Q_{n,d}$, in view of the choice of the subgroup $K$ as described in \eqref{eq:KM} we will take $\|\cdot\|_K$ to be the Euclidean norm.
The light-cone Eisenstein series attached to $Q$ is a function {on} $\G\bk G$ defined for $\Re(s)>n$ by the series 
\begin{equation}
E_Q(s,g):=\|\bm{e}_0\|_K^s \sum_{\bm{v}\in \cV^+_Q(\Z)_{\rm{pr}}}\|\bm{v}g\|_K^{-s},
\end{equation}
where $\cV_Q^+(\Z)_{\rm pr}$ is the set of primitive integer points on $\cV_Q^+$ as before.
In this section we show how the light-cone Eisenstein series is related to the Eisenstein series attached to the cusps of $\G$ and then outline some applications of the analytic properties of $E_Q(s,g)$ to some counting problems.

\subsection{Relation to Eisenstein series}
 The light-cone Eisenstein series can be expressed as a weighted sum of Eisenstein series.
More generally we show this for incomplete Eisenstein series.
\begin{Lem}\label{lem:siegeltrantoincomeisen}
For any function $f$ on $\cV^+_Q$ {satisfying that  $|f(\bm{v})|\ll \|\bm{v}\|^{-\nu}$ with $\nu>n$ for all $\bm{v}\in \cV^+_Q$ with $\|\bm{v}\|\geq 1$}, let 
$$E_Q(g;f):=\sum_{\bm{v}\in \cV^+_Q(\Z)_{\rm pr}}f(\bm{v}g).$$
Let $\kappa$ denote the number of cusps of $\G$.  Then 
	\begin{align}\label{equ:sitrantoincei}
		E_Q(g;f)=\sum_{i=1}^\kappa E_i(g;f_{\lambda_i}),
	\end{align}
	where $f_{\lambda}(\bm{v}):=f(\lambda^{-1}\bm{v})$ is the scaled function, with the scaling factors  $\lambda_i>0$ satisfying that $\lambda_i^n=v_{P_i}$ with $v_{P_i}$ the cusp volume of the $i$-th cusp defined as in \eqref{equ:lamcome}, {and $E_i(g; f_{\lambda_i})$ are the incomplete Eisenstein series defined in \eqref{e:Egf}}.
	\end{Lem}
\begin{proof}
Let $\bm{v}_1,\ldots, \bm{v}_{\kappa}\in \cV_Q^+(\Z)_{\rm pr}$ be such that the orbit decomposition \eqref{equ:orbidecp} in \lemref{lem:cusporbits} holds.
	Then $[\bm{v}_1],\ldots, [\bm{v}_{\kappa}]$ is a full set of representatives for the cusps. Let $P_1=k_1Pk_1^{-1}, \ldots, P_{\kappa}=k_{\kappa}Pk_{\kappa}^{-1}$ (with $k_1,\ldots, k_{\kappa}\in K$) be the corresponding parabolic subgroups. 
	For each $1\leq i\leq \kappa$, let $\tau_i=k_ia_{y_i}$ be the scaling matrix as before and let $\bm{w}_i=\bm{e}_0\tau_i^{-1}$. Note that $[\bm{v}_i]=[\bm{e}_0]k_i^{-1}$; hence
	\begin{align*}
	[\bm{v}_i]=[\bm{e}_0]k_i^{-1}=[\bm{e}_0]a_{y_i}^{-1}k_i^{-1}=[\bm{e}_0]\tau_i^{-1}=[\bm{w}_i].
	\end{align*}
	This implies that there exists some $\lambda_i>0$ such that $\bm{w}_i=\lambda_i\bm{v}_i$. Then \eqref{equ:sitrantoincei} follows by the orbit decomposition \eqref{equ:orbidecp}, and noting that the map from $\G_{P_i}\gamma\in \G_{P_i}\bk \G$ to $\bm{v}_i\gamma\in \bm{v}_i\G$ is an identification (this is so since $\G_{P_i}=\G\cap P_i=\G\cap L_i$ is the stabilizer of $\bm{v}_i$ in $\G$; see \lemref{lem:finindtg}) and the relation $\bm{v}_i=\lambda_i^{-1}\bm{w}_i=\lambda_i^{-1}\bm{e}_0\tau_{i}^{-1}$.	
	
{Finally, to show that $\lambda_i^n=v_{P_i}$, write $\bm{v}_i=\bm{e}_0g_i^{-1}$ where $g_i=\tau_i a_{\lambda_i}^{-1}$} and recall $\tau_i$ is defined so that $\mu_L(\tau_i^{-1}\G_{P_i}\tau_i\bk L)=1$. 
Hence using the conjugating relation $a_yu_{\bm{x}}a_{y}^{-1}=u_{y\bm{x}}$ for any $y>0$ and $\bm{x}\in \R^n$ we have
$$v_{P_i}=\mu_{L}(g_i^{-1}\G_{P_i}g_i\bk L)=\mu_L(a_{{\lambda_i}}\tau_i^{-1}\G_{P_i}\tau_i a_{{\lambda_i}}^{-1}\bk L)=\lambda_i^n\mu_{L}(\tau_i^{-1}\G_{P_i}\tau_i\bk L)=\lambda_i^n,$$
as claimed. 
\end{proof}
In particular, noting that for $g=u_{\bm{x}}a_yk$ 
$$\|\bm{e}_0g\|_K=\|\bm{e}_0a_y\|_K=y(g)^{-1}\|\bm{e}_0\|_K,$$ 
we see that  $y(g)=\|\bm{e}_0\|_K\|\bm{e}_0 g\|_K^{-1}$. Hence taking  $f(\bm v)=\|\bm v\|_K^{-s}$ {with $\Re(s)>n$ in \lemref{lem:siegeltrantoincomeisen} we get}
$$E_i(g; f_{\lambda_i})=\lambda_i^s\|\bm{e}_0\|_K^{-s} \sum_{\g\in \G
_{P_i}\bk \G} y(\tau_i^{-1}\g g)^{s}=\left(\tfrac{\lambda_i}{\|\bm{e}_0\|_K}\right)^sE_i(s,g).$$
We thus get the following.

\begin{Cor}\label{cor:sereisnserre}
For any $\Re(s)>n$ the normalized light-cone Eisenstein series
	satisfies
\begin{align}\label{equ:sereisnserre}
E_Q(s,g)=\sum_{i=1}^\kappa v_{P_i}^{\frac{s}{n}}E_i(s,g).
\end{align}
In particular, $E_Q(s,g)$ has a meromorphic continuation to the whole $s$-plane, and is holomorphic in the half space $\Re(s)\geq \tfrac{n}{2}$ except for a simple pole at $s=n$ with constant residue
\begin{align}\label{equ:omeq}
\omega_Q:=\Res_{s=n}E_Q(s,g)=\mu_G(\G\bk G)^{-1}\sum_{i=1}^{k}v_{P_i},
\end{align} 
and possibly finitely many simple poles on the interval $(\tfrac{n}{2}, n)$.{\qed}
\end{Cor}

\subsection{Application to counting}
For any $T\geq 1$ and $g\in G$ consider the counting function 
$$N(T;g):=\left\{\bm{v}\in \cV_Q^+(\Z)_{\rm pr}: \tfrac{\|\bm{v}g\|_K}{\|\bm{e}_0\|_K}\leq T\right\}.$$
We now outline some general applications, showing how one can use properties of the light-cone Eisenstein series to estimate this counting function.
\begin{rem}\label{rmk:corres}
For the special case where $Q=Q_{n,d}$ we can identify $\bm{v}=(\bm{p},q)\in  \cV_Q^+(\Z)_{\rm pr}$ with the reduced rational point $\frac{\bm{p}}{q}\in dS^n$ on the sphere of radius $d$.
With this identification, and taking $g=1_G$ (the identity in $G$) we see that 
$$N(T; 1_G)=\#\left\{\tfrac{\bm{p}}{q}\in dS^n: (\bm{p},q)\in \Z^{n+2}_{\rm pr},\ 1\leq q\leq T\right\}$$
 counts the number of (reduced) rational points on the (radius $d$) sphere with height bounded by $T$. In this case it is clear that for $T\in \N$ we have that $N(T;1_G)=N(T+\tfrac{1}{2};1_G)$.
\end{rem}

To relate this counting function to the light-cone Eisenstein series we first consider the following smoothed counting.
For any $h\in C^\infty_c(\R^+)$ consider the  
 smoothed counting function
\begin{equation}\label{e:smoothed}N_{h}(T;g):=\sum_{\bm{v}\in \cV^+_{Q}(\Z)_{\rm pr}} h\left(T^{-1} \tfrac{\|\bm{v}g\|_K}{ \|\bm{e}_0\|_K}\right).\end{equation}
The following identity relates this smoothed counting function to the light-cone Eisenstein series.
\begin{Lem}\label{lem:MellinCounting}
For any $h\in C^\infty_c(\R^+)$ and for any $\Re(s)>n$, the Mellin transform
$\hat{N}_h(s;g)=\int_0^\infty N_h(t;g)t^{-s}\frac{\d t}{t}$ is given by  $\hat{N}_h(s;g)=\hat{h}(-s)E_{Q}(s,g)$.
\end{Lem}
\begin{proof}
Let $\rho(y)=h(\tfrac{1}{Ty})$ so its Mellin transform is 
$\hat{\rho}(s)=T^s\hat{h}(-s)$.
We can thus write for any $\sigma>n$,
\begin{align*}
 N_h(T;g)
  &= \sum_{\bm{v}\in \cV_Q^+(\Z)_{\rm pr}}\rho (\|\bm{e}_0\|_K\|\bm{v}g\|_K^{-1}) \\
   &=\sum_{\bm{v}\in \cV_Q^+(\Z)_{\rm pr}}\frac{1}{2\pi i} \int_{(\sigma)} \hat\rho(s)  \|\bm{e}_0\|_K^s\|\bm{v}g\|_K^{-s}\, \d s \\
 &= \frac{1}{2\pi i} \int_{(\sigma)} \hat\rho(s)  E_Q(s,g)\,\d s \\
  &= \frac{1}{2\pi i} \int_{(\sigma)} T^s\hat{h}(-s) E_Q(s,g)\,\d s, 
 \end{align*}
 from which the result follows.
\end{proof}

From this we get the following direct consequences relating the analytic properties of the light-cone Eisenstein series to the smoothed counting function.
\begin{Prop}\label{pr:AC2SC}
Let $g\in G$ and $\alpha\leq \tfrac{n}{2}$. Assume that $E_Q(s,g)$ has an analytic continuation to $\Re(s)\geq \alpha$ with finitely many poles at $n=\sigma_0>\sigma_1> \ldots>  \sigma_k>\tfrac{n}{2}$ and satisfies that for any $\sigma\geq \alpha$ we have that  $\int_{-T}^T|E_Q(\sigma+it,g)|\,\d t$ grows polynomially in $T$. Then  
for any $h\in C^\infty_c(\R^+)$ there are $c_1,\ldots, c_k$ (depending on $h$ and $g$) such that 
\begin{equation}\label{eq:asymp}N_h(T;g)= \omega_Q \hat{h}(-n)T^n+\sum_{j=1}^k  c_j T^{\sigma_j}+O_{h,g}(T^{\alpha}).\end{equation}
Moreover, if \eqref{eq:asymp} holds {for any $h\in C_c^{\infty}(\R^+)$} then $E_Q(s,g)$ has an analytic continuation to $\Re(s)>\alpha$ with the only  poles at $\sigma_0,\ldots,\sigma_k$.
\end{Prop}
\begin{proof}
For a fixed $\sigma>n$ starting from the formula
\begin{align}\label{equ:smoothid}
 N_h(T;g)  &= \frac{1}{2\pi i} \int_{(\sigma)} T^s\hat{h}(-s) E_Q(s,g)\,\d s, 
 \end{align}
 and shifting the contour to the left  up to the line $\Re(s)=\alpha$ 
 and picking up the poles at $s=\sigma_j$ for $0\leq j\leq k$ we get that 
 \begin{eqnarray*}
 N_h(T;g)  =\sum_{j=0}^k c_jT^{\sigma_j}+ \frac{1}{2\pi i} \int_{(\alpha)} T^s\hat{h}(-s) E_Q(s,g)\,\d s, \\
 \end{eqnarray*}
 with $c_j=\hat{h}(-\sigma_j) \Res_{s=\sigma_j}E_Q(s,g)$. In particular, $c_0=\hat{h}(-n)\omega_Q$ (see \eqref{equ:omeq}) and here the shift of contour can be justified using the polynomial growth assumption on $\int_{-T}^T|E_Q(\sigma+it,g)|\,\d t$ and the super-polynomial decay of $|\hat{h}(-\sigma+it)|$ in $t$; see e.g. \cite[Proposition 4.2]{BurrinNevoRuhrWeiss2020}. We can now put absolute values in the last integral and bound it by $O_{h,g}(T^{\alpha})$.

For the other direction, we assume that $N_h(T;g)$ satisfies 
$$N_h(T;g)= \omega_Q \hat{h}(-n)T^n+\sum_{j=1}^k  c_j T^{\sigma_j}+O_{h,g}(T^{\alpha}),$$
and let
$$\cE(T)=N_h(T;g)-\sum_{j=0}^k  c_j T^{\sigma_j}.$$ 
Noting that  there is $T_0>0$  (depending on $h$ and $g$) such that $N_h(T;g)=0$ for all $T\leq T_0$ we have that for $\Re(s)>n$
\begin{align*}
\hat{N}_h(s;g)&=\int_{T_0}^\infty N_h(t;g)t^{-s}\frac{\d t}{t}\\
&= \sum_{j=0}^k  c_j  \frac{T_0^{\sigma_j-s}}{s-\sigma_j}+\int_{T_0}^\infty \cE(T)T^{-s}\frac{\d T}{T}.
\end{align*}
The bound $|\cE(T)|\ll_{h,g} T^\alpha$ implies that $\int_{T_0}^\infty \cE(T)T^{-s}\frac{\d T}{T}$ absolutely converges and hence analytic for $\Re(s)>\alpha$, and hence 
$\hat{N}_h(s;g)=\hat{h}(-s) E_{Q}(s,g)$ is analytic in the half plane $\Re(s)>\alpha$ with the only poles at $n=\sigma_0, \sigma_1,\ldots, \sigma_k$. Since this holds for any {compactly supported} smooth function $h$ we can conclude that $E_{Q}(s,g)$ has no poles with $\alpha<\Re(s)<\tfrac{n}{2}$.
\end{proof}

\begin{rem}
If one considers all integer solutions and not just the primitive ones, given by 
$$N^*_{h}(T;g):=\sum_{\bm{v}\in \cV^+_{Q}(\Z)} h\left(T^{-1} \|\bm{v}g\|_K\right),$$
one can follow the same argument using the completed series 
$$E_Q^*(s,g):=\zeta(s)\|\bm{e}_0\|_K^{-s}E_Q(s,g)=\sum_{\bm{v}\in \cV^+_{Q}(\Z) }\|\bm{v}g\|_K^{-s}.$$
When $n\geq 3$ the argument is identical and the only difference is that the main term is multiplied by $\zeta(n)\|\bm{e}_0\|_K^{-n}$ (and similarly the secondary terms).
For $n=2$, there is a pole at $s=1=\frac{n}{2}$, however, the contribution of this pole can be easily dealt with and it does not change the asymptotic formula. When $n=1$ the situation is a bit different as the extra factor of $\zeta(s)$ gives a double pole at $s=1$ resulting in a main term of order $T\log(T)$ instead of $T$.
\end{rem}

{When $Q=Q_n$ with $Q_n$ given as in \eqref{e:Qn}, \thmref{t:smoothcounting} 
 follows from 
 \propref{pr:AC2SC} (using also the correspondence described in \rmkref{rmk:corres}), but with a weaker bound of $O(T^{\frac{n}{2}})$. For $Q=Q_{n,d}$ as in \eqref{e:Qnd} with $d\in \N$ odd and square-free, we will state in the next section an explicit Fourier expansion formula for $E_{Q_{n,d}}(s,g)$ (see \thmref{thm:fouexpan}) from which we can get more explicit descriptions of  the poles of $E_{Q_{n,d}}(s,g)$ to the left of $\Re(s)=\tfrac{n}{2}$. Assuming this Fourier expansion for now we have the following refinement of  \propref{pr:AC2SC} giving the improved bound on the remainder claimed in \thmref{t:smoothcounting}.}

{ \begin{Prop}\label{pr:AC2SC2}
Let $Q=Q_{n,d}$ with $d$ odd and square-free. For any $h\in C^\infty_c(\R^+)$ and any $k>0$ the smoothed counting function given in \eqref{e:smoothed} satisfies
  $$N_h(T; g)= \omega_{Q_{n,d}}\hat{h}(-n)T^n+c_n\hat{h}\left(-\tfrac{n+1}{2}\right)T^{\frac{n+1}{2}}+O_{h,k}(T^{\frac{n}{2}}\log(T)^{-k}),$$
 where $c_n=0$ unless $n>1$ satisfies  $n\equiv 1\Mod{8}$. 
 \end{Prop}}
 \begin{proof}
Following the same argument as in the proof of \propref{pr:AC2SC} we start with 
$$ N_h(T;g)  = \frac{1}{2\pi i} \int_{(\sigma_0)} T^s\hat{h}(-s) E_Q(s,g)\,\d s,$$
with $\sigma_0>n$. From {the explicit Fourier expansion formula of $E_Q(s,g)$ proved in  \thmref{thm:fouexpan} below}, using known {bounds} for zero-free regions of zeta and Dirichlet $L$-functions we see that $E_Q(s,g)$ has no poles in the region $\{s=\sigma+it:  \frac{n}{2}\geq \sigma>\frac{n}{2}-\frac{c}{\log(|t|+2)}\}$, for some absolute constant $c>0$ (that may depend on $d$). Fix a large parameter $R$ and consider the contour $\cC_R=\cC_R^{-}\cup \cC_R^{0}\cup \cC_R^+$ with 
 $\cC_R^-=\{\frac{n}{2}-\frac{c}{\log R}+it: |t|\leq R\}$, $\cC_R^+= \{\frac{n}{2}+it: |t|\geq R\}$ and $\cC_R^0=\{\sigma\pm iR: \frac{n}{2}-\frac{c}{\log{R}}\leq  \sigma \leq \frac{n}{2}\}$.
 
 We can now shift the contour from $\Re(s)=\sigma_0$ to $\cC_R$ picking up the pole at $s=n$ giving the main term and (when $n\equiv 1\pmod{8}$ also the pole at $\frac{n+1}{2}$ giving the secondary term). It remains to bound the integral over $\cC_R$.
Using the fast decay of ${\hat{h}}(\sigma+it)$ {in $t$} and the fact that $E_Q(\sigma+it,g)$ grows polynomially with $t$ (which follows from the Fourier expansion) we can bound the contribution of $\cC_R^-$ by $O_h(T^{\frac{n}{2} -\frac{c}{\log(R)}})$ and the contribution of $\cC_R^0\cup \cC_R^+$ by $O_{h,k}(T^{\frac{n}{2}}R^{-k})$. Taking $R=\log(T)$  and noting that {for all sufficiently large $T$ we have that $c\log(T)\geq k (\log\log(T))^2$} we get a remainder of order 
 $O_{h,k}(T^{\frac{n}{2}}\log(T)^{-k})$ as claimed. 
 \end{proof}

We now compare these results to the results of \cite{Getz2018,Tran2020} using the new version of the circle method developed in \cite{Brown1996} to prove  \thmref{t:circle}. 
\begin{proof}[{Proof of  \thmref{t:circle}}]
The results of \cite{Getz2018,Tran2020} deal with the more general smooth counting function 
 $$\sum_{\bm{v}\in \cV_Q(\Z)} f(T^{-1}\bm{v}),$$
 for $f\in C^\infty_c(\R^{n+2})$ and $Q(\bm{v})=\bm{v}J\bm{v}^t$ a rational indefinite from.  Applying their result with $f(\bm{v})=h(\|\bm v g\|_K)$  implies that for $n\geq 2$
$$N^*_h(s;g)=c_1T^n+c_2T^{\left\lfloor{\frac{n+2}{2}}\right\rfloor}+O_{h,g}(T^{\frac{n}{2}+\epsilon}),$$
while for $n=1$
$$N^*_h(s;g)=c_1T\log(T)+c_2T+O_{h,g}(T^{\frac{1}{2}+\epsilon}),$$
with $c_1,c_2$ constants that may depend on $h$ and $g$. Moreover, noting that $\det(J)<0$ for any $Q$ of signature $(n+1,1)$, \cite[Theorem 1.1]{Getz2018} implies that $c_2=0$ when $n\equiv 2\Mod{4}$ as well as when $n\equiv 0\Mod{4}$ and $-\det(J)$ is not a square of a rational number. Combining this result with our estimate above concludes the proof.
\end{proof} 

\begin{rem}
We note that for the standard form $Q_n$ we have that $\det(J)=-1$, hence, the result {of \thmref{c:poles1}} that there are no exceptional poles when $n\equiv 0\Mod{4}$ (as well as for odd $n\not\equiv 1\Mod{8}$) does not follow from this more general result  and it is not clear to us how to deduce it without looking at the Fourier expansion. 
\end{rem}

By localizing the smoothing function, we can get results for the sharp cutoff counting function $N(T;g)$ {which in particular implies \thmref{t:supnorm2counting} when choosing $Q=Q_n$; see also \rmkref{rmk:extraassum}}. 
\begin{Prop}\label{p:supnorm2counting}
Let $\nu>0$ be such that  $\int_{-T}^T|E_Q(\tfrac{n}{2}+it,g)|^2\, \d t\ll_g T^{2\nu+1}$ for all $T>0$. {If $n\geq 3$ further assume that $\nu\geq \frac{1}{n-2}$}.
Then
$$N(T;g)=\frac{\omega_Q }{n}T^n+O_g(T^{n-\frac{n}{2(\nu+1)}}).$$
\end{Prop}
\begin{proof}
Fix $h\in C^\infty_c(\R^+)$ supported on $[\tfrac{1}{e},e]$ with $\int_0^{\infty} h(y)\frac{\d y}{y}=1$ (so that $\hat{h}(0)=1$). For a large parameter $M\geq 1$ let $h_M(y)=Mh(y^M)$  and note that
$\hat{h}_M(s)=\hat{h}(\frac{s}{M})$. For any $T>1$ let $\chi_T$ be the indicator function of $(0, T)$. Then we have that $\chi_T*h_M(y)=\chi_1*h_M(\frac{y}{T})=\left\lbrace\begin{array}{ll}1 & y< e^{-1/M}T,\\ 0 & y>e^{1/M}T.\end{array}\right.$ Here $\chi*h(y):=\int_0^{\infty}\chi(x)h(\tfrac{y}{x})\ \frac{\d x}{x}$ is a convolution. 
Hence the corresponding smoothed counting functions satisfy 
\begin{align}\label{equ:uplowbd}
N_{\chi_1*h_M}(Te^{-{1/M}};g)\leq N(T;g)\leq N_{\chi_1*h_M}(Te^{1/M};g).
\end{align}

On the other hand, starting from identity for the smoothed counting function (cf. \eqref{equ:smoothid}) and shifting the contour of integration, picking up the potential pole at $\sigma_1=\left\lfloor{\frac{n+2}{2}}\right \rfloor$ if it exists {(which could only happen when $n\geq 3$)}, we get that 
\begin{align*}
N_{\chi_1*h_M}(T;g)&=\frac{1}{2\pi i}\int_{(\sigma)}\frac{T^s}{s}\hat{h}_M(-s)E_Q(s,g)\,\d s\\
&= \frac{\omega_Q T^n}{n}\hat{h}_M(-n)+ \hat{h}_M(-\sigma_1)\vf_1(g)\frac{T^{\sigma_1}}{\sigma_1}+\frac{1}{2\pi}\int_\R \tfrac{\hat{h}_M(-\tfrac{n}{2}-it) T^{\frac{n}{2}+it}}{\frac{n}{2}+it}E_Q(\tfrac{n}{2}+it,g)\,\d t,
\end{align*}
with $\vf_1(g)=\Res_{s=\sigma_1}E_Q(s,g)$.
We can estimate $\hat{h}_M(-\sigma_j)=\hat{h}(-\frac{\sigma_j}{M})=1+O_h(\frac{1}{M})$ and for the integral, putting absolute values 
we can bound
\begin{align*}
\left|\frac{1}{2\pi}\int_\R \tfrac{\hat{h}_M(-\tfrac{n}{2}-it) T^{\frac{n}{2}+it}}{\frac{n}{2}+it}E_Q(\tfrac{n}{2}+it,g)\,\d t  \right| &\ll T^{\frac{n}{2}}\int_\R \tfrac{\left| \hat{h}(-\tfrac{n}{2M}-\frac{it}{M})\right|}{(\tfrac{n^2}{4}+t^2)^{\frac12}} \left|E_Q(\tfrac{n}{2}+it,g)\right|\,\d t\\
&=T^{\frac{n}{2}}M\int_{\R}\tfrac{\left|\hat{h}(-\tfrac{n}{2M}-it)\right|}{(\frac{n^2}{4}+M^2t^2)^{\frac12}}\left|E_Q(\tfrac{n}{2}+iMt,g)\right|\,\d t\\
&= T^{\frac{n}{2}}M\sum_{k\in \Z}\int_{k}^{k+1}\tfrac{\left|\hat{h}(-\tfrac{n}{2M}-it)\right|}{(\frac{n^2}{4}+M^2t^2)^{\frac12}}\left|E_Q(\tfrac{n}{2}+iMt,g)\right|\,\d t\\
&\ll_{h,N} T^{\frac{n}{2}}\sum_{k\in \Z} (1+|k|)^{-N}\int_{k}^{k+1}\left|E_Q(\tfrac{n}{2}+iMt,g)\right|\,\d t\\
&= T^{\frac{n}{2}}M^{-1}\sum_{k\in \Z} (1+|k|)^{-N}\int_{Mk}^{M(k+1)}\left|E_Q(\tfrac{n}{2}+it,g)\right|\,\d t\\
& \leq T^{\frac{n}{2}}M^{-1}\sum_{k\in \Z} \frac{M^{1/2}}{(1+|k|)^{N}} \left(\int_{Mk}^{M(k+1)}\left|E_Q(\tfrac{n}{2}+it,g)\right|^2\,\d t\right)^{1/2}\\
&\ll T^{\frac{n}{2}}M^\nu, 
\end{align*}
{where in the last line we used the second moment bound for the Eisenstein series and the assumption that $N>\nu+2$ so {that} the series over $k$ converges.}
We thus get that
\begin{align*}
	N_{\chi_1*h_M}(T;g)&=\frac{\omega_Q T^n}{n}+O_g(T^nM^{-1}+T^{\left\lfloor{\frac{n+2}{2}}\right \rfloor}+T^{\frac{n}{2}}M^{\nu}).
\end{align*}
Combining this estimate with \eqref{equ:uplowbd} with the choice of $M=T^{\frac{n}{2(\nu+1)}}$ concludes the proof. 
\end{proof}
\begin{remark}\label{rmk:extraassum}
{The assumption that $\nu\geq \frac{1}{n-2}$ when $n\geq 3$ is only needed when the potential exceptional pole at $s={\left\lfloor{\frac{n+2}{2}}\right \rfloor}$ exists and to ensure that the term $T^{\left\lfloor{\frac{n+2}{2}}\right \rfloor}$ above appearing in the error is dominated by $T^{n-\frac{n}{2(\nu+1)}}$. For $Q=Q_n$ such a pole occurs exactly when $n\geq 9$ with $n\equiv 1\Mod{8}$, thus only in this case we need $n-\frac{n}{2(\nu+1)} \geq \frac{n+1}{2}$ which is equivalent to the assumption stated in \thmref{t:supnorm2counting}.}
\end{remark}
 \begin{remark}\label{r:lowerbound}
 Using the bound $\int_{-T}^T |E_i(\tfrac{n}{2}+it,g)|^2\,\d t \ll_g T^{n+1}$ for general Eisenstein series (see \rmkref{r:WboundE}) we can take $\nu=\frac{n}{2}$ above giving an error term of $O(T^{n-\frac{n}{n+2}})$.
Recall {from \rmkref{rmk:corres}} that for $Q=Q_{n,d}$ we have that $N(T;1_G)=N(T+\tfrac{1}{2};1_G)$ for any $T\in\N$ and hence the best bound one can hope for the remainder is $O(T^{n-1})$. 
 This shows that for these lattices we always have $\nu\geq \frac{n}{2}-1$ giving a lower bound on the (mean square of the) sup-norm of these Eisenstein series. 
 
  \end{remark}
	
\section{Fourier expansion and its applications}\label{s:Fourier}
In this section we state our main result which is the explicit Fourier expansion for the light-cone Eisenstein series.
Here we treat the more general case of the quadratic form $Q=Q_{n,d}$ given in \eqref{e:Qnd} with $d$ odd and square-free (where the case of $d=1$ was stated in the introduction).  
After stating our result we prove some consequences for the location of poles, bounding the sup-norm and counting rational points on the sphere.

{We follow the notation introduced in section \ref{sec:cormea1}.} For $Q=Q_{n,d}$ we denote by $G=\SO_Q^+(\R)$ and by $\G=G\cap \SL_{n+2}(\Z)$ the lattice of integer points. 
We also recall that here $\bm{e}_0=(-d,0,\ldots, 0,1)\in \R^{n+2}$, the group $P\leq G$ is the group stabilizing the line spanned by $\bm{e}_0$
 and  $U<P$ is its unipotent radical. Moreover, 
$$
\G_U=\G\cap U=U(\Z)=\left\{u_{\bm{x}} : \bm{x}\in \Lambda\right\}
	$$
	with 
	$\Lambda$ as in \eqref{equ:lambdalatt}. Let $\Lambda^{\ast}$ be the dual lattice of $\Lambda$. Explicitly, 
	\begin{align*}
		\Lambda^{\ast}=\left\{\tfrac{\bm{x}}{2}: \bm{x}\in \Z^n,\ x_1\equiv x_2\equiv \cdots \equiv x_n \Mod{2}\right\}.
	\end{align*}
	For any $\Re(s)>n$, {$g\in G$} and $\bm{\lambda}\in \Lambda^{\ast}$ let us define
	\begin{align*}
		a_Q(s, g; \bm{\lambda}):=\frac{1}{\vol\left(\R^n/\Lambda\right)}\int_{\R^n/\Lambda}E_Q(s,u_{\bm{x}}g)e(-\bm{\lambda}\cdot\bm{x})\,\d\bm{x},
	\end{align*}
so that 
	\begin{align}\label{equ:fourierexpan}
	E_Q(s,u_{\bm{x}}g)=\sum_{\bm{\lambda}\in \Lambda^{\ast}}a_Q(s,g; \bm{\lambda}) e(\bm{\lambda}\cdot\bm{x}).
\end{align}
Here for any $t\in\R$, $e(t):=e^{2\pi i t}$ {and for any $\bm{x}\in \R^n$, $u_{\bm{x}}$ is as in \eqref{equ:uxay}}. 
Our main result is the following explicit Fourier expansion formulas.

\begin{Thm}\label{thm:fouexpan}
	Assume $d$ is odd and square-free. Let $g=u_{\bm{x}}a_y k$ with $u_{\bm{x}}\in U$, $a_y\in A$, and $k\in K$ {as in section \ref{sec:explicit}}. Then we have
	\begin{align}\label{equ:forcoefcons}
		a_{Q_{n,d}}(s,g; \bm{0})=y^s+\Phi_{n,d}(s)y^{n-s},
	\end{align}
	where
	$$\Phi_{n,d}(s):=\e_{n,d}(s)\left\lbrace\begin{array}{ll}
		\frac{\Gamma(\frac{2s-n+2}{4})^2}{\G(\frac{s+1}{2})\G(\frac{s-n+1}{2})} \frac{\xi(s-n+1)\xi(s-\frac{n}{2})}{\xi(s)\xi(s-\frac{n}{2}+1)} & n\equiv 0\Mod{4},\\
		\frac{\G(\frac{2s-n}{4})\G(\frac{2s-n+4}{4})}{\G(\frac{s+1}{2})\G(\frac{s-n+1}{2}) }\frac{\xi(s-n+1)L^*(s-\frac{n}{2},\chi_{-4})}{\xi(s) L^*(s-\frac{n}{2}+1,\chi_{-4})} & n\equiv 2\Mod{4},\\
		\frac{\G(\frac{2s-n+1}{4})\G(\frac{2s-n+3}{4})}{\G(\frac{s+1}{2}) \Gamma(\frac{s-n+1}{2})}\frac{\xi(s-n+1)\xi(2s-n)}{\xi(s)\xi(2s-n+1)} & n\equiv 1\Mod{2}.
	\end{array}\right.
	$$
	Here $\e_{n,d}(s):=d^{s-n}\prod_{p\mid 2d}\e_{n}^{(p)}(s)$ with $\e_n^{(2)}(s)$ as given in \eqref{equ:e2factor}
	and for any prime $p\mid d$
\begin{align*}
\e_{n}^{(p)}(s):=\left\lbrace\begin{array}{ll}
		\frac{p^{n-1}-p^{\frac{n}{2}-1}+p^{n-s}-p^{2n-1-2s}}{1-p^{\frac{n}{2}-1-s}}	 & n\equiv 0 \Mod{4},\\
		\frac{p^{n-1}-\chi_{-4}(p)p^{\frac{n}{2}-1}+p^{n-s}-p^{2n-1-2s}}{1-p^{\frac{n}{2}-1-s}} &    n\equiv 2 \Mod{4},\\
			\frac{\chi_{-4}(p)p^{\frac{n-1}{2}}+p^{n-1}+p^{n-s}-p^{2n-1-2s}}{1+\chi_{-4}(p)p^{\frac{n-1}{2}-s}}  & n \equiv 1\Mod{4},\\
			\frac{p^{\frac{n-1}{2}}+p^{n-1}+p^{n-s}-p^{2n-1-2s}}{1+p^{\frac{n-1}{2}-s}}  & n \equiv 3\Mod{4}.
		\end{array}\right.
\end{align*}
For the non-constant coefficients we have for any $\bm{\lambda}\in\Lambda^*\setminus\{\bm{0}\}$,
	\begin{align}\label{equ:forcoefnoncons}
		a_{Q_{n,d}}(s,g; \bm{\lambda})=\frac{2^{s-n+1}\pi^s\|\bm{\lambda}\|^{s-\frac{n}{2}}}{d^{\frac{n}{2}}\G(s)\zeta(s)}\Phi_{n,d}(s;\bm{\lambda})y^{\frac{n}{2}}K_{s-n/2}(2\pi\|\bm{\lambda}\|yd^{-1}),
	\end{align}
	where $K_{s}(y)$ is the Bessel function of the second kind and
	\begin{align}\label{equ:phinonconst}
		\Phi_{n,d}(s;\bm{\lambda})&:=\sum_{d_1\mid d}\frac{1}{\varphi(d_1)}\sum_{\chi\Mod{d_1}}\e_n(s;\chi, \bm{\lambda})\left\lbrace\begin{array}{ll}
			\frac{1}{L(s-\frac{n}{2}+1,\chi)} & n\equiv 0\Mod{4},\\
			\frac{1}{L(s-\frac{n}{2}+1,\chi\chi_{-4})} & n\equiv 2\Mod{4},\\
			\frac{L\left(s-\frac{n-1}{2},\chi\chi_{D}\right)}{L(2s-n+1,\chi^2)} & n\equiv 1\Mod{2}.
		\end{array}\right.
	\end{align}
	Here when $n$ is odd $D:=(-1)^{\frac{n-1}{2}}\|2\bm{\lambda}\|^2$ and {$\chi_D$ is the unique primitive character inducing $(\frac{D}{\cdot})$}, $\varphi$ is the Euler's totient function, $\e_n(s;\chi, \bm{\lambda})$ is some function which is holomorphic off the line $\Re(s)= \left \lfloor{\frac{n-1}{2}}\right \rfloor$ with poles contained in a finite union of periodic points on $\Re(s)= \left \lfloor{\frac{n-1}{2}}\right \rfloor$ depending on $d$ and bounded in the half plane $\Re(s)\geq \frac{n}{2}$ by some power of $\|\bm\lambda\|$.
	Moreover, when $n$ is odd and $n\not\equiv 1\Mod{8}$ then $\e_n(\frac{n+1}{2};\chi,\bm{\lambda})=0$ when $\chi\chi_{D}$ is principal.	\end{Thm}

For the special case of $d=1$ we can say more about the Fourier coefficients. 
In this case $\chi\Mod{d}$ is trivial and we abbreviate the notation $\Phi_{n,1}(s), \Phi_{n,1}(s;\bm{\lambda})$ and $\e_n(s;\chi,\bm{\lambda})$ (with $\bm{\lambda}\in \Lambda^*\setminus\{\bm{0}\}$) by $\Phi_n(s), \Phi_{n}(s;\bm{\lambda})$ and $\e_n(s;\bm{\lambda})$ respectively. \begin{Thm}\label{thm:funequbofo}
Keep the notation as in \thmref{thm:fouexpan}. When $d=1$ we have for any $\bm{\lambda}\in \Lambda^*\setminus\{\bm{0}\}$
	\begin{align}\label{equ:phinonconst1}
		\Phi_{n}(s;\bm{\lambda})&=\e_n(s; \bm{\lambda})\left\lbrace\begin{array}{ll}
			\frac{1}{\zeta(s-\frac{n}{2}+1)} & n\equiv 0\Mod{4},\\
			\frac{1}{L(s-\frac{n}{2}+1,\chi_{-4})} & n\equiv 2\Mod{4},\\
			\frac{L\left(s-\frac{n-1}{2},\chi_{D}\right)}{\zeta(2s-n+1)} & n\equiv 1\Mod{2},
		\end{array}\right.
	\end{align}
	 where the factor $\e_n(s; \bm{\lambda})$ 
	satisfies the functional equation 
	\begin{align}\label{equ:funequd1}
	\e_n(n-s;\bm{\lambda})=\mathfrak{s}_n\e_n(s;\bm{\lambda})\left\{\begin{array}{ll}
				\frac{2^{n-2s}(1-2^{\frac{n}{2}-1-s})}{1-2^{s-1-\frac{n}{2}}}\|2\bm{\lambda}\|^{2s-n} & n\equiv 4\Mod{8},\\
				2^{\frac{n}{2}-s}\|2\bm{\lambda}\|^{2s-n} & n\equiv 2\Mod{4},\\
				\frac{\left(1+\mathfrak{s}_n2^{\frac{n-1}{2}-s}\right)2^{\frac{n+1}{2}-s}q^{\frac{n}{2}-s}\|2\bm{\lambda}\|^{2s-n}}{1+\mathfrak{s}_n2^{\frac{n+1}{2}-s}} & n\equiv 1\Mod{2}.
			\end{array}\right.
\end{align}
Here $\mathfrak{s}_n$ is the sign function given by 
\begin{align}\label{equ:sn}
\mathfrak{s}_n&:=\left\{\begin{array}{ll}
				1 & n\equiv 0\Mod{4},\\
				(-1)^{\frac{n-2}{4}} & n\equiv 2\Mod{4},\\
				(-1)^{\frac{n^2-1}{8}} & n\equiv 1\Mod{2}.
			\end{array}\right.
\end{align}
and for the $n$ odd case, $q$ is the modulus of $\chi_{D}$.
Moreover, when $n$ is even, for any $\bm\lambda\in \Lambda^*_{\rm pr}$ and $a\in \N$ we have  
\begin{align}\label{equ:crilinebd}
|\epsilon_n(\tfrac{n}{2}+it; a\bm{\lambda})|\ll_{\e} a^{\frac{n}{2}-1+\e}\|\bm\lambda\|^\e.
\end{align}
\end{Thm}

\subsection{Poles and exceptional poles} 
The poles of the constant term of the light-cone Eisenstein series can be easily read off from our formula.  In particular,  taking residue at the trivial pole, $s=n$, we can explicitly calculate the residue 
$\omega_{Q_{n,d}}=\Res_{s=n}E_{Q_{n,d}}(s,g).$
\begin{Thm}
We have that 
$\omega_{Q_{n,d}}:=\omega_{Q_{n}}\prod_{p|d}\epsilon_n^{(p)}(n)$ with $\omega_{Q_{n}}$ as in \eqref{eq:omegaQ} and $\e_n^{(p)}(s)$ as in \thmref{thm:fouexpan}. {\qed}
\end{Thm}

Regarding the exceptional poles, since in general there is more than one cusp, to find all poles we need make sure that the non-constant terms do not have additional poles, which can be verified from our formulas of the non-constant coefficients. This can be summarized in the following.
\begin{Thm}\label{thm:exppoles}
The set of exceptional poles of $E_{Q_{n,d}}(s,g)$ is empty if $n=1$ or $n\not\equiv 1\Mod{8}$ and equals $\{\frac{n+1}{2}\}$ if  $n\equiv 1\Mod{8}$ with $n\geq 9$.
\end{Thm}
\begin{proof}
Since the residue at any $s\in (\tfrac{n}{2},n)$ is a Laplacian eigenfunction with nonzero eigenvalue, to show there is no pole it is enough to show that all non-constant Fourier coefficients are analytic in this region.
From \eqref{equ:forcoefnoncons} we see that any pole in this range must come from a pole of  $\Phi_{n,d}(s;\bm{\lambda})$. Using  \eqref{equ:phinonconst}, for even $n$ it is not hard to see that $\Phi_{n,d}(s;\bm{\lambda})$ has no poles for $\Re(s)\geq \tfrac{n}{2}$.  For $n\equiv 1\Mod{4} $, we could potentially have a (simple) pole at $s=\tfrac{n+1}{2}$ coming from $L\left(s-\frac{n-1}{2},\chi\chi_{\|2\bm{\lambda}\|^2}\right)$ when $\chi\chi_{\|2 \bm\lambda\|^2}$ is a principal character. Similarly, when $n\equiv 3\Mod{4}$ we could have such a pole when $\chi\chi_{-\|2\bm\lambda\|^2}$ is principal. However, unless $n\equiv 1\Mod{8}$, these poles are cancelled by the zero of $\e_n(s;\chi,\bm{\lambda})$ at $s=\frac{n+1}{2}$; hence there are no exceptional poles in these cases. When $n=1$, this pole at $s=\frac{n+1}{2}=n$ is the trivial pole. When $n\equiv 1\Mod{8}$ with $n\geq 9$ we do in fact have an exceptional pole at $s=\tfrac{n+1}{2}$,
as we can see from looking at residue of the constant term at this point. \end{proof}
Regarding the poles in the region $\Re(s)<\tfrac{n}{2}$ we have the following generalization of \thmref{c:poles2}. 
\begin{Thm}\label{c:poles3}
For $n>1$, $d\geq 1$ {odd and square-free and $\Re(s)<\frac{n}{2}$}, the light-cone Eisenstein series $E_{Q_{n,d}}(s,g)$ has poles in the strip $0<\Re(s)<1$ located at the zeros of the zeta function $\xi(s)$.
 In addition, for any $d_1|d$ and Dirichlet character $\chi$ modulo $d_1$ we can have additional poles in the strip $\tfrac{n}{2}-1<\Re(s)<\tfrac{n}{2}$, corresponding to zeros of $L(2s-n+1,\chi^2)$ (when $n$ is odd) $L(s-\tfrac{n}{2}+1,\chi)$ (when $n\equiv 0\Mod{4}$) or $L(s-\tfrac{n}{2}+1,\chi\chi_{-4})$ (when $n\equiv 2\Mod{4}$), as well as periodic poles on the line $\Re(s)=\left \lfloor{\frac{n-1}{2}}\right \rfloor$ coming from the $\e$-factors {and finitely many real poles from the $\G$-factors}.   When $n=1$, the same holds except there are no poles at the zeros of $\xi(s)$. {\qed}
\end{Thm}

\subsection{Functional Equation}
We now use our results on the Fourier coefficients for the  $d=1$ case (see \thmref{thm:funequbofo}) to prove the functional equation for $E_{Q_{n,1}}(s,g)$ for all $n\not\equiv 0\Mod{8}$. 

\begin{proof}[Proof of \thmref{thm:funequ} assuming \thmref{thm:funequbofo}]
Let $Q=Q_{n,1}$ and let $g=u_{\bm{x}}a_yk$ be as in \thmref{thm:funequbofo}. By comparing Fourier coefficients we see that the functional equation \eqref{e:funeq} is equivalent to the functional equations $\Phi_{n}(n-s)\Phi_{n}(s)=1$ and that for any nonzero $\bm{\lambda}\in \Lambda^*$,
\begin{align}\label{equ:funeqpf}
\Phi_{n}(s)a_Q(n-s,g;\bm{\lambda})=a_Q(s,g;\bm{\lambda}).
\end{align}
For the first functional equation, let 
\begin{align*}
\Upsilon_n(s):=\left\lbrace\begin{array}{ll}
		\frac{\Gamma(\frac{2s-n+2}{4})^2}{\G({\frac{s+1}{2}})\G({\frac{s-n+1}{2}})}& n\equiv 0\Mod{4},\\
		\frac{\G(\frac{2s-n}{4})\G(\frac{2s-n+4}{4})}{\G(\frac{s+1}{2})\G(\frac{s-n+1}{2})} & n\equiv 2\Mod{4},\\
		\frac{\G(\frac{2s-n+1}{4})\G(\frac{2s-n+3}{4})}{\G(\frac{s+1}{2}) \Gamma(\frac{s-n+1}{2})} &n\equiv 1\Mod{2},
	\end{array}\right.
\end{align*}
be the $\G$-factors in the formula of $\Phi_{n}(s)$. From this and the functional equations $\xi(1-s)=\xi(s)$ and $L^*(1-s,\chi_{-4})=L^*(s,\chi_{-4})$ we get that
\begin{align*}
\Phi_{n}(n-s)\Phi_{n}(s)&=\e_n^{(2)}(n-s)\e_n^{(2)}(s)\Upsilon_n(n-s)\Upsilon_n(s).
\end{align*}
One easily checks that $\e_n^{(2)}(n-s)\e_n^{(2)}(s)=1$ for $n\not\equiv 0\Mod{8}$. 
On the other hand, using the recursive relation $\G(s+1)=s\G(s)$ we see that $\Upsilon_n(s)$ has the following alternative expression as a rational function
\begin{align}\label{equ:ratexpalt}
	\Upsilon_{n}(s)=\left\lbrace\begin{array}{ll}
		\prod_{i=1}^{n/4}\frac{s+1-\frac{n}{2}-2i}{s+1-2i} & n\equiv 0\Mod{4},\\
		\prod_{i=1}^{(n-2)/4}\frac{s-\frac{n}{2}-2i}{s+1-2i} & n\equiv 2\Mod{4},\\
		\prod_{i=1}^{\frac{n-1}{2}}\frac{s-2i}{s-i} & n\equiv 1\Mod{2}, 
	\end{array}\right. 
\end{align}
from which one easily checks that $\Upsilon_n(s)$ satisfies $\Upsilon_n(n-s)\Upsilon_n(s)=1$. We have thus finished the proof of the first functional equation.

Next, fix an arbitrary nonzero $\bm{\lambda}\in \Lambda^*$, we need to prove the functional equation \eqref{equ:funeqpf}. For this we show that \eqref{equ:funeqpf} is equivalent to the functional equation \eqref{equ:funequd1}. Let us first treat the case when $n\equiv 4\Mod{8}$. In this case applying formulas \eqref{equ:forcoefnoncons}, \eqref{equ:phinonconst1} and completing the zeta function we get that 
\begin{align*}
a_Q(s,g;\bm{\lambda})&=2^s\pi^{\frac{n-2}{4}}\|\bm{\lambda}\|^{s-\frac{n}{2}}\frac{\G(\tfrac{s}{2})\Gamma(\tfrac{2s+2-n}{4})}{\Gamma(s)}\frac{\e_n(s;\bm{\lambda})}{\xi(s)\xi(s+1-\frac{n}{2})}y^{\frac{n}{2}}K_{s-n/2}(2\pi\|\bm{\lambda}\|y).
\end{align*}
Now using again the relation $\G(s+1)=s\G(s)$ and Legendre duplication formula $\G(s)\G(s+\tfrac12)=2^{1-2s}\sqrt{\pi}\G(2s)$ we have 
\begin{align*}
\frac{\G(\tfrac{s}{2})\Gamma(\tfrac{2s+2-n}{4})}{\Gamma(s)}=2^{1+\frac{n}{4}-s}\sqrt{\pi}\prod_{i=1}^{n/4}(s+1-2i)^{-1}.
\end{align*}
From this formula and using the above rational function expression of $\Upsilon_{n}(s)$  and the functional equation $\xi(1-s)=\xi(s)$ we see that \eqref{equ:funeqpf} is equivalent to 
$$
\e^{(2)}_n(s)\e_n(n-s;\bm{\lambda})=-\|\bm{\lambda}\|^{2s-n}\e_n(s;\bm{\lambda}).
$$
Since in this case $\e^{(2)}_n(s)=\frac{2^{s-\frac{n}{2}-1}-1}{1-2^{\frac{n}{2}-1-s}}$, the above equation is further equivalent to 
\begin{align*}
\e_n(n-s;\bm{\lambda})&=\frac{2^{n-2s}(1-2^{\frac{n}{2}-1-s})}{1-2^{s-\frac{n}{2}-1}}\|2\bm{\lambda}\|^{2s-n}\e_n(s;\bm{\lambda}).
\end{align*}
The case when $n\equiv 2\Mod{4}$ follows similarly using additionally the functional equation $L^*(1-s,\chi_{-4})=L^*(s,\chi_{-4})$. 

When $n$ is odd, let $D=(-1)^{\frac{n-1}{2}}\|2\bm{\lambda}\|^2$ be as above and let $q$ be the modulus of $\chi_D$. Let $a\in\{0,1\}$ be such that $a=0$ if $\chi_D(-1)>0$ (i.e. $n\equiv 1\Mod{4}$) and $a=1$ if $\chi_D(-1)=-1$ (i.e. $n\equiv 3\Mod{4}$) as before. Again applying \eqref{equ:forcoefnoncons} and \eqref{equ:phinonconst1} and completing zeta and corresponding $L$-functions (see \eqref{equ:comlfun}) we have 
\begin{align*}
	a_Q(s,g;\bm{\lambda})&=2^sq^{\frac{n-2s-1-2a}{4}}\pi^{\frac{n-1+2a}{4}}\|\bm{\lambda}\|^{s-\frac{n}{2}}\tfrac{\G(\tfrac{s}{2})\Gamma(s-\tfrac{n-1}{2})}{\Gamma(s)\Gamma(\tfrac{2s-n+1+2a}{4})}\tfrac{L^*(s-\frac{n-1}{2},\chi_D)}{\xi(s)\xi(2s-n+1)}\e_n(s;\bm{\lambda})y^{\frac{n}{2}}K_{s-n/2}(2\pi\|\bm{\lambda}\|y).
\end{align*}
Doing similar computation, this time {using also} the functional equation \eqref{equ:funequlfun}, we see that in this case \eqref{equ:funeqpf} is equivalent to (recalling that when $n$ is odd $\mathfrak{s}_n=(-1)^{\frac{n^2-1}{8}}$) 
$$
\e^{(2)}_n(s)\e_n(n-s;\bm{\lambda})=
\mathfrak{s}_n2^{2s-n}q^{\frac{n}{2}-s}\|\bm{\lambda}\|^{2s-n}\e_n(s;\bm{\lambda}),
$$
which is further equivalent to 
\begin{align*}
	\e_n(n-s;\bm{\lambda})
	&=\frac{\mathfrak{s}_n(1+\mathfrak{s}_n2^{\frac{n-1}{2}-s})2^{\frac{n+1}{2}-s}q^{\frac{n}{2}-s}\|2\bm{\lambda}\|^{2s-n}}{1+\mathfrak{s}_n2^{\frac{n+1}{2}-s}}\e_n(s;\bm{\lambda})
\end{align*}
as claimed. 
\end{proof}

\begin{remark}\label{rmk:fueqfail}
When $n\equiv 0\Mod{8}$, the above analysis implies that
\begin{align*}
\Phi_{n}(n-s)\Phi_{n}(s)&=\e_n^{(2)}(n-s)\e_n^{(2)}(s)=\frac{2^{\frac{n}{2}-1-s}}{1-2^{s-\frac{n}{2}-1}}\times \frac{2^{s-\frac{n}{2}-1}}{1-2^{\frac{n}{2}-1-s}}=\frac{1}{5-2^{\frac{n}{2}+1-s}-2^{s-\frac{n}{2}+1}}\neq 1.
\end{align*}
Nevertheless, for $s=\frac{n}{2}+it$,  we still have the upper bound
$$
|\Phi_{n}(s)|^2=\Phi_{n}(n-s)\Phi_{n}(s)=\left\lbrace\begin{array}{ll}
		1 & n\not\equiv 0\Mod{8},\\
		\frac{1}{5-2^{1+it}-2^{1-it}}\leq 1 & n\equiv 0\Mod{8}.
	\end{array}\right. 
$$ 
\end{remark}

\subsection{Supremum norms}
For the special case of $d=1$ we can prove the convexity bound for the {sup-norm}.
We give two different proofs, one using the functional equation (which is valid for all $n\not\equiv 0\Mod{8}$) and another using bounds for the Fourier coefficients (which is valid for even $n$).
Combining the two arguments we get the sup-norm bounds {on} $E_{Q_n}(\tfrac{n}{2}+it,{g})$ for all $n$ {claimed in \thmref{t:supnorm}}.
\begin{Prop}\label{p:supodd}
Let $g\in G$ be fixed. For any {even} $n\not\equiv 0\Mod{8}$ we have for all $|t|\geq 2$
$$|E_{Q_n}(\tfrac{n}{2}+it,g)|\ll_{g} |t|^{\frac{n}{2}}\log^2(|t|),$$
while for odd $n\geq 3$ we have
$$|E_{Q_n}(\tfrac{n}{2}+it,g)|\ll_{g} |t|^{\frac{3n}{4}}\log^2(|t|).$$
\end{Prop}
\begin{proof}
Given $n\geq 2$ {with $n\not\equiv 0\Mod{8}$}, for a fixed $g\in G$ consider the function
$F(s):=H(s)E_{Q_n}(s,g)$, where 
$$H(s):=\left\lbrace\begin{array}{ll}
		\frac{(1-2^{\frac{n}{2}-1-s})\xi(s)\xi(s-\frac{n}{2}+1)}{\Gamma(\frac{2s-n+2}{4})^2} & n\equiv 0\Mod{4},\\
		\frac{\xi(s) L^*(s-\frac{n}{2}+1,\chi_{-4})}{\G(\frac{2s-n}{4})\G(\frac{2s-n+4}{4})} & n\equiv 2\Mod{4},\\
		\frac{(1-2^{n-1-2s})\xi(s)\xi(2s-n+1)}{\G(\frac{3s-n}{2})} & n\equiv  1\Mod{2}.\\
	\end{array}\right.$$
Note that $H(s)$ has no poles off the real line and its zeros exactly cancel the poles of $E_{Q_n}(s,g)$ so that $F(s)$ has no poles for {$\Im(s)\neq 0$}.
Moreover, using the functional equation \eqref{e:funeq} we get the following functional equation for $F(s)$ 
\begin{align*}
F(n-s) &= F(s) \e_n^{(2)}(n-s) \left\lbrace\begin{array}{ll}
		  \frac{\Gamma(\frac{2s-n+2}{4})^2}{\G(\frac{n-s+1}{2})\G(\frac{1-s}{2})}\frac{1-2^{s-1-\frac{n}{2}}}{1-2^{\frac{n}{2}-1-s}}  & n\equiv 0\Mod{4},\\
		\frac{\G(\frac{2s-n}{4})\G(\frac{2s-n+4}{4})}{\G(\frac{n-s+1}{2})\G(\frac{1-s}{2}) } & n\equiv 2\Mod{4},\\
	{\frac{\G(\frac{n-2s+1}{4})\G(\frac{n-2s+3}{4})\Gamma(\frac{3s-n}{2})}{\G(\frac{n-s+1}{2}) \Gamma(\frac{1-s}{2})\G(n-\frac{3s}{2})}\frac{1-2^{2s-1-n}}{1-2^{n-1-2s}}}  & n\equiv 1\Mod{2}.\\
	\end{array}\right.
	\end{align*}
In particular, using Stirling's approximation for $s=\sigma+it$ with $|t|\geq 2$ and $n\leq \sigma\leq 2n$ we can bound 
\begin{align}\label{e:AFE}
|F(n-\sigma-it)| &\ll |F(\sigma+it)|  \left\lbrace\begin{array}{ll}
		  |t|^{2\sigma-n}  & n\equiv 0\Mod{2},\\
		|t|^{3(\sigma-\frac{n}{2})} & n\equiv 1\Mod{2}.	
	\end{array}\right.
	\end{align}
Now, for $\Re(s)=\sigma>n$ using the series expansion we can bound 
$$|E_{Q_n}(\sigma+it,g)|\leq |E_{Q_n}(\sigma,g)|\ll_g \frac{1}{\sigma-n}.$$
For the function $H(s)$, using Stirling's approximation for the  $\G$-functions and noting that $\zeta(s)$ and $L(s,\chi_{-4})$ are uniformly bounded  $\Re(s)\geq 2$ we can bound for $|t|\geq 2$
\begin{align}\label{equ:nebd}
|F(n+\e+it)|  \ll \e^{-1} |H(n+\e+it) | \ll \e^{-1}|t|^{\nu_n}
\end{align}
with 
$$\nu_n:= \left\lbrace\begin{array}{ll}
		 \frac{n-2}{4} & n\equiv 0\Mod{4},\\
		\frac{n}{4} & n\equiv 2\Mod{4},\\
	0 & n\equiv 1\Mod{2}.
	\end{array}\right.
	$$
Next, to the left of the line $\Re(s)=-\e<0$ we can use \eqref{e:AFE} with $\sigma=n+\e$ and \eqref{equ:nebd} to get that for $|t|>2$ 
$$|F(-\e+it)|\ll  \e^{-1}\left\lbrace\begin{array}{ll}
		  |t|^{n+2\e+\nu_n}  & n\equiv 0\Mod{2},\\
		|t|^{\frac{3n}{2}+3\e+\nu_n} & n\equiv 1\Mod{2}.	
	\end{array}\right.
	 $$
Using the Phragm\'en-Lindel\"of principle we get that for any $\sigma\in (0,n)$ and any $|t|\geq 2$ we can bound 
$$|F(\sigma+it)|\ll  \e^{-1} |t|^{a(\sigma)+{3}\e},$$
with 
$$a(\sigma):=\left\lbrace\begin{array}{ll}
		 {\nu_n+n-\sigma} & n\equiv 0\Mod{2},\\
	{\nu_n+\frac{3}{2}(n-\sigma)}& n\equiv 1\Mod{2}.\\
	\end{array}\right.
	$$
	In particular, 
	$$a(\tfrac{n}{2})=\left\lbrace\begin{array}{ll}
		 {\nu_n+\frac{n}{2}} & n\equiv 0\Mod{2},\\
	{\nu_n+\frac{3n}{4}}& n\equiv 1\Mod{2}.\\
	\end{array}\right.
	$$
Since $|\zeta(1+it)|\gg \frac{1}{\log(|t|)}$ for $|t|\geq 2$ we {can again apply Stirling's formula to} get the lower bound 
$$|H(\tfrac{n}{2}+it)|\gg  \frac{1}{\log(|t|)} {|t|^{\nu_n}}. 
	$$
Combining these bounds we see that 
$$|E_{Q_n}(\tfrac{n}{2}+it,g)|=\frac{|F(\tfrac{n}{2}+it)|}{|H(\tfrac{n}{2}+it)|}\ll \e^{-1}\log(|t|) {|t|^{a(\frac{n}{2})-\nu_n+3\e}}= \e^{-1}\log(|t|) \left\lbrace\begin{array}{ll}
		 |t|^{\frac{n}{2}+3\e} & n\equiv 0\Mod{2},\\
		|t|^{\frac{3n}{4}+3\e} & n\equiv 1\Mod{2}.\\
	\end{array}\right.$$
Taking $\e=\frac{1}{\log(|t|)}$ we {then get the desired sup-norm bounds on $E_{Q_n}(\tfrac{n}{2}+it,g)$}.
\end{proof}
To deal with the case of $n\equiv 0\Mod{8}$ we give an alternative proof which follows the strategy of Young \cite[Lemma 3.1]{Young2018}.
\begin{Prop}\label{p:supeven}
Let $g\in G$ be fixed. 
For any $n\equiv 0\Mod{4}$ we have for all $|t|\geq 2$
$$|E_{Q_n}(\tfrac{n}{2}+it,g)|\ll_{g,\varepsilon} |t|^{\frac{n}{2}+\varepsilon}. $$

\end{Prop}
\begin{proof}
Since $g\in G$ is fixed we can write $g=u_{\bm{x}}a_yk$ with $y\asymp_g 1$. 
Let $n\equiv 0\Mod{4}$. Expanding $E_{Q_n}(s,g)$ into its Fourier expansion and taking absolute values on all terms we can bound 
$$|E_{Q_n}(\tfrac{n}{2}+it,g)|\leq 2y^{\frac{n}{2}}+\frac{2^{1-\frac{n}{2}}\pi^{\frac{n}{2}}y^{\frac{n}{2}}}{|\G(\tfrac{n}{2}+it)\zeta(\tfrac{n}{2}+it)\zeta(1+it)|}\sum_{a=1}^\infty \sum_{{\bm\lambda}\in \Lambda^*_{\rm pr}} |\epsilon_n(\tfrac{n}{2}+it,a\bm \lambda)K_{it}(2\pi a\|\bm \lambda\|y)|.$$
Here we used the estimate $|\Phi_n(\frac{n}{2}+it)|\leq 1$; see \rmkref{rmk:fueqfail}. Using the Stirling's approximation we can bound $\frac{1}{|\G(\frac{n}{2}+it)|}\ll |t|^{\frac{1-n}{2}}\cosh(\frac{\pi t}{2})$, while for the zeta function we can bound $\frac{1}{|\zeta(1+it)|}\ll |t|^\e$ and $\frac{1}{|\zeta(\tfrac{n}{2}+it)|}\ll 1$. Moreover, using the estimate  
$|\epsilon_n(\tfrac{n}{2}+it; a\bm\lambda)|\ll a^{\frac{n}{2}-1+\e}\|\bm\lambda\|^\e$ {(see \eqref{equ:crilinebd})} we get
$$|E_{Q_n}(\tfrac{n}{2}+it,g)|\ll  y^{\frac{n}{2}}+ y^{\frac{n}{2}}|t|^{\frac{1-n}{2}+\e}\sum_{a=1}^\infty \sum_{\bm\lambda\in \Lambda^*_{\rm pr}} a^{\frac{n}{2}-1+\e}\|\bm\lambda\|^\e \cosh\left(\tfrac{\pi t}{2}\right)|K_{it}(2\pi a\|\bm \lambda\|y)|.$$
{Since the function on the right is even in $t$, we may
assume that  $t > 0$.} We {now} use the following bound for the Bessel function extracted in \cite[equation (3.4)]{Young2018} from the asymptotic expansion of  \cite{Balogh1967}.  There are constants $C,c>0$ such that 
$$\cosh(\tfrac{\pi t}{2})|K_{it}(u)|\ll\left\lbrace\begin{array}{ll} 
t^{-1/4}(t-u)^{-1/4} & 0<u< t-Ct^{1/3},\\
t^{-1/3} & |u-t|\leq Ct^{1/3},\\
u^{-1/4}(u-t)^{-1/4} \exp(-c\frac{u^{3/2}(u-t))^{3/2}}{t^2}) & u>t+Ct^{1/3}.
\end{array}\right.
$$
We split the sum into different regions corresponding to these bounds. First when $2\pi a\|\bm\lambda\| y\leq \tfrac{t}{2}$ we can bound 
$\cosh(\frac{\pi t}{2})|K_{it}(2\pi a\|\bm \lambda\|y)|\ll t^{-1/2}$ to get that the contribution from this range is bounded by 
\begin{align*}
y^{\frac{n}{2}}t^{-\frac{n}{2}+\e}\sum_{a=1}^{\frac{t}{2\pi y}}  a^{\frac{n}{2}-1+\e}\sum_{\substack{\bm{\lambda}\in \Lambda^*_{\rm pr}\\ \|\bm\lambda\|\leq \frac{t}{4\pi ay}}}\|\bm{\lambda}\|^\e&\ll 
y^{\frac{n}{2}}t^{-\frac{n}{2}+\e}\sum_{a=1}^{\frac{t}{2\pi y}}  a^{\frac{n}{2}-1+\e}\left(\tfrac{t}{ay}\right)^{n+\e}\ll t^{\frac{n}{2}+2\e}y^{-\frac{n}{2}-\e}{\ll_{g,\e} t^{\frac{n}{2}+2\e}}.
 \end{align*}
 Here for the first inequality we used that for all $T>1$,
 \begin{align}\label{equ:countinglapt}
 \#\{\bm{\lambda}\in \Lambda^*_{\rm pr}: \|\bm{\lambda}\|\leq T\}=\frac{2}{\zeta(n)}T^n+O(T^{n-1}).
 \end{align}
Next we split the region $\frac{t}{2}<2\pi a\|\bm\lambda\|y<t-Ct^{1/3}$ into dyadic intervals of the form $t-2A<2\pi a\|\bm\lambda\|y\leq t-A$ with  $A$ a power of $2$ in the range $Ct^{1/3}<A<\frac{t}{4}$, and we denote by $S(A)$ the contribution of terms from this dyadic interval. For each such region we can estimate 
$\cosh(\frac{\pi t}{2})|K_{it}(2\pi a\|\bm \lambda\|y)|\ll t^{-1/4}A^{-1/4}$ so the contribution of these terms is bounded by 
\begin{align*}
S(A)&\ll y^{\frac{n}{2}}A^{-1/4}t^{\frac{1}{4}-\frac{n}{2}+\e}\sum_{a=1}^{\frac{t-A}{2\pi y}} \mathop{\sum_{\bm\lambda\in \Lambda^*_{\rm pr}}}_{\tfrac{t-2A}{2\pi ay}<\|\bm\lambda\|\leq \tfrac{t-A}{2\pi ay}} a^{\frac{n}{2}-1+\e}\|\bm\lambda\|^\e.
\end{align*}
Using again the estimate \eqref{equ:countinglapt}
we can bound
$$\#\left\{\bm\lambda\in \Lambda^*_{\rm pr}: \tfrac{t-2A}{2\pi ay}< \|\bm\lambda\|\leq \tfrac{t-A}{2\pi ay}\right\}\ll (\tfrac{t}{ay})^{n-1}(\tfrac{A}{ay}+1),$$
implying that 
\begin{align*}
S(A) &\ll y^{\frac{n}{2}}A^{-1/4}t^{\frac{1}{4}-\frac{n}{2}+\e}\sum_{a=1}^{\frac{t-A}{2\pi y}} a^{\frac{n}{2}-1+\e}  \left(\tfrac{t}{ay}\right)^{n-1+\e}\left(\tfrac{A}{ay}+1\right) 
\ll_{\e} A^{\frac{3}{4}}t^{-\frac{3}{4}+\frac{n}{2}+\e}.
\end{align*}
Summing over $Ct^{1/3}<A<\frac{t}{4}$ in dyadic intervals, we see that the contribution of these terms is also bounded by $O_{g,\e}(t^{\frac{n}{2}+\e})$. 
A similar argument shows that the contribution of $2\pi a\|\bm\lambda\|y-t|\leq  Ct^{1/3}$ is even smaller and can be bounded by $O(t^{\frac{n}{2}})$ 
as well as the contribution of $2\pi a\|\bm\lambda\|y\geq t+Ct^{1/3}$ where we again partition into dyadic intervals and then use the exponential decay of the Bessel function.
\end{proof}
\begin{rem}
A similar argument could also work when $n\equiv 2\Mod{4}$, however, when $n$ is odd this direct approach does not work due to the following  two main obstacles: First, the factors $\epsilon_n(s;\bm\lambda)$ are no longer expressed in terms of divisor function and we do not have the bound \eqref{equ:crilinebd}. The second obstacle is that when $n$ is odd, the Fourier coefficients $\Phi_n(\tfrac{n}{2}+it;\bm\lambda)$ also contain a factor of 
$L(\tfrac{1}{2}+it,\chi_D)$ with $D=(-1)^{\frac{n-1}{2}}\|2\bm\lambda\|^2$. Hence, in order for the same argument to work we need a bound of the form $L(\tfrac{1}{2}+it,\chi_D)\ll |t|^{\e}$. While such a bound follows from GRH (or from the Lindel\"{o}f hypothesis) it is not known unconditionally. It might still be possible to make this approach work by using such bounds on average, and we leave this as an open problem.
\end{rem}
Combining these results we get the following.
\begin{proof}[{Proof of \thmref{t:supnorm}}]
The upper bound $\nu_n\leq \tfrac{n}{2}$ when $n\equiv 0\pmod{8}$ follows from \propref{p:supeven}  while the bound $\nu_n\leq \tfrac{n}{2}$ for even $n\not\equiv 0\pmod{8}$ and the bound $\nu_n\leq \tfrac{3n}{4}$ for odd $n\geq 3$ follows from \propref{p:supodd}. Next, since we can bound the {$L^2$-norm by the supremum norm} we have that  $\tilde\nu_n\leq \nu_n$. The lower bound $\tilde\nu_n\geq \tfrac{n}{2}-1$ follows from the argument in \rmkref{r:lowerbound}. \end{proof}

\section{Calculation of Fourier coefficients}\label{sec:FourierCalc}
{In this section we prove our Fourier expansion formulas stated in the previous section.}
\subsection{A preliminary identity}
We first prove a preliminary formula for these Fourier coefficients which reduces the problem to evaluating certain counting functions and exponential sums which we compute later.
\begin{Prop}\label{prop:preid}
Keep the notation and assumptions as in \thmref{thm:fouexpan}. Then we have
\begin{align*}
a_{Q_{n,d}}(s,g;\bm{0})=y^s+2^{s-n}d^{s-n}\pi^{\frac{n}{2}} \frac{\G\left(\frac{2s-n}{2}\right)}{\zeta(s)\G(s)}Z_{n,d}(s;\bm{0})y^{n-s}.
\end{align*}
For any nonzero $\bm{\lambda}\in \Lambda^{\ast}$ we have
\begin{align*}
a_{Q_{n,d}}(s,g;\bm{\lambda})=
\frac{2^{s-n+1}\pi^s\|\bm{\lambda}\|^{s-\frac{n}{2}} }{d^{\frac{n}{2}}} \frac{Z_{n,d}(s;\bm{\lambda})}{\zeta(s)\G(s)} y^{n/2}K_{s-n/2}(\tfrac{2\pi\|\bm{\lambda}\|y}{d}),
\end{align*}
where for any $\bm{\lambda}\in \Lambda^{\ast}$ and $ \Re(s)>n$, 
\begin{align}\label{equ:arifunc}
Z_{n,d}(s;\bm{\lambda}):=\sum_{t=1}^{\infty}\frac{\phi_{n,d}(t;\bm{\lambda})}{t^s},
\end{align}
with 
$$\phi_{n,d}(t;\bm{\lambda}):=\sum_{\substack{\bm{h}\in (\Z/2td\Z)^n\\ \|\bm{h}\|^2\equiv -t^2\Mod{2td}}}e\left(\tfrac{2\bm{\lambda}\cdot \bm{h}}{2td}\right).$$
\end{Prop}

\begin{proof}
Abbreviate $Q_{n,d}$ by $Q$. Fix $\cF\subseteq \R^n$ a fundamental domain for $\Lambda$ and note that $\vol(\cF)=2$. Since every vector in $\cV^+_Q(\Z)$ can be written uniquely as a positive multiple of a vector in $\cV^+_Q(\Z)_{\rm pr}$, we can rewrite that for any $\Re(s)>n$ and for any $\bm{\lambda}\in \Lambda^{\ast}$
\begin{align*}
\zeta(s)a_Q(s,g;\bm{\lambda})&=\frac{1}{\vol(\cF)}\int_{\cF}\sum_{\bm{v}\in \cV^+_Q(\Z)}(1+d^2)^{s/2}\|\bm{v}u_{\bm{x}}g\|^{-s}e(-\bm{\lambda}\cdot\bm{x})\,\d\bm{x}\\
&=2^{-2}(1+d^2)^{s/2}\int_{\cF}\sum_{\bm{v}\in \cV_Q(\Z)}\|\bm{v}u_{\bm{x}}g\|^{-s}e(-\bm{\lambda}\cdot\bm{x})\,\d\bm{x},
\end{align*}
where 
$$
\cV_Q(\Z):=\left\{\bm{v}\in \Z^{n+2}\setminus\{\bm{0}\}: Q(\bm{v})=0\right\}
$$ 
is the set of integral points in the two-sheeted light cone of $Q$.
Next, for $\bm{v}\in\R^{n+2}$ let us write 
$\bm{v}=(v_1,\bm{w},v_2)\in \R\times \R^n\times \R$ so that
$$
Q(\bm{v})=\|\bm{w}\|^2+v_1^2-d^2v_2^2=\|\bm{w}\|^2+(v_1+dv_2)(v_1-dv_2).
$$
Note that if $v_1+dv_2=0$, 
then $Q(\bm{v})=0 \Leftrightarrow \|\bm{w}\|^2=0$. Thus for $\bm{v}\in \cV_Q(\Z)$, $v_1+dv_2=0$ implies that $\bm{w}=\bm{0}$. Moreover, by direct computation,
\begin{align*}
(-dv_2, \bm{0}, v_2)u_{\bm{x}}a_y=y^{-1}v_2(-d, \bm{0}, 1).
\end{align*}
Thus the contribution of the terms with $v_1+dv_2=0$ is
\begin{align*}
2^{-2}(1+d^2)^{s/2}\sum_{v_2\neq 0}\|y^{-1}(-dv_2, \bm{0}, v_2)\|^{-s}\int_{\cF}e(-\bm{\lambda}\cdot\bm{x})\,\d\bm{x}=\left\{\begin{array}{ll} \zeta(s)y^s  & \bm{\lambda}=\bm{0},\\
0 & \bm{\lambda}\neq \bm{0}.
\end{array}\right.
\end{align*}
If $v_1+dv_2\neq 0$ then $Q(\bm{v})=0$ implies $v_1-dv_2=-\tfrac{\|\bm{w}\|^2}{v_1+dv_2}$.
Using this relation we can compute
\begin{align*}
(v_1,\bm{w},v_2)u_{\bm{x}}a_y=
(v_1+dv_2)\left(\tfrac{y}{2}-\tfrac{1}{2y}\| \bm{x}d+\tfrac{\bm{w}}{v_1+dv_2}\|^2, \bm{x}d+\tfrac{\bm{w}}{v_1+dv_2},  \tfrac{y}{2d}+\tfrac{1}{2dy}\| \bm{x}d+\tfrac{\bm{w}}{v_1+dv_2}\|^2\right),
\end{align*}
and taking the norm square we notice that 
\begin{align*}
\|(v_1,\bm{w},v_2)u_{\bm{x}}a_y\|^2=
\tfrac{(v_1+dv_2)^2(1+\tfrac{1}{d^2})}{4y^2}\left(y^2+\|\bm{x}d+\tfrac{\bm{w}}{v_1+dv_2}\|^2\right)^2.
\end{align*}

Let us denote $(t_1,t_2):=(v_1+dv_2,v_1-dv_2)$ so that $Q(\bm{v})=\|\bm{w}\|^2+t_1t_2$. We note that the pairs $(v_1,v_2)\in \Z^2$ are in one-to-one correspondence with the pairs $(t_1,t_2)\in \Z^2$ satisfying $t_1\equiv t_2\Mod{2d}$. Moreover, the points $(v_1, \bm{w},v_2)\in \cV_Q(\Z)$ with $v_1+dv_2\neq 0$ are in one-to-one correspondence with the points $(\bm{w}, t_1)\in \Z^{n}\times (\Z\setminus\{0\})$ satisfying $\|\bm{w}\|^2\equiv -t_1^2\Mod{2t_1d}$ (with $t_2$ uniquely determined by $(\bm{w}, t_1)$ via $t_2=-\frac{\|\bm{w}\|^2}{t_1}$). Thus the contribution of the remaining terms with $v_1+dv_2\neq 0$ is

\begin{align*}
2^{-2}(2dy)^s\sum_{t_1\neq 0}|t_1|^{-s}
\int_{\cF}\sum_{\|\bm{w}\|^2\equiv -t_1^2\Mod{2t_1d}}
\left(y^2+\|\bm{x}d+\tfrac{\bm{w}}{t_1}\|^2 \right)^{-s}
e(-\bm{\lambda}\cdot\bm{x})\,\d\bm{x}.
\end{align*}
We can further split the above inner sum over congruence classes modulo $t_1d\Lambda$, noting that the congruence condition $\|\bm{w}\|^2\equiv -t_1^2\Mod{2t_1d}$ is preserved modulo $t_1d\Lambda$. 

Thus writing elements of a fixed congruence class modulo $t_1d\Lambda$ as $\bm{w}=\bm{h}+t_1d\bm{v}$ with $\bm{v}\in\Lambda$ for some fixed $\bm{h}$ in this congruence class 
we get the contribution of each such a congruence class to the above integral is 

\begin{align*}
\cI_{\bm{h}}(\bm\lambda)&:=
\int_{\cF}\sum_{\bm{v}\in \Lambda}
\left(y^2+\|(\bm{x}+\bm{v})d+\tfrac{\bm{h}}{t_1}\|^2 \right)^{-s}
e(-\bm{\lambda}\cdot\bm{x})\,\d\bm{x}\\
&=\int_{\R^n}
\left(y^2+d^2\|\bm{x}+\tfrac{\bm{h}}{t_1d}\|^2 \right)^{-s}
e(-\bm{\lambda}\cdot\bm{x})\,\d\bm{x}\\
&=e\left(\tfrac{\bm{\lambda}\cdot \bm{h}}{t_1 d}\right)\int_{\R^n}
\left(y^2+\|\d\bm{x}\|^2 \right)^{-s}
e(-\bm{\lambda}\cdot\bm{x})\,\d\bm{x}\\
&=d^{-n}y^{n-2s}e\left(\tfrac{\bm{\lambda}\cdot \bm{h}}{t_1 d}\right)\int_{\R^n}
\frac{e(-\frac{\bm{\lambda}\cdot\bm{x} y}{d})\,\d\bm{x}}{\left(1+\|\bm{x}\|^2 \right)^{s}}.
\\
\end{align*}
When $\bm{\lambda}= \bm{0}$ we get
\begin{align*}
\cI_{\bm{h}}(\bm{0})=d^{-n}y^{n-2s}\int_{\R^n}\frac{\d\bm{x}}{(1+\|\bm{x}\|^2)^{s}}=d^{-n}y^{n-2s}\pi^{n/2}\frac{\G\left(\frac{2s-n}{2}\right)}{\G(s)}.
\end{align*}
Thus in this case, the contribution of the terms with $v_1+dv_2\neq 0$ is
\begin{align*}
&2^{-2}(2dy)^s \sum_{t_1\neq 0}|t_1|^{-s} \sum_{\substack{\bm{h}\in \Z^n/t_1d\Lambda\\ \|\bm{h}\|^2\equiv -t_1^2\Mod{2t_1d}}}d^{-n}y^{n-2s}\pi^{n/2}\frac{\G\left(\frac{2s-n}{2}\right)}{\G(s)}\\
&=2^{s-1}d^{s-n}y^{n-s}\pi^{n/2}  \frac{\G\left(\frac{2s-n}{2}\right)}{\G(s)}
\sum_{t_1= 1}^\infty \frac{\tilde{\phi}_{n,d}(t_1;\bm{0})}{t_1^{s}},
\end{align*}
where 
$$\tilde{\phi}_{n,d}(t;\bm{0}):=\#\{\bm{h}\in \Z^n/t d\Lambda:  \|\bm{h}\|^2\equiv -t^2\Mod{2td}\}.$$
When $\bm{\lambda}\neq \bm{0}$, let $k\in \SO_n(\R)$ be the rotation such that $\bm{\lambda}k=\|\bm{\lambda}\|(1,0,\ldots,0)$. Then applying the change of variables $\bm{x}\mapsto \bm{x}k^t$ and using the identity $\|\bm{x}k^t\|=\|\bm{x}\|$ we get for any $t>0$,
\begin{align*}
\int_{\R^n}\frac{e(-t\bm{\lambda}\cdot\bm{x})\,\d\bm{x}}{(1+\|\bm{x}\|^2)^{s}}&=\int_{\R^n}\frac{e(-t\bm{\lambda}k\cdot\bm{x})\,\d\bm{x}}{(1+\|\bm{x}\|^2)^{s}}=\int_{\R^n}\frac{e(-t\|\bm{\lambda}\|x_1)\,\d\bm{x}}{(1+\|\bm{x}\|^2)^{s}}\\
&=\int_{\R}\frac{e(-t\|\bm{\lambda}\|x)}{(1+x^2)^s}\int_{\R^{n-1}}\frac{\d\bm{t}}{\left(1+\frac{\|\bm{t}\|^2}{1+x^2}\right)^s}\,\d x\\
&=\int_{\R}\frac{e(-t\|\bm{\lambda}\|x)\,\d x}{(1+x^2)^{s-(n-1)/2}}\int_{\R^{n-1}}\frac{\d\bm{t}}{(1+\|\bm{t}\|^2)^s}\\
&=\frac{2\pi^{s-\frac{n-1}{2}}}{\G\left(s-\frac{n-1}{2}\right)}(t\|\bm{\lambda}\|)^{s-n/2}K_{s-n/2}(2\pi\|\bm{\lambda}\|t)\times \pi^{(n-1)/2}\frac{\G\left(s-\frac{n-1}{2}\right)}{\G(s)}\\
&=\frac{2\pi^s(\|\bm{\lambda}\|t)^{s-n/2}K_{s-n/2}(2\pi\|\bm{\lambda}\|t)}{\G(s)},
\end{align*}
implying that
\begin{align*}
\cI_{\bm{h}}(\bm{\lambda})=d^{-n}y^{n-2s}e\left(\tfrac{\bm{\lambda}\cdot \bm{h}}{t_1d}\right) \frac{2\pi^s(\|\bm{\lambda}\|yd^{-1})^{s-n/2}K_{s-n/2}(2\pi\|\bm{\lambda}\|yd^{-1})}{\G(s)}.
\end{align*}
Hence in this case, the contribution of the terms with $v_1+dv_2\neq 0$ is
\begin{align*}
&2^{s}d^{-\frac{n}{2}}y^{\frac{n}{2}} \frac{\pi^s\|\bm{\lambda}\|^{s-\frac{n}{2}}K_{s-n/2}(2\pi\|\bm{\lambda}\|yd^{-1})}{\G(s)}\sum_{t_1=1}^\infty\frac{ \tilde{\phi}_{n,d}(t_1;\bm{\lambda})}{t_1^{s} },
\\
\end{align*}
where 
$$\tilde{\phi}_{n,d}(t;\bm{\lambda}):=\sum_{\substack{\bm{h}\in \Z^n/td\Lambda\\ \|\bm{h}\|^2\equiv -t^2\Mod{2td}}}e\left(\tfrac{\bm{\lambda}\cdot \bm{h}}{td}\right).$$
Finally noting that $2td\Z^n\leq td\Lambda$ is a subgroup of index $2^{n-1}$ we have
$$\tilde{\phi}_{n,d}(t;\bm{\lambda})=2^{1-n}\sum_{\substack{\bm{h}\in (\Z/2td\Z)^n\\ \|\bm{h}\|^2\equiv -t^2\Mod{2td}}}e\left(\tfrac{2\bm{\lambda}\cdot \bm{h}}{2td}\right)=2^{1-n}\phi_{n,d}(t; \bm{\lambda}).$$
Plugging this relation into the previous expression finishes the proof.
\end{proof}

\subsection{Euler product expansion of $Z_{n,d}(s; \bm\lambda)$}
For the remaining of this section we assume $d\geq 1$ is odd and square-free. Fix $\bm{\lambda}\in \Lambda^*$. In order to decompose the function $Z_{n,d}(s;\bm\lambda)$ as an Euler product we need to express the function 
$\phi_{n,d}(t;\bm{\lambda})$ in terms of a multiplicative function (of $t$). For this we define two associated (almost) multiplicative functions:
For $\bm{m}\in \Z^n$ and $t\in \N$ define 
\begin{equation}\label{equ:phinfun1}
\phi_n(t;\bm{m}):=\sum_{\substack{\bm{h}\in (\Z/t\Z)^n\\ \|\bm{h}\|^2\equiv 0\Mod{t}}}e\left(\tfrac{\bm{m}\cdot \bm{h}}{t}\right),\end{equation}
and 
\begin{equation}\label{equ:varphifun2}
\vf_n(t;\bm{m}):=\sum_{\substack{\bm{h}\in (\Z/t\Z)^n\\ \|\bm{h}\|^2\equiv -1\Mod{t}}}e\left(\tfrac{\bm{m}\cdot \bm{h}}{t}\right).\end{equation}
Note that $\phi_n(t;\bm{m})$ is multiplicative in $t$, but $\varphi_n(t;\bm{m})$ is not, instead, satisfying the relation that for any co-prime pair $(t_1,t_2)\in \N^2$,
\begin{align}\label{equ:almultivp}
\varphi_n(t_1t_2;\bm{m})=\varphi_n(t_1;\bar{t}_2\bm{m})\varphi_n(t_2;\bar{t}_1\bm{m}),
\end{align}
where $\bar{t}_1$ is the multiplicative inverse of $t_1$ modulo $t_2$ and similarly $t_2\bar{t}_2\equiv 1\Mod{t_1}$. For any fixed $t\in\N$, $\phi_n(t;\cdot)$ and $\varphi_n(t;\cdot)$ can be viewed as functions in $\bm{m}\in (\Z/t\Z)^n$. For later reference, we also note that it is easy to see that for any $\bm{\lambda}\in \Lambda^*$,
\begin{align}\label{equ:evaat2}
\phi_n(2;2\bm{\lambda})=2^{n-1}\quad \text{and}\quad \varphi_n(2;2\bm{\lambda})=(-1)^{2\lambda_1}2^{n-1}.
\end{align}
\begin{Lem}\label{lem:firsdeclem}
For any $t\in \N$ 
we have 
$$\phi_{n,d}(t; \bm\lambda)=\phi_n(at;2\bm{\lambda})\vf_n(\tfrac{2d}{a}; 2t\bm{\lambda}),$$
where $a=\gcd(t,2d)$. 
\end{Lem}
\begin{proof}
Write $t=t_1a$ and $2d=d_1a$ with $\gcd(t_1,d_1)=1$. Since we assume $d$ is odd and square-free we also have $\gcd(t_1a^2,d_1)=1$ so writing $2td=a^2t_1d_1$ and using the Chinese remainder theorem we get
\begin{align*}
\phi_{n,d}(t;\bm{\lambda})&=\sum_{\substack{\bm{h}\in (\Z/2td\Z)^n\\ \|\bm{h}\|^2\equiv -t^2\Mod{2td}}}e\left(\tfrac{2\bm{\lambda}\cdot \bm{h}}{2td}\right)\\
&= \sum_{\substack{\bm{h}\in (\Z/a^2t_1\Z)^n\\ \|\bm{h}\|^2\equiv 0\Mod{a^2t_1}}}e\left(\tfrac{2\bm{\lambda}\cdot \bm{h}}{a^2t_1}\right)\sum_{\substack{\bm{h}\in (\Z/d_1\Z)^n\\ \|\bm{h}\|^2\equiv -t^2\Mod{d_1}}}e\left(\tfrac{2\bm{\lambda}\cdot \bm{h}}{d_1}\right).\\
\end{align*}
{Noting that $a^2t_1=at$ we see that the first sum is 
$$ \sum_{\substack{\bm{h}\in (\Z/a^2t_1\Z)^n\\ \|\bm{h}\|^2\equiv 0\Mod{a^2t_1}}}e\left(\tfrac{2\bm{\lambda}\cdot \bm{h}}{a^2t_1}\right)=\sum_{\substack{\bm{h}\in (\Z/at\Z)^n\\ \|\bm{h}\|^2\equiv 0\Mod{a t}}}e\left(\tfrac{2\bm{\lambda}\cdot \bm{h}}{a t}\right)=\phi_n(at; 2\bm\lambda).$$
For the second sum, note that $\gcd(t,d_1)=1$ and make a change of variables $\bm{h}\mapsto t\bm{h}$ to get that 
$$\sum_{\substack{\bm{h}\in (\Z/d_1\Z)^n\\ \|\bm{h}\|^2\equiv -t^2\Mod{d_1}}}e\left(\tfrac{2\bm{\lambda}\cdot \bm{h}}{d_1}\right)=\sum_{\substack{\bm{h}\in (\Z/d_1\Z)^n\\ \|\bm{h}\|^2\equiv -1\Mod{d_1}}}e\left(\tfrac{2t\bm{\lambda}\cdot \bm{h}}{d_1}\right)=\vf_n(d_1; 2t\bm\lambda),$$
implying that $\phi_{n,d}(t; \bm\lambda)=\phi_n(at;2\bm{\lambda})\vf_n(\tfrac{2d}{a}; 2t\bm{\lambda})$ as claimed.}

\end{proof}
In particular, when $\bm\lambda=\bm{0}$ we get that $\phi_{n,d}(t;\bm{0})=\phi_n(at;\bm0)\vf_n(\tfrac{2d}{a};\bm0)$ is a product of multiplicative functions leading to the following Euler product formula.
\begin{Prop}\label{prop:eulprodconst}
For any $\Re(s)>n$,
$$Z_{n,d}(s;\bm{0})=\prod_{p\mid 2d}(\vf_n(p; \bm 0)+ \tilde{Z}_n^{(p)}(s;\bm{0}))\prod_{(p,2d)=1} Z_n^{(p)}(s;\bm{0}),$$
where 
$$Z_n^{(p)}(s;\bm{0}):=\sum_{k=0}^\infty \frac{\phi_n(p^k;\bm{0})}{p^{ks}}\quad\text{and}\quad \tilde{Z}_n^{(p)}(s;\bm{0}):=\sum_{k=0}^\infty \frac{\phi_n(p^{k+2};\bm{0})}{p^{(k+1)s}}.
$$
\end{Prop}
\begin{proof}
Separating the sum over $t$ based on the $\gcd(t,2d)$ we get 
\begin{align*}
Z_{n,d}(s; \bm0)&= \sum_{t=1}^\infty \frac{\phi_{n,d}(t;\bm0)}{t^s}\\
&= \sum_{a\mid 2d}  \sum_{\gcd(t,2d)=a}  \frac{\phi_n(at; \bm 0) \vf_n(\tfrac{2d}{a}; \bm 0) }{t^s}\\
&= \sum_{a\mid 2d}  \vf_n(\tfrac{2d}{a}; \bm 0)  \sum_{{\substack{t\in \N\\ \gcd(t,\tfrac{2d}{a})=1}}}  \frac{\phi_n(a^2t;\bm 0)}{(at)^s}.
\end{align*}
The second sum can be factored into an Euler product
\begin{align*}
\sum_{{\substack{t\in \N\\ \gcd(t,\tfrac{2d}{a})=1}}}  \frac{\phi_n(a^2t; \bm 0)}{(at)^s}&=\prod_{p\mid a}\left(\sum_{k=0}^\infty \frac{\phi_n(p^{k+2};\bm0)}{p^{(k+1)s}}\right)\prod_{(p,2d)=1}\left(\sum_{k=0}^\infty \frac{\phi_n(p^k;\bm0)}{p^{ks}}\right)\\
&= \prod_{p\mid a}\tilde{Z}_n^{(p)}(s;\bm0) \prod_{(p,2d)=1}Z_n^{(p)}(s;\bm{0}),
\end{align*}
and taking $ \prod_{(p,2d)=1}Z^{(p)}(s;\bm0)$ as a common factor, using that $\vf_n(\cdot;\bm{0})$ is multiplicative we get 
\begin{displaymath}
Z_{n,d}(s;\bm{0})=   \prod_{(p,2d)=1}Z_n^{(p)}(s;\bm0)\prod_{p\mid 2d}(\vf_n(p; \bm 0)+ \tilde{Z}_n^{(p)}(s;\bm0)).\qedhere 
\end{displaymath}
\end{proof}

When $\bm\lambda\neq \bm{0}$ this decomposition is more complicated and is described in the following.
\begin{Prop}\label{prop:nonconsprform}
For any $\bm\lambda\in \Lambda^*\setminus\{\bm{0}\}$ and for any $\Re(s)>n$, 
\begin{align*}
Z_{n,d}(s;\bm{\lambda})&= \sum_{d_1|d} \frac{1}{\vf(d_1)}\sum_{\chi\Mod{d_1}} \bigg(\tilde{Z}_n^{(2)}(s;\chi, \bm\lambda)+ (-1)^{2\lambda_1}2^{n-1}\bigg)\\
&\overline\chi(\tfrac{2d}{d_1}) 
\cS(d_1,\chi, 2\bm\lambda)\prod_{p|\frac{d}{d_1}}\tilde{Z}_n^{(p)}(s;\chi, \bm\lambda)\prod_{(p,2d)=1}Z^{(p)}_n(s; \chi, \bm\lambda),
 \end{align*}
where {$\varphi$ is the Euler's totient function} and 
for any Dirichlet character $\chi$ of modulus $d_1\in \N$, for any prime $p$ and for any $\Re(s)>n$,
\begin{equation}\label{eq:Znp}
Z_n^{(p)}(s; \chi, \bm{\lambda}):=\sum_{k=0}^\infty \frac{\chi(p^k)\phi_n(p^k; 2\bm \lambda)}{p^{ks}},\quad \tilde{Z}_n^{(p)}(s;\chi, \bm\lambda):=\sum_{k=0}^\infty \frac{\chi(p^{k+2})\phi_n(p^{k+2}; 2\bm \lambda)}{p^{(k+1)s}},
\end{equation}
and for any $\bm{m}\in \Z^n$,
\begin{equation}\label{eq:Sdchi}
\cS(d_1,\chi, \bm{m}):=\sum_{t\in (\Z/d_1\Z)^{\times}} \vf_n(d_1; t\bm{m}) \overline{\chi}(t).
\end{equation}
\end{Prop}
\begin{remark}\label{rmk:d1case}
 To simplify notation when $\chi\equiv 1$ we abbreviate $Z_{n}^{(p)}(s;\chi,\bm{\lambda})$ and $\tilde{Z}_{n}^{(p)}(s;\chi,\bm{\lambda})$ by $Z_n^{(p)}(s;\bm{\lambda})$ and $\tilde{Z}_{n}^{(p)}(s;\bm{\lambda})$ respectively. Note that when $d=1$ the above formula reads as
\begin{align*}
Z_{n,1}(s;\bm{\lambda})&=\left(\tilde{Z}_n^{(2)}(s; \bm\lambda)+ (-1)^{2\lambda_1}2^{n-1}\right)\text{$\prod$}'_{p}Z_n^{(p)}(s;\bm{\lambda}),
\end{align*}
where the product $\prod'_p$ is over all odd prime numbers.
\end{remark}
\begin{proof}[Proof of \propref{prop:nonconsprform}]
By \lemref{lem:firsdeclem} we have
\begin{align*}
Z_{n,d}(s;\bm{\lambda})&= \sum_{a\mid 2d}\sum_{\substack{t\in\N\\ \gcd(t,2d)=a}}\frac{\phi_n(at; 2\bm{\lambda})\vf_n(\frac{2d}{a}; 2t\bm{\lambda})}{t^s}\\
&= \sum_{a\mid 2d}\sum_{\substack{t\in\N\\ \gcd(t,2d/a)=1}}\frac{\phi_n(a^2t; 2\bm{\lambda})\vf_n(\tfrac{2d}{a}; 2at\bm{\lambda})}{(at)^s}\\
&= \sum_{a\mid 2d}\sum_{b\in (\Z/\tfrac{2d}{a}\Z)^{\times}}\vf_n(\tfrac{2d}{a}; 2ab\bm{\lambda})\mathop{\sum_{t\in \N}}_{t\equiv b\Mod{2d/a}}\frac{\phi_n(a^2t; 2\bm{\lambda})}{(at)^s}\\
&{= \sum_{a\mid 2d}\sum_{b\in (\Z/\tfrac{2d}{a}\Z)^{\times}}\vf_n(\tfrac{2d}{a}; 2b\bm{\lambda})\mathop{\sum_{t\in \N}}_{t\equiv \bar{a}b\Mod{2d/a}}\frac{\phi_n(a^2t; 2\bm{\lambda})}{(at)^s}},
\end{align*}
{where in the last line we made a change of variable $ab\mapsto b$, and $\bar{a}$ is the inverse of $a$ in $(\Z/\frac{2d}{a}\Z)^{\times}$.}
Now for the inner sum we can capture the condition that $t\equiv {\bar{a}}b\Mod{2d/a}$ by summing over Dirichlet characters modulo $2d/a$ to get 
\begin{align*}
\mathop{\sum_{t\in \N}}_{t\equiv {\bar{a}}b\Mod {2d/a}}\frac{\phi_n(a^2t; 2\bm{\lambda})}{(at)^s}&=\frac{1}{\vf(\tfrac{2d}{a})}\sum_{\chi\Mod{2d/a}} \overline\chi({\bar{a}}b)\sum_{t=1}^\infty \frac{\chi(t)\phi_n(a^2t; 2\bm{\lambda})}{(at)^s}\\
&=\frac{1}{\vf(\tfrac{2d}{a})}\sum_{\chi\Mod{2d/a}} \overline{\chi}(ba)\sum_{t=1}^\infty \frac{\chi(a^2t)\phi_n(a^2t;2\bm{\lambda})}{(at)^s}.
\end{align*}
Since the function $\chi(t)\phi_n(t; 2\bm{\lambda})$ is multiplicative in $t$ the innermost sum factors as  
$$\sum_{t=1}^\infty \frac{\chi(a^2t)\phi_n(a^2t; 2\bm{\lambda})}{(at)^s}=\prod_{p\mid a}\tilde{Z}_n^{(p)}(s; \chi, \bm\lambda)\prod_{(p,2d)=1}Z^{(p)}_n(s; \chi,\bm\lambda),$$
where we used that $Z_n^{(p)}(s;\chi,\bm{\lambda})=1$ for any prime $p\mid \frac{2d}{a}$ and for any Dirichlet character $\chi$ of modulus $\frac{2d}{a}$.
Plugging this back we get that 
\begin{align*}
\mathop{\sum_{t\in \N}}_{t\equiv {\bar{a}}b\Mod {2d/a}}\frac{\phi_n(a^2t; 2\bm{\lambda})}{(at)^s}&= \frac{1}{\vf(\tfrac{2d}{a})}\sum_{\chi\Mod{2d/a}} \overline\chi(ba) \prod_{p|a}\tilde{Z}_n^{(p)}(s; \chi, \bm\lambda)\prod_{(p,2d)=1}Z^{(p)}_n(s; \chi,\bm\lambda),
\end{align*}
and hence
\begin{align*}
Z_{n,d}(s;\bm{\lambda})
&= \sum_{a\mid 2d}\sum_{b\in (\Z/\tfrac{2d}{a}\Z)^{\times}} \vf_n(\tfrac{2d}{a};2b\bm{\lambda})\frac{1}{\vf(\tfrac{2d}{a})}\sum_{\chi\Mod{2d/a}} \overline\chi(ba) \prod_{p|a}\tilde{Z}_n^{(p)}(s; \chi, \bm\lambda)\prod_{(p,2d)=1}Z^{(p)}_n(s; \chi,\bm\lambda)\\
&= \sum_{d_1\mid 2d} \frac{1}{\vf(d_1)}\sum_{\chi\Mod{d_1}} \overline\chi(\tfrac{2d}{d_1})\sum_{b\in (\Z/d_1\Z)^{\times}} \vf_n(d_1;2b\bm{\lambda}) \overline{\chi}(b)\prod_{p\mid \tfrac{2d}{d_1}}\tilde{Z}_n^{(p)}(s; \chi, \bm\lambda)\prod_{(p,2d)=1}Z^{(p)}_n(s; \chi,\bm\lambda)\\
&= \sum_{d_1\mid 2d} \frac{1}{\vf(d_1)}\sum_{\chi\Mod{d_1}} \overline\chi(\tfrac{2d}{d_1}) \cS(d_1,\chi, 2\bm\lambda)\prod_{p\mid \tfrac{2d}{d_1}}\tilde{Z}_n^{(p)}(s; \chi, \bm\lambda)\prod_{(p,2d)=1}Z^{(p)}_n(s; \chi,\bm\lambda),
 \end{align*}
 where $\cS(d_1,\chi,2\bm{\lambda})$ is as given in \eqref{eq:Sdchi}.
We now group together corresponding even and odd divisors of $2d$. {When the divisor is even, we rewrite it as $2d_1$ so that $d_1\mid d$.} 
We can 
identify the set of Dirichlet characters modulo $d_1$ with the set of Dirichlet characters modulo $2d_1$ as follows: Any $\chi\Mod{d_1}$ gets identified with $\chi':=\chi\tilde{\chi}\Mod{2d_1}$, where $\tilde{\chi}$ is the unique Dirichlet character modulo $2$. In particular, $\chi$ and $\chi'$ agree at odd integers and hence {$\tilde{Z}_n^{(p)}(s; \chi', \bm\lambda)=\tilde{Z}_n^{(p)}(s; \chi, \bm\lambda)$ and $Z^{(p)}_n(s; \chi',\bm\lambda)=Z^{(p)}_n(s; \chi,\bm\lambda)$ for any odd prime $p$.} With this identification, using that for any odd $b$, 
$$\vf_n(2d_1; 2b\bm{\lambda})=\vf_n(d_1; \overline{2}b2\bm{\lambda})\vf_n(2; 2\bm{\lambda})=(-1)^{2\lambda_1}2^{n-1}\vf_n(d_1; \overline{2}b2\bm{\lambda})
$$ 
with $\overline{2}$ the multiplicative inverse of $2$ modulo $d_1$ (see \eqref{equ:almultivp} and \eqref{equ:evaat2}) we get for any $\chi'\Mod{2d_1}$ (identified with $\chi\Mod{d_1}$)
\begin{align*}
\cS(2d_1,\chi', 2\bm\lambda)&=\sum_{b\in (\Z/2d_1\Z)^{\times}} \vf_n(2d_1; 2b\bm{\lambda}) \overline{\chi'}(b)\\
&= \vf_n(2; 2\bm{\lambda})\sum_{b\in (\Z/d_1\Z)^{\times}} \vf_n(d_1; \overline{2}b2\bm{\lambda}) \overline{\chi}(b)\\
&= (-1)^{2\lambda_1}2^{n-1}\overline{\chi}(2)\sum_{b\in (\Z/d_1\Z)^{\times}} \vf_n(d_1; 2b\bm{\lambda}) \overline{\chi}(b)\\
&=(-1)^{2\lambda_1}2^{n-1}\overline{\chi}(2)\cS(d_1,\chi,2\bm{\lambda}).
\end{align*}
Using this and collecting together contribution from even and odd divisors of $2d$ we get 
\begin{align*}
Z_{n,d}(s;\bm{\lambda})&= \sum_{d_1|d} \frac{1}{\vf(d_1)}\sum_{\chi\Mod{d_1}} \left(\tilde{Z}_n^{(2)}(s;\chi, \bm\lambda)+ (-1)^{2\lambda_1}2^{n-1}\right)\\
&\overline\chi(\tfrac{2d}{d_1}) 
\cS(d_1,\chi, 2\bm\lambda)\prod_{p|\frac{d}{d_1}}\tilde{Z}_n^{(p)}(s; \chi, \bm\lambda)\prod_{(p,2d)=1}Z^{(p)}_n(s; \chi,\bm\lambda),
 \end{align*}
 finishing the proof. 
\end{proof}

\subsection{Vanishing of exponential sums}
Fix a nonzero $\bm{m}\in \Z^n$ and let $D=\|\bm{m}\|^2$ or $-\|\bm{m}\|^2$ so that $D\not\equiv 3\Mod{4}$. For later purpose, we record in this section the following vanishing statement regarding the exponential sum $\cS(d,\chi,\bm m)$ defined in \eqref{eq:Sdchi} for $\chi$ a Dirichlet character modulo $d$.  
 \begin{Prop}\label{prop: vfnp}
Let $\bm{m}\in \Z^n\setminus\{\bm{0}\}$, $D\in \{\pm \|\bm{m}\|^2\}$ and $\chi$ be as above. If $\chi\chi_{D}$ is principal but $\chi$ is non-principal, then $\cS(d,\chi,\bm m)=0$.
\end{Prop}

We will reduce the proof to the case when $d=p$ is a prime using the relation \eqref{equ:almultivp}. First we prove the following formula for $\varphi_n(p;\bm{m})$ appearing in the definition of $\cS(p,\chi,\bm{m})$. 
\begin{Lem}\label{lem:indoft}
For any odd prime $p$ and $\bm{m}\in \Z^n$, if  $\bm{m}\not\equiv \bm{0}\Mod{p}$ then
$$\vf_n(p;\bm{m})=\frac{\cG_p(1)^n}{p}\sum_{b\in (\Z/p\Z)^{\times}}\left(\tfrac{b}{p}\right)^ne\left(\tfrac{b-\bar{4}\bar{b}\|\bm{m}\|^2}{p}\right),$$
while for $\bm{m}\equiv \bm{0}\Mod{p}$
$$\vf_n(p;\bm{m})=p^{n-1}+\frac{\cG_p(1)^n}{p}\sum_{b\in (\Z/p\Z)^{\times}}\left(\tfrac{b}{p}\right)^n e\left(\tfrac{b}{p}\right),$$
where 
{for any $a\in\N$ and $b\in (\Z/a\Z)^{\times}$, $\cG_a(b):=\sum_{v\in \Z/a\Z}e\left(\frac{bv^2}{p}\right)$ is the generalized quadratic Gauss sum}, and $\bar{b}$ denotes the multiplicative inverse of $b$ modulo $p$.
\end{Lem}
\begin{proof}
{We first assume $\bm{m}\not\equiv \bm{0}\Mod{p}$.} Using the identity $\frac{1}{q}\sum_{a\in \Z/q\Z}e\left(\frac{a v}{q}\right)=\left\lbrace\begin{array}{cc}
1 & q\mid v\\
0 & q\nmid v \end{array}\right.$ for any $q\in \N$ we have
\begin{align*}
\vf_n(p; \bm{m})
&=\frac{1}{p}\sum_{b\in \Z/p\Z} \sum_{\bm{h}\in (\Z/p\Z)^n}e\left( \tfrac{\bm{m}\cdot\bm{h} +b(\|\bm{h}\|^2+1)}{p}\right)\\
&=\frac{1}{p}\sum_{b\in (\Z/p\Z)^{\times}}e\left( \tfrac{b-\bar{4}\bar{b}\|\bm{m}\|^2}{p}\right) \sum_{\bm{h}\in (\Z/p\Z)^n}e\left( \tfrac{b\|\bm{h}+\bar{2}\bar{b}\bm{m}\|^2}{p}\right)\\
&= \frac{1}{p}\sum_{b\in (\Z/p\Z)^{\times}}e\left( \tfrac{b-\bar{4}\bar{b}\|\bm{m}\|^2}{p}\right)\cG_p(b)^n.
\end{align*}
{Now recall the following formula for the generalized Gauss sum that
\begin{align}\label{equ:qgsodd}
\cG_a(b)=\varepsilon_a\sqrt{a}\left(\tfrac{b}{a}\right), \quad\forall\ \gcd(a,2b)=1.
\end{align}
Here $\varepsilon_a=1$ if $a\equiv 1\Mod{4}$ and $\varepsilon_a=i$ if $a\equiv 3\Mod{4}$; see e.g. \cite[p. 52]{IwaniecKowalski2004}.} 
In particular, this implies that $\cG_p(b)=\cG_p(1)\left(\tfrac{b}{p}\right)$ for any $b\in (\Z/p\Z)^{\times}$. Hence
\begin{align*}
\vf_n(p; \bm{m})&=\frac{\cG_p(1)^n}{p}\sum_{b\in (\Z/p\Z)^{\times}} \left(\tfrac{b}{p}\right)^ne\left( \tfrac{b-\bar{4}\bar{b}\|\bm{m}\|^2}{p}\right)
\end{align*}
as desired.
The formula for $\bm{m}\equiv \bm{0}\Mod{p}$ follows from the same argument, noting that in this case we have a nontrivial contribution from $b=0$.
\end{proof}
In particular if $p\mid \|\bm{m}\|^2$ we get that $\vf_n(p; t\bm{m})$ does not depend on $t\in (\Z/p\Z)^{\times}$. 
Using this observation we get the following.
\begin{Cor}\label{c:vanish}
Let $p$ be an odd prime number. If $p\mid \|\bm{m}\|^2$, then $\cS(p,\chi, \bm{m})=0$ for any non-principal character $\chi$ modulo $p$.
\end{Cor}
We can now give the

\begin{proof}[Proof of \propref{prop: vfnp}]
Writing $d=p_1\cdots p_k$ into product of distinct odd primes, we can write $\chi=\chi_1\cdots\chi_k$ with each $\chi_i$ a character of modulus $p_i$, then \eqref{equ:almultivp} and the Chinese Remainder theorem imply that 
$$\cS(d,\chi,\bm{m})=\prod_{i=1}^k \cS(p_i, \chi_i,\overline{d/p_i}\bm{m}),$$
where $\overline{d/p_i}$ denotes the multiplicative inverse of $\frac{d}{p_i}$ modulo $p_i$.
Hence to show that $\cS(d,\chi,\bm{m})=0$ it is enough to show that  $\cS(p_i, \chi_i,\overline{d/p_i}\bm{m})=0$ for some $1\leq i\leq k$, or in view of \corref{c:vanish} that for some factor $p_i$ we have that $p_i\mid \|\bm{m}\|^2$ and $\chi_i$ is non-principal. 

The condition that $\chi\chi_{D}$ is a principal character implies that $\chi$ is a quadratic character induced from $\chi_{D}$. (Recall that $\chi_D$ is the unique primitive character inducing $(\frac{D}{\cdot})$.)
Since we assume $\chi$ is non-principal, $\chi_D$ has modulus $q>1$ dividing $d$, implying that $q$ is odd and square-free. This then further implies that $q\mid |D|=\|\bm{m}\|^2$. 
Now if $p_i\mid q$ is any prime dividing the modulus of $\chi_D$, then $\chi_i$ is non-principal (as it is the quadratic character modulo $p_i$) and $p_i\mid \|\bm{m}\|^2$ as needed.
\end{proof}


\subsection{Evaluation of exponential sums}\label{sec:eoes}
Before we proceed with computing the local zeta functions we need to evaluate certain exponential sums corresponding to the arithmetic function $\phi_n(p^k;\bm{m})$ that we now introduce. Fix $\bm{m}\in \Z^n$ and a prime number $p$. For any $k\in \N$ and $0\leq i\leq k$ we define the exponential sum
\begin{align}\label{eq:Fki}
F(k,i):=p^{-k}\sum_{b\in (\Z/p^{k-i}\Z)^{\times}}\sum_{\bm{h}\in (\Z/p^k\Z)^n}e\left(\tfrac{b\|\bm{h}\|^2}{p^{k-i}}\right)e\left(\tfrac{\bm{h}\cdot\bm{m}}{p^k}\right).
\end{align}
We also introduce some notation that we fix throughout the remaining computations. For fixed $p$ and $\bm{m}$ as above, we denote by $\alpha_p:=\nu_p(\|\bm{m}\|^2)$ and $\ell_p:=\nu_p(\gcd(\bm{m}))$, where $\nu_p: \N\to \Z_{\geq 0}$ is the $p$-adic valuation function. When $\bm{m}=\bm{0}$, $\alpha_p$ and $\ell_p$ are both understood as $\infty$. When $\bm{m}\neq \bm{0}$ we denote by $t_p:=\frac{\|\bm{m}\|^2}{p^{\alpha_p}}$. When our discussion only involves one prime and there is no ambiguity,
we will omit the subscript from these parameters. 
Finally, for any $k\in \Z$ we denote by
$a_k:=(1+i)^k+(1-i)^k$ (with $i=\sqrt{-1}$ the imaginary unit) and $b_k:=\sum_{v\in (\Z/8\Z)^{\times}}e\left(\tfrac{vk}{8}\right)$.

The relevance of the sums $F(k,i)$ is given in the following simple formula for $\phi_n(p^k;\bm{m})$.
\begin{Lem}\label{l:phinFki}
For any $k\in\N$ we have
\begin{align*}
\phi_n(p^k;\bm{m})=\sum_{i=0}^{k}F(k,i).
\end{align*}
\end{Lem}
\begin{proof}
Similar to the proof of \lemref{lem:indoft} we have
\begin{align*}
\phi_n(p^k;\bm{m})&=\frac{1}{p^k}\sum_{a\in \Z/p^k\Z}\sum_{\bm{h}\in (\Z/p^k\Z)^n}e\left(\tfrac{a\|\bm{h}\|^2}{p^k}\right)e\left(\tfrac{\bm{h}\cdot\bm{m}}{p^k}\right)\\
&=\sum_{i=0}^kp^{-k}\sum_{\substack{a\in \Z/p^k\Z\\ \gcd(a,p^k)=p^i}}\sum_{\bm{h}\in (\Z/p^k\Z)^n}e\left(\tfrac{a\|\bm{h}\|^2}{p^k}\right)e\left(\tfrac{\bm{h}\cdot\bm{m}}{p^k}\right)\\
&=\sum_{i=0}^kp^{-k}\sum_{b\in (\Z/p^{k-i}\Z)^{\times}}\sum_{\bm{h}\in (\Z/p^k\Z)^n}e\left(\tfrac{b\|\bm{h}\|^2}{p^{k-i}}\right)e\left(\tfrac{\bm{h}\cdot\bm{m}}{p^k}\right),
\end{align*}
finishing the proof.
\end{proof} 

Before we evaluate these sums explicitly we give the following prelimimary formula. 

\begin{Lem}\label{lem:redufirst}
For any $k\in \N$ and $0\leq i\leq k$ we have 
\begin{align*}
F(k,i)&=\left\{\begin{array}{ll}   
				0 & i>\ell,\\
				p^{ni-k}\sum_{b\in (\Z/p^{k-i}\Z)^{\times}}\prod_{j=1}^n\left(\sum_{v\in \Z/p^{k-i}\Z}e\left(\tfrac{bv^2+vm_j/p^i}{p^{k-i}}\right)\right) & i\leq \ell.
			\end{array}\right.
\end{align*}\end{Lem}
\begin{proof}
First note that from definition of $F(k,i)$ we have
\begin{align*}
F(k,i)&=p^{-k}\sum_{b\in (\Z/p^{k-i}\Z)^{\times}}\prod_{j=1}^n\left(\sum_{h\in \Z/p^k\Z}e\left(\tfrac{bh^2}{p^{k-i}}\right)e\left(\tfrac{hm_j}{p^k}\right)\right)\\
&=p^{-k}\sum_{b\in (\Z/p^{k-i}\Z)^{\times}}\prod_{j=1}^n\left(\sum_{v\in \Z/p^{k-i}\Z}e\left(\tfrac{bv^2}{p^{k-i}}\right)\sum_{\substack{h\in \Z/p^k\Z\\ h\equiv v\Mod{p^{k-i}}}}e\left(\tfrac{hm_j}{p^k}\right)\right).
\end{align*}
Note that by writing $h=v+p^{k-i}l$ with $0\leq l\leq p^i-1$, the above innermost sum equals
\begin{align*}
e\left(\tfrac{vm_j}{p^k}\right)\sum_{l=0}^{p^{i}-1}e\left(\tfrac{l m_j}{p^i}\right).
\end{align*}
If $i>\ell$, then there exists $1\leq j\leq n$ such that $p^i\nmid m_j$. For such $j$ the above innermost sum vanishes, implying that in this case $F(k,i)=0$. If $i\leq \ell$, the above innermost sum equals $p^i$, implying that in this case
\begin{align*}
F(k,i)
&=p^{ni-k}\sum_{b\in (\Z/p^{k-i}\Z)^{\times}}\prod_{j=1}^n\left(\sum_{v\in \Z/p^{k-i}\Z}e\left(\tfrac{bv^2+vm_j/p^i}{p^{k-i}}\right)\right).
\end{align*}
This finishes the proof.
\end{proof}

	
{We now proceed to computing $F(k,i)$ explicitly. We first treat the case when $p$ is odd.} 
\begin{Prop}\label{prop:qusumgaoddp}
Assume further that $p$ is odd.
Then  $F(k,i)=0$ for all  $i>\ell$ and for  $0\leq i\leq \ell$,
\begin{enumerate}
\item[(i)] If $n$ is even, then
\begin{align*}
F(k,i)&=\left\{\begin{array}{ll}   
				\chi_{-4}(p)^{\frac{n(k-i)}{2}}p^{\frac{nk+ni-2k}{2}}\varphi(p^{k-i}) &  i\leq k\leq \alpha-i,\\
				-{\chi_{-4}(p)}^{\frac{n(k-i)}{2}}p^{\frac{nk+(n-2)i-2}{2}} & k= \alpha-i+1,\\
				0 & k\geq \alpha-i+2,
			\end{array}\right.
\end{align*}
where $\varphi$ is the Euler's totient function.
\item[(ii)] If $n$ is odd and $k-i$ is even, then 
\begin{align*}
F(k,i)&=\left\{\begin{array}{ll}   
				p^{\frac{nk+ni-2k}{2}}\varphi(p^{k-i})  &  i\leq k\leq \alpha-i,\\
				-p^{\frac{nk+(n-2)i-2}{2}} & k= \alpha-i+1,\\
				0 & k\geq \alpha-i+2.
			\end{array}\right.
\end{align*}
\item[(iii)] If both $n$ and $k-i$ are odd, then $F(k,i)=0$ except when $k=\alpha-i+1$ 
in which case 
\begin{align*}
F(k, i)&=\left(\tfrac{(-1)^{\frac{n-1}{2}}t}{p}\right)p^{\frac{nk+(n-2)i-1}{2}}.
\end{align*}
\end{enumerate}
\end{Prop}
\begin{proof}
The assertion that $F(k,i)=0$ when $i>\ell$ already follows from \lemref{lem:redufirst}; we thus assume $0\leq i\leq \ell$.
Since $i\leq \ell$, $m_j/p^i$ is an integer for each $1\leq j\leq n$. Applying the second formula in \lemref{lem:redufirst} and completing the square we get
\begin{align*}
F(k,i)&=p^{ni-k}\sum_{b\in (\Z/p^{k-i}\Z)^{\times}}\prod_{j=1}^n\left(\cG_{p^{k-i}}(b)e\left(\tfrac{-\overline{4b}(m_j/p^i)^2}{p^{k-i}}\right)\right)\\
&=p^{ni-k}\sum_{b\in (\Z/p^{k-i}\Z)^{\times}}\cG_{p^{k-i}}(b)^ne\left(-\tfrac{\overline{4b}\|\bm{m}/p^i\|^2}{p^{k-i}}\right)\\
&=p^{ni-k}\sum_{b\in (\Z/p^{k-i}\Z)^{\times}}\cG_{p^{k-i}}(b)^ne\left({-\tfrac{\overline{4b}t}{p^{k+i-\alpha}}}\right).
\end{align*}
Here $\overline{4b}$ is the multiplicative inverse of $4b$ modulo $p^{k-i}$ and $t=\|\bm m\|^2/p^{\alpha}$. 
From the quadratic Gauss sum formula \eqref{equ:qgsodd} we see that $\cG_{p^{k-i}}(b)=\cG_{p^{k-i}}(1)\left(\tfrac{b}{p}\right)^{k-i}$ for any $b\in (\Z/p^{k-i}\Z)^{\times}$. Hence
\begin{align*}
F(k,i)&=p^{ni-k}\cG_{p^{k-i}}(1)^n\sum_{b\in (\Z/p^{k-i}\Z)^{\times}}\left(\tfrac{b}{p}\right)^{(k-i)n}e\left(-{\tfrac{\overline{4b}t}{p^{k+i-\alpha}}}\right).
\end{align*}
Let us first assume $n$ is even so that
\begin{align*}
	F(k,i)&=p^{ni-k}\cG_{p^{k-i}}(1)^n\sum_{b\in (\Z/p^{k-i}\Z)^{\times}}e\left(-{\tfrac{\overline{4b}t}{p^{k+i-\alpha}}}\right)\\
	&=\left(\tfrac{-1}{p}\right)^{\frac{n(k-i)}{2}}p^{ni-k+\frac{n(k-i)}{2}}\sum_{b\in (\Z/p^{k-i}\Z)^{\times}}e\left(-{\tfrac{\overline{4b}t}{p^{k+i-\alpha}}}\right).
\end{align*}
Using the formula $\cG_{p^{k-i}}(1)^2=\left(\tfrac{-1}{p}\right)^{k-i}p^{k-i}$ (which follows from \eqref{equ:qgsodd}) we have if $i\leq \alpha-k$ (which is void if $k>\alpha$) then 
\begin{align*}
	F(k,i)&=
	\left(\tfrac{-1}{p}\right)^{\frac{n(k-i)}{2}}p^{\frac{nk+ni-2k}{2}}\varphi(p^{k-i}).
\end{align*}
If $i=\alpha-k+1$, then 
\begin{align*}
	F(k,i)&=\left(\tfrac{-1}{p}\right)^{\frac{n(k-i)}{2}}p^{ni-k+\frac{n(k-i)}{2}}\sum_{b\in (\Z/p^{k-i}\Z)^{\times}}e\left(-\tfrac{\overline{4b}t}{p}\right)\\
	&=\left(\tfrac{-1}{p}\right)^{\frac{n(k-i)}{2}}p^{ni-k+\frac{n(k-i)}{2}}\times p^{k-i-1}\sum_{b\in (\Z/p\Z)^{\times}}e\left(-\tfrac{\overline{4b}t}{p}\right)\\
	&=-\left(\tfrac{-1}{p}\right)^{\frac{n(k-i)}{2}}p^{\frac{nk+(n-2)i-2}{2}}.
\end{align*}
When $i\geq \alpha-k+2$ we have
\begin{align*}
	F(k,i)&=\left(\tfrac{-1}{p}\right)^{\frac{n(k-i)}{2}}p^{ni-k+\frac{n(k-i)}{2}}\times p^{\alpha-2i}\sum_{b\in (\Z/p^{k-\alpha+i}\Z)^{\times}}e\left(-\tfrac{\overline{4b}t}{p^{k-\alpha+i}}\right)=0.
\end{align*}
Here for the last equality we used the assumption that $k\geq 2+\alpha$ 
and the equality that
\begin{align*}
	\sum_{a\in (\Z/p^l\Z)^{\times}}e\left(\tfrac{a}{p^l}\right)=0
\end{align*}
whenever $l\geq 2$. {We can then finish the proof of this case by noting that $\chi_{-4}(p)=\left(\frac{-1}{p}\right)$}. The case when $n$ is odd and $k-i$ is even is identical to the previous case except now we have the formula $\cG_{p^{k-i}}(1)=p^{\frac{k-i}{2}}$ (since $k-i$ is even implies that $p^{k-i}\equiv 1\Mod{4}$). Next, we treat the case when both $n$ and $k-i$ are {odd}. In this case we have
\begin{align*}
F(k,i)&=p^{ni-k}\cG_{p^{k-i}}(1)^n\sum_{b\in (\Z/p^{k-i}\Z)^{\times}}\left(\tfrac{b}{p}\right)e\left(-\tfrac{\overline{4b}t }{p^{k+i-\alpha}}\right).
\end{align*} 
Now,  if $i\leq \alpha-k$ then
\begin{align*}
F(k,i)&=p^{ni-k}\cG_{p^{k-i}}(1)^n\sum_{b\in (\Z/p^{k-i}\Z)^{\times}}\left(\tfrac{b}{p}\right)=0.
\end{align*}
If $i\geq \alpha-k+1$ 
then 
\begin{align}\label{equ:innmossum}
F(k,i)&=p^{ni-k+\alpha-2i}\cG_{p^{k-i}}(1)^n\sum_{b\in (\Z/p^{k-\alpha+i}\Z)^{\times}}\left(\tfrac{b}{p}\right)e\left(-\tfrac{\overline{4b}t}{p^{k-\alpha+i}}\right)\\ 
&=p^{ni-k+\alpha-2i}\cG_{p^{k-i}}(1)^n\left(\tfrac{-t}{p}\right)\sum_{b\in (\Z/p^{k-\alpha+i}\Z)^{\times}}\left(\tfrac{b}{p}\right)e\left(\tfrac{b}{p^{k-\alpha+i}}\right),\nonumber
\end{align}
where for the second equality we applied a change of variable $-\overline{4b}t\mapsto b$ and used the relation $\left(\tfrac{-\overline{4b}t}{p}\right)=\left(\tfrac{-bt}{p}\right)=\left(\tfrac{-t}{p}\right)\left(\tfrac{b}{p}\right)$. If $i\geq \alpha-k+2$ we can group the above innermost sum by the congruence class of $b$ modulo $p$ to get
\begin{align*}
\sum_{b\in (\Z/p^{k-\alpha+i}\Z)^{\times}}\left(\tfrac{b}{p}\right)e\left(\tfrac{b}{p^{k-\alpha+i}}\right)&=\sum_{v\in (\Z/p\Z)^{\times}}\left(\tfrac{v}{p}\right)\sum_{\substack{b\in (\Z/p^{k-\alpha+i}\Z)^{\times} \\ b\equiv v\Mod{p}}}e\left(\tfrac{b}{p^{k-\alpha+i}}\right)=0.
\end{align*}
If $i= \alpha-k+1$ (which is possible only when $\alpha$ is even), using the identity \eqref{equ:aintermide} {below}
the innermost sum in \eqref{equ:innmossum} equals $\cG_p(1)$.
Hence we have in this case
\begin{align*}
F(k,i)&=p^{ni-k+k-i-1}\cG_{p^{k-i}}(1)^n\cG_p(1)\left(\tfrac{-t}{p}\right)=p^{ni-i-1}\varepsilon_p^{n+1}p^{\frac{n(k-i)+1}{2}}\left(\tfrac{-t}{p}\right)\\
&=p^{\frac{nk+(n-2)i-1}{2}}\left(\tfrac{-1}{p}\right)^{\frac{n+1}{2}}\left(\tfrac{-t}{p}\right)=p^{\frac{nk+(n-2)i-1}{2}}\left(\tfrac{(-1)^{\frac{n-1}{2}}t}{p}\right).\qedhere
\end{align*}
\end{proof}

When $p=2$ the formulas for these exponential sums are more involved. {For this we need the following quadratic Gauss sum formula when the moduli is a power of $2$: For any odd $b$,}
\begin{align}\label{equ:qgs2}
		\cG_{2^k}(b)=\left\{\begin{array}{ll}   
			2^{\frac{k}{2}}(1+e(\tfrac{b}{4})) & \text{$k$ is even},\\
			2^{\frac{k+1}{2}}e(\tfrac{b}{8}) & \text{$k>1$ is odd},\\
			0 & k=1.
		\end{array}\right.
	\end{align}

\begin{Prop}\label{prop:fkifortwo}
{Assume $p=2$}. 
We have that $F(k,i)=0$ for $i>\ell$. For $i=\ell$ we have $F(\ell,\ell)=2^{(n-1)\ell}$, $F(k,\ell)=0$ for all 
$k> \ell+1$ and $F(\ell+1,\ell)=0$ unless $2^{-\ell}\bm{m}$ has only odd entries in which case $F(\ell+1,\ell)=2^{(n-1)(\ell+1)}$.
Next {recall that} $a_k=(1+i)^k+(1-i)^k$ and $b_k=\sum_{v\in (\Z/8\Z)^{\times}}e\left(\tfrac{vk}{8}\right)$. 
When $\ell\geq 1$ for any $0\leq i\leq \ell-1$ 
we have
\begin{enumerate}
\item[(i)]  If $i\leq k\leq \alpha-i-2$, then
\begin{align*}
F(k,i)&=
\left\{\begin{array}{ll}
				2^{(n-1)k} & k=i,\\
				0 & k=i+1,\\
				2^{\frac{nk+(n-2)i-4}{2}}a_n & \text{$k-i\geq 2$ is even},\\
				2^{\frac{n(k+1)+(n-2)i-6}{2}}b_n & \text{$k-i\geq 3$ is odd}.
			\end{array}\right.
\end{align*}
\item[(ii)] If $k=\alpha-i-1$ then
\begin{align*}
F(k,i)&=
\left\{\begin{array}{ll}
                                     0 & (\alpha, i)=(2\ell, \ell-1),\\
				-2^{\frac{nk+(n-2)i+n-6}{2}}b_n & (\alpha, i)\neq (2\ell, \ell-1)\ \text{and $\alpha$ is even}, \\
				-2^{\frac{nk+(n-2)i-4}{2}}a_n & \text{$\alpha$ is odd}.
			\end{array}\right.
		\end{align*}

\item[(iii)] If $k=\alpha-i$  then 
$$F(k,i)=\left\{\begin{array}{ll}
				2^{\frac{nk+(n-2)i-6}{2}}(-1)^{\frac{t+1}{2}} a_{n+2}
				& \text{$\alpha$ is even},\\
				2^{\frac{n(k+1)+(n-2)i-6}{2}} b_{n-2t} 
				 & \text{$\alpha$ is odd}.
			\end{array}\right.$$
\item[(iv)] If $k=\alpha-i+1$ then 
$$F(k,i)=\left\{\begin{array}{ll}
				2^{\frac{n(k+1)+(n-2)i -6}{2}}b_{n-t}
				 & \text{$\alpha$ is even},\\
				0 & \text{$\alpha$ is odd}.
			\end{array}\right.$$
\item[(v)] If $k\geq \alpha-i+2$, then $F(k,i)=0$.
\end{enumerate}
\end{Prop}
\begin{proof}
Again the assertion that $F(k,i)=0$ when $i>\ell$ already follows from \lemref{lem:redufirst} and we assume $0\leq i\leq \ell$. In this case by \lemref{lem:redufirst} we have
\begin{align*}
F(k,i)&=2^{ni-k}\sum_{b\in (\Z/2^{k-i}\Z)^{\times}}\prod_{j=1}^n\left(\sum_{v\in \Z/2^{k-i}\Z}e\left(\tfrac{bv^2+vm_j/2^i}{2^{k-i}}\right)\right).
\end{align*}
Now we treat the case when $i=\ell$. First it is clear from the above expression that $F(\ell, \ell)=2^{(n-1)\ell}$. Next, for $k\geq \ell+2$, since {$\ell=\nu_2(\gcd(\bm{m}))$}, there exists $1\leq j\leq n$ such that $m_j/2^{\ell}$ is odd. For such $j$, the above inner sum (with $i=\ell$) equals
\begin{align*}
\sum_{h\in \Z/2^{k-\ell-1}\Z}e\left(\tfrac{bh^2}{2^{k-\ell}}\right)\sum_{\substack{v\in \Z/2^{k-\ell}\Z\\ v\equiv h\Mod{2^{k-\ell-1}}}}e\left(\tfrac{vm_j/2^{\ell}}{2^{k-\ell}}\right)=0,
\end{align*}
implying that $F(k,\ell)=0$.
Here we used the observation that for any $l\in \N$ and $v\in \Z/2^l\Z$, $v^2$ is a well-defined element in $\Z/2^{l+1}\Z$. Next, assume $k=\ell+1$. In this case {$b=1$ and} the above inner sum equals
\begin{align*}
\sum_{v\in \Z/2\Z}e\left(\tfrac{v^2+vm_j/2^{\ell}}{2}\right)=1+e\left(\tfrac{1+m_j/2^{\ell}}{2}\right),\quad 1\leq j\leq n.
\end{align*}
If there exists some $1\leq j\leq n$ such that $m_j/2^{\ell}$ is even, then the above sum vanishes, implying that $F(k,\ell)=0$ as in the previous case. If all $m_j/2^{\ell}$ are odd, then the above sum equals $2$, implying that in this case
\begin{align*}
F(k,\ell)=2^{n\ell-k}\times 2^n=2^{(n-1)k}=2^{(n-1)(\ell+1)}. 
\end{align*} 
This proves the formulas for $F(k,\ell)$. Finally, we assume $0\leq i\leq \ell-1$ in which case we can compute $F(k,i)$ using the completing-square trick as in the proof of \propref{prop:qusumgaoddp}.
Since $i\leq \ell-1$, we can rewrite $vm_j/2^i=2vm_j/2^{i+1}$ for each $1\leq j\leq n$ to get
\begin{align*}
F(k,i)
&=2^{ni-k}\sum_{b\in (\Z/2^{k-i}\Z)^{\times}}\prod_{j=1}^n\left(\sum_{v\in \Z/2^{k-i}\Z}e\left(\tfrac{bv^2+2vm_j/2^{i+1}}{2^{k-i}}\right)\right).
\end{align*}
Completing the square we get
\begin{align*}
F(k,i)&=2^{ni-k}\sum_{b\in (\Z/2^{k-i}\Z)^{\times}}\cG_{2^{k-i}}(b)^ne\left(-\tfrac{\overline{b}\|\bm{m}/2^{i+1}\|^2}{2^{k-i}}\right).
\end{align*}
{Recall that when $\bm{m}\neq \bm{0}$ we have $t=\|\bm{m}\|^2/2^{\alpha}$. Hence if $\bm{m}\neq\bm{0}$ we have}
\begin{align*}
\frac{\|\bm{m}/2^{i+1}\|^2}{2^{k-i}}=\frac{t}{2^{k+i-\alpha+2}},
\end{align*}
implying that 
\begin{align*}
F(k,i)&=2^{ni-k}\sum_{b\in (\Z/2^{k-i}\Z)^{\times}}\cG_{2^{k-i}}(b)^ne\left(-\tfrac{\overline{b}t}{2^{k+i-\alpha+2}}\right).
\end{align*}
Hence if $i\leq k\leq \alpha-i-2$, by applying \eqref{equ:qgs2} we have in this case
\begin{align*}
F(k,i)&=2^{ni-k}\sum_{b\in(\Z/2^{k-i}\Z)^{\times}}\cG_{2^{k-i}}(b)^n=\left\{\begin{array}{ll}
				2^{(n-1)k} & k=i,\\
				0 & k=i+1,\\
				2^{\frac{nk+(n-2)i-4}{2}}a_n & \text{$k-i\geq 2$ is even},\\
				2^{\frac{n(k+1)+(n-2)i-6}{2}}b_n & \text{$k-i\geq 3$ is odd}.
			\end{array}\right.
\end{align*}
If $k=\alpha-i-1$, then $e\left(-\tfrac{\overline{b}t}{2^{k+i-\alpha+2}}\right)=e\left(\tfrac{-\overline{b}t}{2}\right)=-1$ and by \eqref{equ:qgs2} we have
\begin{align*}
F(k,i)&=-2^{ni-k}\sum_{b\in (\Z/2^{k-i}{\Z})^{\times}}\cG_{2^{k-i}}(b)^n=\left\{\begin{array}{ll}
                                     0 & (\alpha, i)=(2\ell, \ell-1),\\
				-2^{\frac{n(k+1)+(n-2)i-6}{2}}b_n & (\alpha, i)\neq (2\ell, \ell-1) \ \text{and}\ \text{$\alpha$ is even},\\
				-2^{\frac{nk+(n-2)i-4}{2}}a_n & \text{$\alpha$ is odd}.
			\end{array}\right.
\end{align*}
{Here we also used that $\alpha\geq 2\ell$.}
If $k=\alpha-i$  then by applying \eqref{equ:qgs2}  we have (noting that in this case $k-i\geq 2$)
\begin{align*}
F(k,i)&=\left\{\begin{array}{ll}
				(-1)^{\frac{t+1}{2}}2^{\frac{nk+(n-2)i-4}{2}}i\left((1+i)^n-(1-i)^n\right) & \text{$\alpha$ is even},\\
				2^{\frac{n(k+1)+(n-2)i-6}{2}}\sum_{b\in (\Z/8\Z)^{\times}}e\left(\tfrac{b(n-2t)}{8}\right) & \text{$\alpha$ is odd}.
			\end{array}\right.
\end{align*}
We then get the desired formula by noting that 
$$
a_{n+2}=(1+i)^{n+2}+(1-i)^{n+2}=2i\left((1+i)^n-(1-i)^n\right).
$$
If $i=\alpha-k+1$ then again by applying \eqref{equ:qgs2} we have
\begin{align*}
F(k,i)&=\left\{\begin{array}{ll}
				2^{\frac{n(k+1)+(n-2)i-6}{2}}\sum_{b\in (\Z/8\Z)^{\times}}e\left(\tfrac{b(n-t)}{8}\right) & \text{$\alpha$ is even},\\
				0 & \text{$\alpha$ is odd}.
			\end{array}\right.
\end{align*}
Finally, for $k\geq \alpha-i+2$, making a change of variable $b\mapsto \overline{b}$ and noting that \eqref{equ:qgs2} implies that $\cG_{2^{k-i}}(b)$ is well-defined for $b$ modulo $8$ we get
\begin{align*}
F(k,i)&=2^{ni-k+\alpha-2i-2}\sum_{v\in (\Z/8\Z)^{\times}}\cG_{2^{k-i}}(\overline{v})^n\sum_{\substack{b\in (\Z/2^{k+i-\alpha+2}\Z)^{\times}\\ b\equiv v\Mod{8}}}e\left(\tfrac{-bt}{2^{k+i-\alpha+2}}\right)=0.
\end{align*}
Here for the last equality we used that $k+i-\alpha+2\geq 4$. This finishes the proof.
\end{proof}

\begin{remark}\label{rmk:abn}
For later reference, we note that one easily sees that 
\begin{align*}
b_k&=\left\{\begin{array}{ll}
				(-1)^{\frac{k}{4}}4 & k\equiv 0\Mod{4},\\
				0 & k\not\equiv 0\Mod{4},
			\end{array}\right.
\end{align*}
and using induction and the identity $(1+i)^4=(1-i)^4=-4$ we have for any $k\geq 0$
		\begin{align*}
			a_k=\left\{\begin{array}{ll}
			(-1)^{\frac{k}{4}}2^{\frac{k}{2}+1} & k\equiv 0\Mod{4},\\
				0 & k\equiv 2\Mod{4},\\
				\mathfrak{s}_k 2^{\frac{k+1}{2}} & k\equiv 1\Mod{2}.
			\end{array}\right.
		\end{align*} 
	Here recall that when $k$ is odd, $\mathfrak{s}_k=(-1)^{\frac{k^2-1}{8}}$.
\end{remark}

\subsection{Computing the local zeta functions}\label{sec:localzeta}
In this subsection we conclude the calculation of the Fourier coefficients by calculating the local zeta functions appearing in \propref{prop:eulprodconst} and \propref{prop:nonconsprform}.
{As before $p$ is a fixed prime and starting from now we take $\bm{m}=2\bm{\lambda}$ for some fixed $\bm{\lambda}\in \Lambda^*$.}
We first consider the case of $\bm\lambda=\bm 0$ {so that} $\alpha=\ell=\infty$. 

\begin{Prop}\label{prop:loczetconst}
For any odd prime $p$ and $\Re(s)>n$ we have 
$$Z_n^{(p)}(s;\bm{0})=\left\lbrace\begin{array}{cc}
\frac{1-\chi_{-4}(p)^{\frac{n}{2}}p^{\frac{n}{2}-1-s}}{(1-p^{n-1-s})(1-\chi_{-4}(p)^{\frac{n}{2}}p^{\frac{n}{2}-s})} & \text{$n$ even},\\
\frac{1-p^{n-1-2s}}{(1-p^{n-2s})(1-p^{n-1-s})} & \text{$n$ odd}, \end{array}\right.$$
and
\begin{align*}
\varphi_n(p;\bm{0})+\tilde{Z}_n^{(p)}(s;\bm{0})&=\left\lbrace\begin{array}{cc}
\frac{p^{n-1}+p^{n-s}-\chi_{-4}(p)^{\frac{n}{2}}p^{\frac{n}{2}-1}-p^{2n-1-2s}}{(1-p^{n-1-s})(1-\chi_{-4}(p)^{\frac{n}{2}}p^{\frac{n}{2}-s})} & \text{$n$ even},\\
\frac{(1-\chi_{-4}(p)^{\frac{n+1}{2}}p^{\frac{n-1}{2}-s})(p^{n-s}+p^{n-1}+\chi_{-4}(p)^{\frac{n+1}{2}}p^{\frac{n-1}{2}}-p^{2n-1-2s})}{(1-p^{n-2s})(1-p^{n-1-s})} & \text{$n$ odd}. \end{array}\right.
\end{align*}
For $p=2$ we have
\begin{align*}
\varphi_n(2;\bm{0})+\tilde{Z}_n^{(2)}(s;\bm0)=2^{n-1}\times \left\lbrace \begin{array}{ll}
	\frac{1-\left(1-(-1)^{\frac{n}{4}}\right)2^{\frac{n}{2}-s}}{(1-2^{n-1-s})(1-2^{\frac{n}{2}-s})}& n\equiv 0 \Mod{4},\\
\frac{1}{1-2^{n-1-s}} & n\equiv 2\Mod{4},\\
\frac{1-2^{n-2s}+(-1)^{\frac{n-1}{4}}2^{\frac{n-1}{2}-s}}{(1-2^{n-1-s})(1-2^{n-2s})} & n\equiv  1\Mod{4},\\
\frac{1-2^{n-2s}+(-1)^{\frac{n+1}{4}}2^{\frac{n-1}{2}-s}}{(1-2^{n-1-s})(1-2^{n-2s})} & n\equiv 3\Mod{4}.
  \end{array}\right.
\end{align*}
\end{Prop}
For the computation of the above second formula we also need a formula for the difference $\varphi_n(p;\bm{0})-\phi_n(p;\bm{0})$
 given in the following. 

\begin{Lem}\label{lem:diffcomp}
Let $p$ be an odd prime. Then
\begin{align*}
\varphi_n(p;\bm{0})-\phi_n(p;\bm{0})&=\left\{\begin{array}{ll}   
			-{\chi_{-4}(p)}^{\frac{n}{2}}p^{\frac{n}{2}} & \text{$n$ even},\\
			{\chi_{-4}(p)}^{\frac{n+1}{2}}p^{\frac{n-1}{2}} & \text{$n$ odd}.
		\end{array}\right.
\end{align*}
\end{Lem}
\begin{proof}
By \propref{prop:qusumgaoddp} (recalling that for $\bm{\lambda}=\bm{0}$, $\alpha=\ell=\infty$) we have
\begin{align*}
\phi_n(p;\bm{0})&=F(1,0)+F(1,1)=\left\{\begin{array}{ll}   
			p^{n-1}+{\chi_{-4}(p)}^{\frac{n}{2}}p^{\frac{n-2}{2}}(p-1) & \text{$n$ even},\\
			p^{n-1} & \text{$n$ odd},
		\end{array}\right.
\end{align*}
and by  \lemref{lem:indoft} we have
\begin{align*}
\vf_n(p;\bm 0)
&=p^{n-1}+\frac{\cG_p(1)^n}{p}\sum_{b\in (\Z/p\Z)^{\times}}e\left(\tfrac{b}{p}\right)\left(\tfrac{b}{p}\right)^n.
\end{align*}
If $n$ is even then
\begin{align*}
\varphi_n(p;\bm{0})&=p^{n-1}-\frac{\cG_p(1)^n}{p}=p^{n-1}-\left(\tfrac{-1}{p}\right)^{\frac{n}{2}}p^{\frac{n-2}{2}},
\end{align*}
while for odd $n$ 
\begin{align*}
\varphi_n(p;\bm{0})&
=p^{n-1}+\frac{\cG_p(1)^n}{p}\sum_{{b}\in \Z/p\Z}e\left(\tfrac{b}{p}\right)\left(\tfrac{b}{p}\right).
\end{align*}
Noting that  $\#\{v\in \Z/p\Z: v^2=b\}=1+\left(\frac{b}{p}\right)$ for any $b\in \Z/p\Z$ we see that 
\begin{align}\label{equ:aintermide}
\sum_{b\in \Z/p\Z}\left(\tfrac{b}{p}\right)e\left(\tfrac{b}{p}\right)&=\sum_{b\in \Z/p\Z}\left(1+\left(\tfrac{b}{p}\right)\right)e\left(\tfrac{b}{p}\right)=\cG_p(1),
\end{align}
and the formula for odd $n$ can be evaluated as 
\begin{align*}
\varphi_n(p;\bm{0})&=p^{n-1}+\frac{\cG_p(1)^{n+1}}{p}=p^{n-1}+{\chi_{-4}(p)}^{\frac{n+1}{2}}p^{\frac{n-1}{2}}.
\end{align*}
Subtracting the above formulas for $\varphi_n(p;\bm{0})$ and $\phi_n(p;\bm{0})$ concludes the proof.
\end{proof}

We now give the 
\begin{proof}[Proof of \propref{prop:loczetconst}]
From the relation $\phi_n(p^k;\bm{0})=\sum_{i=0}^kF(k,i)$ we have
\begin{align*}
Z_n^{(p)}(s;\bm{0})&=\sum_{k=0}^{\infty}\sum_{i=0}^{k}F(k,i)p^{-ks}=\sum_{i=0}^{\infty}\sum_{k=i}^{\infty}F(k,i)p^{-ks}.
\end{align*}
Here we can change of the order of summation since the above first double series absolute converges for $\Re(s)>n$ which can be verified using the trivial bound $|F(k,i)|\leq p^{(n-1)k}\varphi(p^{k-i})$. We first treat the case when $p$ is odd. Since $\bm{\lambda}=\bm{0}$ we have $\ell=\infty$. Applying case (i) of \propref{prop:qusumgaoddp} 
we have that if $n$ is even 
\begin{align*}
Z_n^{(p)}(s;\bm{0})&=\sum_{i=0}^{\infty}\sum_{k=i}^{\infty}\chi_{-4}(p)^{\frac{n(k-i)}{2}}p^{\frac{nk+ni-2k}{2}}\varphi(p^{k-i})\times p^{-ks}\\ 
&=\sum_{i=0}^{\infty}p^{(n-1-s)i}+(1-p^{-1})\sum_{i=0}^{\infty}\chi_{-4}(p)^{\frac{ni}{2}}p^{\frac{(n-2)i}{2}}\sum_{k=i+1}^{\infty}\chi_{-4}(p)^{\frac{nk}{2}}p^{(\frac{n}{2}-s)k}\\
&=\frac{1}{1-p^{n-1-s}}+\frac{(1-p^{-1})\chi_{-4}(p)^{\frac{n}{2}}p^{\frac{n}{2}-s}}{1-\chi_{-4}(p)^{\frac{n}{2}}p^{\frac{n}{2}-s}}\sum_{i=0}^{\infty}p^{(n-1-s)i}\\
&=\frac{1-\chi_{-4}(p)^{\frac{n}{2}}p^{\frac{n}{2}-s}+(1-p^{-1})\chi_{-4}(p)^{\frac{n}{2}}p^{\frac{n}{2}-s}}{(1-p^{n-1-s})(1-\chi_{-4}(p)^{\frac{n}{2}}p^{\frac{n}{2}-s})}\\
&=\frac{1-\chi_{-4}(p)^{\frac{n}{2}}p^{\frac{n}{2}-1-s}}{(1-p^{n-1-s})(1-\chi_{-4}(p)^{\frac{n}{2}}p^{\frac{n}{2}-s})}.
\end{align*}
This proves the case when $n$ is even.
Next, we treat the case when $n$ is odd in which case, applying case (ii) and (iii) of \propref{prop:qusumgaoddp} we get
\begin{align*}
Z_n^{(p)}(s;\bm{0})&=\sum_{i=0}^{\infty}\sum_{\substack{k\geq i\\ k\equiv i\Mod{2}}}p^{\frac{nk+ni-2k}{2}}\varphi(p^{k-i})\times p^{-ks}\\
&=\frac{1}{1-p^{n-1-s}}+(1-p^{-1})\sum_{i=0}^{\infty}p^{(n-1-s)i}\sum_{k=1}^{\infty}p^{(n-2s)k}\\
&=\frac{1}{1-p^{n-1-s}}+\frac{(1-p^{-1})p^{n-2s}}{(1-p^{n-1-s})(1-p^{n-2s})}\\
&=\frac{1-p^{n-1-2s}}{(1-p^{n-1-s})(1-p^{n-2s})},
\end{align*}
as claimed. For the second formula note that 
$$
\tilde{Z}_n^{(p)}(s;\bm{0})=\sum_{k=0}^\infty \frac{\phi_n(p^{k+2};\bm{0})}{p^{(k+1)s}}=p^s(Z_n^{(p)}(s;\bm{0})-1)-\phi_n(p;\bm{0}),
$$ we have
\begin{align*}
\vf_n(p; \bm 0)+ \tilde{Z}_n^{(p)}(s;\bm{0})&=p^s(Z_n^{(p)}(s;\bm{0})-1)+\vf_n(p;\bm 0)-\phi_n(p;\bm 0).
\end{align*}
The formula of $\varphi_n(p;\bm{0})+\tilde{Z}_n^{(p)}(s;\bm{0})$ then immediately follows from the formula for $Z_n^{(p)}(s;\bm{0})$ and the formula for $\varphi_n(p;\bm{0})-\phi_n(p;\bm{0})$ given in \lemref{lem:diffcomp}.

Next, we treat the case when $p=2$. Applying \propref{prop:fkifortwo} we have
\begin{align*}
Z_n^{(2)}(s;\bm{0})
&=\sum_{i=0}^{\infty}\left(2^{(n-1-s)i}+2^{-2}a_n\sum_{\substack{k\geq i+2\\ k\equiv i\Mod{2}}}2^{\frac{nk+(n-2)i}{2}-ks}+2^{\frac{n-6}{2}}b_n\sum_{\substack{k\geq i+2\\ k\not\equiv i\Mod{2}}}2^{\frac{nk+(n-2)i}{2}{-ks}}\right)\\
&=:\frac{1}{1-2^{n-1-s}}+2^{-2}a_nI_1+2^{\frac{n-6}{2}}b_nI_2.
\end{align*} 
For $I_1$ we have
\begin{align*}
I_1&=2^{n-2s}\sum_{i=0}^{\infty}2^{(n-1-s)i}\sum_{k=0}^{\infty}2^{(n-2s)k}=\frac{2^{n-2s}}{(1-2^{n-1-s})(1-2^{n-2s})}.
\end{align*}
Similarly, for $I_2$ we have
\begin{align*}
I_2&=2^{\frac{3n-6s}{2}}\sum_{i=0}^{\infty}2^{(n-1-s)i}\sum_{k=0}^{\infty}2^{(n-2s)k}=\frac{2^{\frac{3n-6s}{2}}}{(1-2^{n-1-s})(1-2^{n-2s})}.
\end{align*}
Collecting all the terms we get
\begin{align*}
Z_n^{(2)}(s;\bm{0})&=
\frac{1-2^{n-2s}+2^{n-2-2s}a_n+2^{2n-3-3s}b_n}{(1-2^{n-1-s})(1-2^{n-2s})}.
\end{align*}
Plugging the formulas of $a_n$ and $b_n$ (see \rmkref{rmk:abn}) into the above formula and using the relation $\varphi_n(2;\bm{0})+\tilde{Z}_n^{(2)}(s;\bm{0})=2^s(Z_n^{(2)}(s;\bm{0})-1)$
we get the desired formula for $\varphi_n(2;\bm{0})+\tilde{Z}_n^{(2)}(s; \bm{0})$. 
Here for the above relation we used the identity $\varphi_n(2;\bm{0})=\phi_n(2;\bm{0})$ (since both equal $2^{n-1}$; see \eqref{equ:evaat2}).
\end{proof}

For the non-constant coefficients we let  $\bm{m}=2\bm{\lambda}$ where $\bm\lambda\in\Lambda^*$ is  nonzero, so  that both $\alpha=\nu_p(\|2\bm{\lambda}\|^2)$ and $\ell=\nu_p(\gcd(2\bm\lambda))$ are finite and $t=\frac{\|2\bm\lambda\|^2}{p^{\alpha}}$ is defined. 

\begin{remark}\label{rmk:dichadsqu}
{Recall that for odd $n$ we have $D=(-1)^{\frac{n-1}{2}}\|2\bm{\lambda}\|^2$ and that $\chi_D$ is the primitive character inducing the Kronecker symbol $\left(\frac{D}{\cdot}\right)$, and note that when $\alpha$ is even (so that $D/t=p^{\alpha}$ is an even power of an odd prime), we have $\chi_D(p)=\left(\frac{(-1)^{\frac{n-1}{2}}t}{p}\right)$.}
\end{remark}

\begin{Lem}\label{lem:localzeta}
The local zeta function $Z_n^{(p)}(s;\chi,\bm{\lambda})$ (resp. $\tilde{Z}_n^{(p)}(s; \chi, \bm{\lambda})$) is a polynomial in $p^{-s}$ of degree $\leq \alpha+1$ (resp. $\leq \alpha$) and satisfies that for all $\sigma<n$
$$|Z_n^{(p)}(\sigma+it;\chi,\bm{\lambda})|\leq (\alpha+2)p^{(\alpha+1)(n-\sigma)}\quad \text{and}\quad |\tilde{Z}_n^{(p)}(\sigma+it;\chi, \bm{\lambda})|\leq (\alpha+1)p^{n+(\alpha+1)(n-\sigma)}.$$
\end{Lem}
\begin{proof}
The observation that  $F(k,i)=0$ for $k\geq \alpha-i+2$, together with \lemref{l:phinFki}, implies that $\phi_n(p^k;2\bm{\lambda})=0$ when $k\geq \alpha+2$ so that  
\begin{align}\label{equ:apreforodp}
Z_n^{(p)}(s; \chi, \bm{\lambda})=\sum_{k=0}^{\alpha+1}\frac{\chi(p^k)\phi_n(p^k;2\bm{\lambda})}{p^{ks}}.
\end{align}
The bound then follows from the trivial bound $|\chi(p^k)\phi_n(p^k;2\bm{\lambda})|\leq p^{nk}$. The statements for $\tilde{Z}_n^{(p)}(s; \chi, \bm{\lambda})$ follow similarly.
\end{proof}
\begin{Lem}\label{lem:lozetanonzero}
For any odd prime $p$ not dividing $\|2\bm{\lambda}\|^2$ we have
\begin{align*}
Z_n^{(p)}(s; \chi, \bm{\lambda})&=
\left\lbrace \begin{array}{ll}
1-\chi_{-4}(p)^{\frac{n}{2}}\chi(p)p^{\frac{n-2}{2} -s} & \text{$n$ even},\\ 
1+\chi\chi_{D}(p)p^{\frac{n-1}{2} -s} & \text{$n$ odd}.
\end{array}\right. 
\end{align*}
\end{Lem}
\begin{proof}
{When $p\nmid \|2\bm{\lambda}\|^2$ we have $\alpha=0$ and  by \eqref{equ:apreforodp} } $Z_n^{(p)}(s;\chi,\bm{\lambda})=1+\chi(p)\phi_n(p;2\bm\lambda)p^{-s}$. Applying case (i) and (iii) of \propref{prop:qusumgaoddp} respectively {(see also \rmkref{rmk:dichadsqu})} we have that for such $p$
\begin{align*}
\phi_n(p;2\bm{\lambda})&=F(1,0)=\left\lbrace \begin{array}{ll}
-{\chi_{-4}(p)}^{\frac{n}{2}}p^{\frac{n-2}{2}} & \textrm{$n$ even},\\ 
{\chi_D(p)}p^{\frac{n-1}{2}} & \textrm{$n$ odd}. 
\end{array}\right.  \qedhere
\end{align*}
\end{proof}
Next we treat the case of odd primes $p$ dividing $\|2\bm\lambda\|^2$. Here the formula is more involved and we only calculate it for the case when $\chi$ is the trivial character (though a similar calculation will also work for general $\chi$; see \rmkref{rmk:genchi} below). 
{In the following computations, we will use freely \propref{prop:qusumgaoddp} and the notation therein. Especially, by case (i) we mean case (i) in \propref{prop:qusumgaoddp}, and the same for other cases.} 
As before, in this case we abbreviate $Z_n^{(p)}(s;\chi,\bm{\lambda})$ by $Z_n^{(p)}(s;\bm{\lambda})$. 
\begin{Prop}\label{prop:lozetanonzeroramified}
Suppose $p$ is an odd prime dividing $\|2\bm\lambda\|^2$. 
Then we have that for $n$ even
\begin{align*}
	Z_n^{(p)}(s;\bm{\lambda})&=\frac{1-\chi_{-4}(p)^{\frac{n}{2}}p^{\frac{n}{2}-1-s}}{1-\chi_{-4}(p)^{\frac{n}{2}}p^{\frac{n}{2}-s}}\left(\frac{1-p^{(n-1-s)(\ell+1)}}{1-p^{n-1-s}}-\frac{\chi_{-4}(p)^{\frac{n(\alpha+1)}{2}}p^{(\frac{n}{2}-s)(\alpha+1)}(1-p^{(s-1)(\ell+1)})}{1-p^{s-1}}\right).
	\end{align*}
For $n$ odd, if $\alpha$ is odd then 
\begin{align*}
	Z_n^{(p)}(s;\bm{\lambda}) =&
			\frac{1-p^{n-1-2s}}{1-p^{n-2s}}\left(\frac{1-p^{(n-1-s)(\ell+1)}}{1-p^{n-1-s}} -\frac{p^{(\frac{n}{2}-s)(\alpha+1)}(1-p^{(s-1)(\ell+1)})}{1-p^{s-1}}\right),
	\end{align*} 
and if $\alpha$ is even, then 
	\begin{align*}
	Z_n^{(p)}(s;\bm{\lambda}) =&
			\tfrac{1+\chi_D(p) p^{\frac{n-1}{2}-s}}{1-p^{n-2s}}\left(\tfrac{(1-\chi_D(p) p^{\frac{n-1}{2}-s})(1-p^{(n-1-s)(\ell+1)})}{1-p^{n-1-s}}+ \tfrac{\chi_D(p) p^{(\frac{n}{2}-s)(\alpha+1)-\frac12}(1-\chi_D(p) p^{\frac{n+1}{2}-s})(1-p^{(s-1)(\ell+1)})}{1-p^{s-1}}\right).
	\end{align*}
	\end{Prop}

\begin{proof}
We have $F(k,i)=0$ for $i>\ell$ as well as for $0\leq i\leq \ell$ whenever $k\geq \alpha-i+2$. Hence
\begin{align*}
Z_n^{(p)}(s;\bm{\lambda})&=\sum_{i=0}^{\ell}\sum_{k=i}^{\alpha-i+1}F(k,i)p^{-ks}\\
&=\sum_{i=0}^{\ell}\left(F(i,i)p^{-is}+\sum_{k=i+1}^{\alpha-i}F(k,i)p^{-ks}+F(\alpha-i+1,i)p^{-(\alpha-i+1)s}\right)=:I_1+I_2+I_3.
\end{align*} 
For $I_1$ we have $F(i,i)=p^{(n-1)i}$ for all $0\leq i\leq \ell$, hence
\begin{align*}
I_1&=\sum_{i=0}^{\ell}p^{(n-1-s)i}=\frac{1-p^{(n-1-s)(\ell+1)}}{1-p^{n-1-s}}.
\end{align*}
For $I_2$ and $I_3$ the result depends on the parity of $n$. First we assume $n$ is even. For $I_2$ 
we have
\begin{align*}
I_2&=(1-p^{-1})\sum_{i=0}^{\ell}\chi_{-4}(p)^{\frac{ni}{2}}p^{\frac{(n-2)i}{2}}\sum_{k=i+1}^{\alpha-i}\chi_{-4}(p)^{\frac{nk}{2}}p^{(\frac{n}{2}-s)k}\\
&=\frac{1-p^{-1}}{1-\chi_{-4}(p)^{\frac{n}{2}}p^{\frac{n}{2}-s}}\left(\chi_{-4}(p)^{\frac{n}{2}}p^{\frac{n}{2}-s}\sum_{i=0}^{\ell}p^{n-1-s}-\chi_{-4}(p)^{\frac{n(\alpha+1)}{2}}p^{(\frac{n}{2}-s)(\alpha+1)}\sum_{i=0}^{\ell}p^{(s-1)i}\right)\\
&=\frac{1-p^{-1}}{1-\chi_{-4}(p)^{\frac{n}{2}}p^{\frac{n}{2}-s}}\left(\frac{\chi_{-4}(p)^{\frac{n}{2}}p^{\frac{n}{2}-s}(1-p^{(n-1-s)(\ell+1)})}{1-p^{n-1-s}}-\frac{\chi_{-4}(p)^{\frac{n(\alpha+1)}{2}}p^{(\frac{n}{2}-s)(\alpha+1)}(1-p^{(s-1)(\ell+1)})}{1-p^{s-1}}\right).
\end{align*}
Similarly, for $I_3$ we have
\begin{align*}
I_3&=-\sum_{i=0}^{\ell}\chi_{-4}(p)^{\frac{(\alpha-2i+1)n}{2}}p^{\frac{n(\alpha-i+1)+(n-2)i-2}{2}-(\alpha-i+1)s}\\
&=-\chi_{-4}(p)^{\frac{n(\alpha+1)}{2}}p^{(\frac{n}{2}-s)(\alpha+1)-1}\sum_{i=0}^{\ell}p^{(s-1)i}=-\frac{\chi_{-4}(p)^{\frac{n(\alpha+1)}{2}}p^{(\frac{n}{2}-s)(\alpha+1)-1}(1-p^{(s-1)(\ell+1)})}{1-p^{s-1}}.
\end{align*}
Collecting all the terms we get the desired formula for $Z_n^{(p)}(s;\bm{\lambda})$ when $n$ is even.
Next, we treat the case when $n$ is odd. In this case we have
\begin{align*}
I_2&=(1-p^{-1})\sum_{i=0}^{\ell}p^{\frac{(n-2)i}{2}}\sum_{\substack{i+1\leq k\leq \alpha-i\\ k\equiv i\Mod{2}}}p^{(\frac{n}{2}-s)k}\\
&=(1-p^{-1})\sum_{i=0}^{\ell}p^{(n-1-s)i}\sum_{k=1}^{\left \lfloor{\frac{\alpha-2i}{2}}\right \rfloor}p^{(n-2s)k}\\
&=\frac{1-p^{-1}}{1-p^{n-2s}}\left(p^{n-2s}\sum_{i=0}^{\ell}p^{(n-1-s)i}-p^{(n-2s)\left \lfloor{\frac{\alpha+2}{2}}\right \rfloor}\sum_{i=0}^{\ell}p^{(s-1)i}\right)\\
&=\frac{1-p^{-1}}{1-p^{n-2s}}\left(\frac{p^{n-2s}(1-p^{(n-1-s)(\ell+1)})}{1-p^{n-1-s}}-\frac{p^{(n-2s)\left \lfloor{\frac{\alpha+2}{2}}\right \rfloor}(1-p^{(s-1)(\ell+1)})}{1-p^{s-1}}\right).
\end{align*} 
For $I_3$, if $\alpha$ is odd using case (ii) we have
\begin{align*}
I_3&=-\sum_{i=0}^{\ell}p^{\frac{n(\alpha-i+1)+(n-2)i-2}{2}-(\alpha-i+1)s}=-\frac{p^{(\frac{n}{2}-s)(\alpha+1)-1}(1-p^{(s-1)(\ell+1)})}{1-p^{s-1}},
\end{align*}
and if $\alpha$ is even, using instead case (iii) {and \rmkref{rmk:dichadsqu}} 
we get that
\begin{align*}
I_3&=\chi_D(p) \sum_{i=0}^{\ell}p^{\frac{n(\alpha-i+1)+(n-2)i-1}{2}-(\alpha-i+1)s}=\chi_D(p) \frac{p^{(\frac{n}{2}-s)(\alpha+1)-\frac{1}{2}}(1-p^{(s-1)(\ell+1)})}{1-p^{s-1}}.
\end{align*}
Collecting all the terms concludes the proof of this case. 
\end{proof}

\begin{remark}\label{rmk:genchi}
Keep the notation and assumptions as in \propref{prop:lozetanonzeroramified}. 
For a general Dirichlet character $\chi$, if $\chi(p)=0$ then $Z_n^{(p)}(s;\chi,\bm{\lambda})=1$. Otherwise, by rewriting $Z_n^{(p)}(s;\chi,\bm{\lambda})=\sum_{k=0}^{\infty}\phi_n(p^k;\bm{\lambda})(\chi(p)p^{-s})^k$
and doing the same computation as above we can get the same formula for $Z_n^{(p)}(s;\chi,\bm{\lambda})$ as in \propref{prop:lozetanonzeroramified} but with $\chi(p)p^{-s}$ in place of $p^{-s}$. In particular, by examining the formulas we see that $Z_n^{(p)}(s;\chi,\bm{\lambda})$ (as a polynomial in $p^{-s}$, cf. \lemref{lem:localzeta}) is divisible by $1-\chi_{-4}(p)^{\frac{n}{4}}\chi(p)p^{\frac{n}{2}-1-s}$ when $n$ is even, by $1-\chi(p)^2p^{n-1-2s}$ when $n$ and $\alpha$ are both odd and by $1+\chi(p)\chi_D(p)p^{\frac{n-1}{2}-s}$ when $n$ is odd and $\alpha$ is even.
\end{remark}

Finally, we treat the case of $ Z_n^{(2)}(s;\bm{\lambda})$, where, while the arguments are very similar, the formulas are more complicated. 
\begin{Prop}\label{prop:Zn2}
{Assume $p=2$}. 
For $\alpha$ odd we have 
  \begin{align*}
Z_n^{(2)}(s; \bm\lambda)&= \frac{1-2^{(n-1-s)(\ell+1)}}{1-2^{n-1-s}}+\delta_{\bm\lambda} 2^{(n-1-s)(\ell+1)}+\frac{2^{n-2s}(2^{-2}a_n+2^{n-3-s}b_n)(1-2^{(n-1-s)\ell})}{(1-2^{n-2s})(1-2^{n-1-s})}\\
&-\frac{2^{(\frac{n}{2}-s)(\alpha-1)}(1-2^{(s-1)\ell})}{1-2^{s-1}}\left(2^{-2}a_n-2^{n-3-s}b_{n-2t}+\frac{2^{-2}a_n+2^{n-3-s}b_n}{1-2^{n-2s}} \right),
\end{align*}
while for $\alpha$ even
\begin{align*}
Z_n^{(2)}(s; \bm\lambda)&= \frac{1-2^{(n-1-s)(\ell+1)}}{1-2^{n-1-s}}+\delta_{\bm\lambda} 2^{(n-1-s)(\ell+1)}+\frac{2^{n-2s}(2^{-2}a_n+2^{n-3-s}b_n)(1-2^{(n-1-s)(\ell-\delta_{\alpha, 2\ell})})}{(1-2^{n-2s})(1-2^{n-1-s})}\\
&+\frac{2^{(\frac{n}{2}-s)\alpha}(1-2^{(s-1)\ell})}{1-2^{s-1}}\left((-1)^{\frac{t+1}{2}}2^{-3}a_{n+2}+2^{n-3-s}b_{n-t}\right)\\
&-\frac{2^{(\frac{n}{2}-s)\alpha}(1-2^{(s-1)(\ell-\delta_{\alpha, 2\ell})})}{1-2^{s-1}}\left(2^{s-2}b_n+\frac{2^{-2}a_n+2^{n-3-s}b_n}{1-2^{n-2s}} \right),
\end{align*}
where $\delta_{\alpha,2\ell}=1$ if $\alpha=2\ell$ and zero otherwise, {and $\delta_{\bm{\lambda}}\in \{0,1\}$ is such that $\delta_{\bm{\lambda}}=1$ if all entries of $(2\bm{\lambda})/2^{\ell}$ are odd and $\delta_{\bm{\lambda}}=0$ otherwise.}
\end{Prop}
\begin{proof}
{For this proof we also use freely \propref{prop:fkifortwo} and the notation therein.}
Using case (v) and the vanishing of $F(k,i)$ for $i> \ell$ and $F(k,\ell)$ for $k\geq \ell+2$ we can start as before with 
\begin{align*}
Z_n^{(2)}(s; \bm\lambda)&=\sum_{i=0}^{\ell-1} \sum_{k=i}^{\alpha-i+1} F(k,i)2^{-sk}+F(\ell,\ell)2^{-\ell s}+ F(\ell+1,\ell)2^{-s(\ell+1)}\\
&=\sum_{i=0}^{\ell-1} \sum_{k=i+1}^{\alpha-i+1} F(k,i)2^{-sk}+\frac{1-2^{(n-1-s)(\ell+1)}}{1-2^{n-1-s}}+\delta_{\bm\lambda} 2^{(n-1-s)(\ell+1)}.
\end{align*}
Here for the second equality we used that $F(i,i)=2^{(n-1)i}$ for all $0\leq i\leq \ell$ and that $F(\ell+1,\ell)=\delta_{\bm{\lambda}}2^{(n-1)(\ell+1)}$.
It thus remains to compute $I:=\sum_{i=0}^{\ell-1} \sum_{k=i+1}^{\alpha-i+1} F(k,i)2^{-sk}$. Rewrite
\begin{align*}
I&=\sum_{i=0}^{\ell-1}\left(\sum_{k=i+1}^{\alpha-i-1}F(k,i)2^{-ks}+F(\alpha-i,i)2^{-(\alpha-i)s}+F(\alpha-i+1,i)2^{-(\alpha-i+1)s}\right)=:I_1+I_2+I_3.
\end{align*}
For $I_2$ applying case (iii) and evaluating geometric sums we get
\begin{align*}
	I_2&=2^{(\frac{n}{2}-s)\alpha}\frac{1-2^{(s-1)\ell}}{1-2^{s-1}}\left\{\begin{array}{ll}
				(-1)^{\frac{t+1}{2}} 2^{-3}a_{n+2}
				& \text{$\alpha$ even},\\
				2^{\frac{n-6}{2}} b_{n-2t} 
				 & \text{$\alpha$ odd}.
			\end{array}\right.
\end{align*}
Here for the $\ell=0$ case we extended the geometric sum formula 
$
\sum_{k=0}^{a}q^k=\frac{1-q^{a+1}}{1-q}$ ($q\neq 1$) for any $a\geq -1$. (When $a=-1$ this is a formal identity with the left hand side void and the right hand side vanishing.)
Similarly, applying case (iv) we get
\begin{align*}
I_3&=2^{(\frac{n}{2}-s)(\alpha+1)}\frac{1-2^{(s-1)\ell}}{1-2^{s-1}}\left\{\begin{array}{ll}
				 2^{\frac{n-6}{2}}b_{n-t}
				& \text{$\alpha$ even},\\
				0
				 & \text{$\alpha$ odd}.
			\end{array}\right.
\end{align*}
For $I_1$ 
we first treat the case when $\alpha=2\beta+1$ is odd. Note that $\beta\geq \ell$. Making a change of variable $k\mapsto k+i$ and splitting the sum into even and odd values of $k$ we get
\begin{align*}
I_1&=\sum_{i=0}^{\ell-1}\left(\sum_{k=1}^{\beta-i}F(2k+i,i)2^{-(2k+i)s}+\sum_{k=1}^{\beta-i}F(2k-1+i, i)2^{-(2k-1+i)s}\right).
\end{align*}
Applying case (i) and (ii) 
we have
\begin{align*}
I_1&=\sum_{i=0}^{\ell-1}\left(\sum_{k=1}^{\beta-i-1}F(2k+i,i)2^{-(2k+i)s}-2^{\frac{n(2\beta-i)+(n-2)i-4}{2}}a_n2^{-(2\beta-i)s}+\sum_{k=2}^{\beta-i}F(2k-1+i, i)2^{-(2k-1+i)s}\right)\\
&=\sum_{i=0}^{\ell-1}\left(\sum_{k=1}^{\beta-i-1}2^{\frac{n(2k+i)+(n-2)i-4}{2}}a_n2^{-(2k+i)s}-a_n2^{(n-2s)\beta-2+i(s-1)}+\sum_{k=2}^{\beta-i}2^{\frac{n(2k+i)+(n-2)i-6}{2}}b_n2^{-(2k-1+i)s}\right)\\
&=\frac{a_n2^{n-2s-2}+b_n2^{2n-3s-3}}{1-2^{n-2s}}\sum_{i=0}^{\ell-1}\left(2^{(n-1-s)i}-2^{(n-2s)(\beta-1)+(s-1)i}\right)-\frac{a_n2^{(n-2s)\beta-2}(1-2^{(s-1)\ell})}{1-2^{s-1}}\\
&=\frac{2^{n-2s}(2^{-2}a_n+2^{n-s-3}b_n)(1-2^{(n-1-s)\ell})}{(1-2^{n-2s})(1-2^{n-1-s})}-\frac{1-2^{(s-1)\ell}}{1-2^{s-1}}\left(\tfrac{2^{(n-2s)\beta}(2^{-2}a_n+2^{n-s-3}b_n)}{1-2^{n-2s}}+a_n2^{(n-2s)\beta-2}\right).
\end{align*}
Collecting all the terms gives the desired formula of $Z_n^{(2)}(s;\bm{\lambda})$ for this case.

Next, we treat the case when $\alpha=2\beta$ is even with $\beta>\ell$. Similarly in this case we can apply case (i) and (ii) 
to compute
\begin{align*}
I_1&=\sum_{i=0}^{\ell-1}\left(\sum_{k=1}^{\beta-i-1}F(2k+i,i)2^{-(2k+i)s}+\sum_{k=1}^{\beta-i}F(2k-1+i, i)2^{-(2k-1+i)s}\right)\\
&=\sum_{i=0}^{\ell-1}\left(\sum_{k=1}^{\beta-i-1}2^{\frac{n(2k+i)+(n-2)i-4}{2}}a_n2^{-(2k+i)s}+\sum_{k=2}^{\beta-i-1}2^{\frac{n(2k+i)+(n-2)i-6}{2}}b_n2^{-(2k-1+i)s}-2^{(n-2s)\beta+i(s-1)+s-3}b_n\right)\\
&=\sum_{i=0}^{\ell-1}\left(\sum_{k=1}^{\beta-i-1}2^{\frac{n(2k+i)+(n-2)i-4}{2}}a_n2^{-(2k+i)s}+\sum_{k=2}^{\beta-i}2^{\frac{n(2k+i)+(n-2)i-6}{2}}b_n2^{-(2k-1+i)s}-2^{(n-2s)\beta+i(s-1)+s-2}b_n\right)\\
&=\frac{2^{n-2s}(2^{-2}a_n+2^{n-s-3}b_n)(1-2^{(n-1-s)\ell})}{(1-2^{n-2s})(1-2^{n-1-s})}-\frac{1-2^{(s-1)\ell}}{1-2^{s-1}}\left(\tfrac{2^{(n-2s)\beta}(2^{-2}a_n+2^{n-s-3}b_n)}{1-2^{n-2s}}+b_n2^{(n-2s)\beta+s-2}\right).
\end{align*}
Collecting all terms we get the desired formula in this case. Finally we treat the case when $\alpha=2\ell$. We note that the only difference between this case and the previous case is that $F(\alpha-i-1,\ell-1)$ is given by different values (cf. case (ii)). 
In this case applying the first equation in case (ii) 
we get
\begin{align*}
I_1&=\sum_{i=0}^{\ell-2}\sum_{k=i+1}^{\alpha-i-1}F(k,i)2^{-ks}+F(\ell, \ell-1)2^{-(\ell-1)s}=\sum_{i=0}^{\ell-2}\sum_{k=i+1}^{\alpha-i-1}F(k,i)2^{-ks}.
\end{align*}
Then doing the same calculation as in the previous case with $\ell-1$ in place of $\ell$ gives us the desired formula for $I_1$ and hence also for $Z_n^{(2)}(s;\bm{\lambda})$ in this case.
\end{proof}

We actually need to compute the term $\varphi_n(2;2\bm{\lambda})+\tilde{Z}_n^{(2)}(s;\bm{\lambda})=(-1)^{2\lambda_1}2^{n-1}+\tilde{Z}_n^{(2)}(s;\bm{\lambda})$. 
Define
\begin{align}\label{def:z2}
\cZ_n^{(2)}(s;\bm{\lambda})&:=(-1)^{2\lambda_1}+2^{1-n}\tilde{Z}_n^{(2)}(s;\bm\lambda)
\end{align}
so that $\varphi_n(2;2\bm{\lambda})+\tilde{Z}_n^{(2)}(s;\bm{\lambda})=2^{n-1}\cZ_n^{(2)}(s;\bm{\lambda})$, or equivalently (noting also that $\phi_n(2;2\bm{\lambda})=2^{n-1}$)
\begin{align}\label{equ:czzrela}
\cZ_n^{(2)}(s;\bm{\lambda})=2^{1-n+s}\left(Z_n^{(2)}(s;\bm{\lambda})-1\right)+(-1)^{2\lambda_1}-1.
\end{align}
Applying \propref{prop:Zn2} together with the formulas for $a_n$ and $b_n$ (see \rmkref{rmk:abn}) we get the following formulas for $\cZ_n^{(2)}(s;\bm{\lambda})$ for which we omit the proof.
\begin{Prop}\label{prop:efactocomp2}
Keep the notation {and assumptions} as in \propref{prop:Zn2}. When $\ell=0$ we have that $\cZ_n^{(2)}(s;\bm{\lambda})=-1$.
When  $\ell\geq 1$, 
for $n\equiv 2\Mod{4}$,
\begin{align*}
\cZ_n^{(2)}(s;\bm{\lambda})&=\frac{1-2^{(n-1-s)\ell}}{1-2^{n-1-s}}+\delta_{\bm{\lambda}}2^{(n-1-s)\ell}+\frac{(-1)^{\frac{n-2t}{4}}2^{(\frac{n}{2}-s)(\alpha-1)}(1-2^{(s-1)\ell})}{1-2^{s-1}}.
\end{align*}
For $n\equiv 0\Mod{4}$,
\begin{align*} 
\cZ_n^{(2)}(s;\bm{\lambda})&=\frac{1-2^{(n-1-s)\ell}}{1-2^{n-1-s}}+\delta_{\bm\lambda} 2^{(n-1-s)\ell}+\frac{(-1)^{\frac{n}{4}}2^{\frac{n}{2}-s}(1-2^{(n-1-s)(\ell-\delta_{\alpha,2\ell})})}{(1-2^{\frac{n}{2}-s})(1-2^{n-1-s})}\\
&-\frac{(-1)^{\frac{n}{4}}2^{(\frac{n}{2}-s)(\alpha-2)+1}(1-2^{\frac{n}{2}-1-s})(1-2^{(s-1)(\ell-\delta_{\alpha,2\ell})})}{(1-2^{\frac{n}{2}-s})(1-2^{s-1})}.
\end{align*}
When $n$ is odd we have for $\alpha$ odd,
\begin{align*}
\cZ_n^{(2)}(s;\bm{\lambda})&=\frac{(1+\mathfrak{s}_n2^{\frac{n+1}{2}-s})(1-\mathfrak{s}_n2^{\frac{n-1}{2}-s})(1-2^{(n-1-s)\ell})}{(1-2^{n-2s})(1-2^{n-1-s})}-\frac{\mathfrak{s}_n2^{(\frac{n}{2}-s)(\alpha-2)+\frac12}(1-2^{n-1-2s})(1-2^{(s-1)\ell})}{(1-2^{n-2s})(1-2^{s-1})},
\end{align*}
while for $\alpha$ even,
\begin{align*}
\cZ_n^{(2)}(s; \bm\lambda)&= \frac{1-2^{(n-1-s)\ell}}{1-2^{n-1-s}}+\delta_{\bm\lambda} 2^{(n-1-s)\ell}+\frac{\mathfrak{s}_n2^{\frac{n-1}{2}-s}(1-2^{(n-1-s)(\ell-\delta_{\alpha, 2\ell})})}{(1-2^{n-2s})(1-2^{n-1-s})}\\
&+\tfrac{2^{(\frac{n}{2}-s)(\alpha-1)-\frac{n}{2}+1}(1-2^{(s-1)\ell})}{1-2^{s-1}}\left(2^{n-3-s}b_{n-t}-(-1)^{\frac{n+t}{2}}\mathfrak{s}_n 2^{\frac{n-3}{2}}\right)-\tfrac{\mathfrak{s}_n2^{(\frac{n}{2}-s)(\alpha-1)-\frac{1}{2}}(1-2^{(s-1)(\ell-\delta_{\alpha, 2\ell})})}{(1-2^{n-2s})(1-2^{s-1})}.
\end{align*}
Here recall that when $n$ is odd, $\mathfrak{s}_n=(-1)^{\frac{n^2-1}{8}}$. \qed
\end{Prop}

\subsection{Putting it all together}\label{sec:piat}
Combining our preliminary formulas together with the Euler product expansion and the calculation of the local zeta functions we get our Fourier expansion formulas. 
\begin{proof}[Proof of \thmref{thm:fouexpan}]
We first prove the constant term formula. Combining the Euler product formula proved in \propref{prop:eulprodconst} and the local zeta function formulas in \propref{prop:loczetconst} and the definitions of $\e_n^{(p)}(s)$ we have
	\begin{align*}
	Z_{n,d}(s;\bm{0})&=\left\lbrace\begin{array}{cc}
2^{\frac{3}{2}n-s}\prod_{p\mid 2d}\e_n^{(p)}(s)\frac{\zeta(s-n+1)\zeta(s-\frac{n}{2})}{\zeta(s-\frac{n}{2}+1)} & n\equiv 0\Mod{4},\\
2^{\frac{3}{2}n-s-1}\prod_{p\mid 2d}\e_n^{(p)}(s)\frac{\zeta(s-n+1)L(s-\frac{n}{2},\chi_{-4})}{L(s-\frac{n}{2}+1,\chi_{-4})} & n\equiv 2\Mod{4},\\
2^{\frac{3n-1}{2}-s}\prod_{p\mid 2d}\e_n^{(p)}(s)\frac{\zeta(2s-n)\zeta(s-n+1)}{\zeta(2s-n+1)} & n\equiv 1\Mod{2}. \end{array}\right.
	\end{align*}
Plugging this into the relation (cf. \propref{prop:preid})
\begin{align*}
\Phi_{n,d}(s)&=2^{s-n}d^{s-n}\pi^{\frac{n}{2}} \frac{\G\left(\frac{2s-n}{2}\right)}{\zeta(s)\G(s)}Z_{n,d}(s;\bm{0}),
\end{align*} 
using the definition $\e_{n,d}(s)=d^{s-n}\prod_{p\mid 2d}\e_{n}^{(p)}(s)$ and applying the relations $\zeta(s)=\pi^{s/2}\G\left(\tfrac{s}{2}\right)^{-1}\xi(s)$  and $L(s,\chi_{-4})=\left(\tfrac{\pi}{4}\right)^{\frac{s+1}{2}} \G\left(\tfrac{s+1}{2}\right)^{-1}L^*(s,\chi_{-4})$ to complete the zeta and $L$-function (the latter is only needed for the case of $n\equiv 2\Mod{4}$) and the relation 
\begin{align*}
\frac{\G(s-m)}{\G(s)}=2^{-m}\frac{\G(\frac{s-m}{2})\G(\frac{s-m+1}{2})}{\G(\frac{s}{2})\G(\frac{s+1}{2})},\qquad\forall\ m\in\Z_{\geq 0}
\end{align*}
to further simplify the $\G$-factors we get the desired formulas for $\Phi_{n,d}(s)$.

Next, for non-constant terms fix $\bm{\lambda}\in \Lambda^{*}\setminus\{\bm{0}\}$. Combining
\propref{prop:nonconsprform} and \lemref{lem:lozetanonzero}
we have that $Z_{n,d}(s;\bm{\lambda})$ equals
\begin{align*}
&\sum_{d_1|d} \frac{1}{\vf(d_1)}\sum_{\chi\Mod{d_1}} \bigg(\tilde{Z}_n^{(2)}(s;\chi, \bm\lambda)+ (-1)^{2\lambda_1}2^{n-1}\bigg)\overline\chi(\tfrac{2d}{d_1}) 
\cS(d_1,\chi, \bm\lambda)\prod_{p|\frac{d}{d_1}}\tilde{Z}_n^{(p)}(s;\chi, \bm\lambda)\\
&\times \prod_{\substack{(p,2d)=1\\ p| \|2\bm\lambda\|^2 }}Z^{(p)}_n(s; \chi, \bm\lambda)\prod_{\substack{(p,2d)=1\\ p\nmid \|2\bm\lambda\|^2 }}\left\lbrace \begin{array}{ll}
1-\chi_{-4}(p)^{\frac{n}{2}}\chi(p)p^{\frac{n}{2}-1 -s} & \text{$n$ even},\\ 
1+\chi(p)\chi_{D}(p)p^{\frac{n-1}{2} -s} & \text{$n$ odd},
\end{array}\right. 
 \end{align*}
 with $D=(-1)^{\frac{n-1}{2}}\|2\bm{\lambda}\|^2$ when $n$ is odd as before. 
Let $N=2d\|2\bm\lambda\|^2$ and note that since for any odd prime $p$, $\chi_{-4}(p)^{\frac{n}{2}}=\left\lbrace \begin{array}{ll}
1 & n\equiv 0\Mod{4},\\ 
\chi_{-4}(p) & n\equiv 2\Mod{4},
\end{array}\right. 
$ the above last product factor equals
\begin{align*}
\left\lbrace \begin{array}{ll}
\prod_{p|N}(1-\chi(p)p^{\frac{n}{2}-1 -s})^{-1} \times \frac{1}{L(s-\frac{n}{2}+1,\chi)} & n\equiv 0\Mod{4},\\ 
\prod_{p|N}(1-\chi\chi_{-4}(p)p^{\frac{n-2}{2} -s})^{-1}\times \frac{1}{L(s-\frac{n}{2}+1,\chi\chi_{-4})} & n\equiv 2\Mod{4},\\ 
\prod_{p|N}\frac{1-\chi\chi_{D}(p)p^{\frac{n-1}{2} -s}}{1-\chi(p)^2p^{n-1-2s}}\times \frac{L(s-\frac{n-1}{2},\chi\chi_{D})}{L(2s-n+1,\chi^2)} & n\equiv 1\Mod{2}.
\end{array}\right. 
\end{align*}
Here when $n$ is odd, we used the identity that  
$$
1+\chi(p)\chi_{D}(p)p^{\frac{n-1}{2}-s}=\frac{1-\chi^2(p)p^{n-1-2s}}{1-\chi(p)\chi_{D}(p)p^{\frac{n-1}{2}-s}}
$$ 
for $p$ not dividing $N$ (since for such $p$, $\chi_{D}(p)^2=1$).
Thus in view of \propref{prop:preid} for any $d_1\mid d$ and any Dirichlet character $\chi$ modulo $d_1$ if we define 
\begin{align}\label{def:epsilonfac}
\e_n(s; \chi, \bm{\lambda})&:=\left((-1)^{2\lambda_1}+ 2^{1-n}\tilde{Z}_n^{(2)}(s;\chi, \bm\lambda) \right)\overline\chi(\tfrac{2d}{d_1}) 
\cS(d_1,\chi, 2\bm\lambda)\prod_{p|\frac{d}{d_1}}\tilde{Z}_n^{(p)}(s;\chi, \bm\lambda)\nonumber\\
&\times \left\lbrace \begin{array}{ll}
\prod_{\substack{(p,2d)=1\\ p| \|2\bm\lambda\|^2}} \frac{Z^{(p)}_n(s; \chi,\bm\lambda)}{(1-\chi(p)p^{\frac{n}{2}-1 -s})}\prod_{p\mid 2d}(1-\chi(p)p^{\frac{n}{2}-1 -s})^{-1}  & n\equiv 0\Mod{4},\\ 
\prod_{\substack{(p,2d)=1\\ p| \|2\bm\lambda\|^2}}\frac{ Z^{(p)}_n(s; \chi,\bm\lambda)}{(1-\chi\chi_{-4}(p)p^{\frac{n-2}{2} -s})}\prod_{p|2d}(1-\chi\chi_{-4}(p)p^{\frac{n-2}{2} -s})^{-1} & n\equiv 2\Mod{4},\\ 
\prod_{\substack{(p,2d)=1\\ p| \|2\bm\lambda\|^2}}\frac{(1-\chi\chi_{D}(p)p^{\frac{n-1}{2} -s})Z^{(p)}_n(s; \chi,\bm\lambda)}{1-\chi(p)^2p^{n-1-2s}} \prod_{p|2d}\frac{1-\chi\chi_{D}(p)p^{\frac{n-1}{2} -s}}{1-\chi(p)^2p^{n-1-2s}} & n\equiv 1\Mod{2},
\end{array}\right.  
\end{align}
then we have $\Phi_{n,d}(s;\bm{\lambda}):=2^{1-n}Z_{n,d}(s;\bm{\lambda})$ satisfies the desired formula \eqref{equ:phinonconst}.
Moreover, by \lemref{lem:localzeta} and \rmkref{rmk:genchi} each of the finitely many local zeta functions $\tilde{Z}_n^{(p)}(s;\chi, \bm\lambda)$ and all the above factors with $(p,2d)=1$ appearing in the above definition are polynomials in $p^{-s}$, and hence holomorphic everywhere.
Since the factors with $p\mid 2d$ are also holomorphic off the line $\Re(s)=\left \lfloor{\frac{n-1}{2}}\right \rfloor$ we get that $\e_n(s;\chi, \bm{\lambda})$ is holomorphic off this line. 
Moreover, the poles are coming from these finitely many prime factors $p\mid 2d$ which {lie} in a finite union of periodic points on the line $\Re(s)=\left \lfloor{\frac{n-1}{2}}\right \rfloor$. 
The statement regarding the bound for $|\e_n(s;\chi, \bm{\lambda})|$ in the half space $\Re(s)\geq \frac{n}{2}$ follows easily from  \lemref{lem:localzeta}.

Finally we prove the vanishing statement. Assume $n$ is odd and $n\not\equiv 1\Mod{8}$ and fix a nonzero $\bm{\lambda}\in \Lambda^*$. In view of the above definition of $\e_n(s;\chi,\bm{\lambda})$ it suffices to show that for any Dirichlet character $\chi$ modulo $d_1$ for some $d_1\mid d$, the product
$$
\left((-1)^{2\lambda_1}+2^{1-n}\tilde{Z}_n^{(2)}(s;\chi, \bm\lambda)\right)
\cS(d_1,\chi, 2\bm\lambda)
$$ 
vanishes at $s=\frac{n+1}{2}$ whenever $\chi\chi_{D}$ is principal. First if $\chi$ is not principal, then by \propref{prop: vfnp} we have that $\cS(d_1,\chi, 2\bm\lambda)=0$. 
If $\chi$ is principal, since $d_1$ is odd, $\chi(2^k)=1$ for all $k\in\N$. Hence $\tilde{Z}_n^{(2)}(s;\chi, \bm\lambda)=\tilde{Z}_n^{(2)}(s; \bm\lambda)$ and in view of the relation \eqref{def:z2}, the first factor in the above product equals $\cZ_n^{(2)}(s;\bm{\lambda})$. It thus suffices to show $\cZ_n^{(2)}(\tfrac{n+1}{2};\bm{\lambda})=0$. {For the remainder of the proof we fix the prime $p=2$ and omit the subscripts from $\alpha=\nu_2(\|2\bm{\lambda}\|^2)$, $\ell=\nu_2(\gcd(2\bm{\lambda}))$ and $t=\|2\bm{\lambda}\|^2/2^{\alpha}$}. Now,  since $\chi\chi_D$ is principal, $\chi_{D}$ is also principal, but this is possible only when $n\equiv 5\Mod{8}$ and $\|2\bm{\lambda}\|^2$ is a square, and thus 
$\alpha$ is even. In particular, we must have entries of $2\bm{\lambda}$ are all even since otherwise 
we would have $\|2\bm{\lambda}\|^2\equiv n\equiv 5\Mod{8}$, but $5 \Mod{8}$ is not a square residue. Hence $\ell\geq 1$. Moreover, we also have $t$ is an odd square, implying that $t\equiv 1\Mod{8}$ (so that $t\equiv n\Mod{4}$ but $t\not\equiv n\Mod{8}$). This then implies that $\delta_{\bm{\lambda}}=0$ since otherwise we would have $\alpha=2\ell$ (since $n$ is odd) and $t=\|2\bm{\lambda}\|^2/2^{2\ell}\equiv n\equiv 5\Mod{8}$, contradicting the fact that $t$ is an odd square. Hence by \propref{prop:efactocomp2} we have (using also that in this case $\mathfrak{s}_n=(-1)^{\frac{n^2-1}{8}}=-1=(-1)^{\frac{n+t}{{2}}}$ and $b_{n-t}=-4$) 
\begin{align*}
\cZ_n^{(2)}(s; \bm\lambda)&= \frac{1-2^{(n-1-s)\ell}}{1-2^{n-1-s}}-\frac{2^{\frac{n-1}{2}-s}(1-2^{(n-1-s)(\ell-\delta_{\alpha, 2\ell})})}{(1-2^{n-2s})(1-2^{n-1-s})}\\
&-\frac{2^{(\frac{n}{2}-s)(\alpha-1)-\frac{1}{2}}(1+2^{\frac{n+1}{2}-s})(1-2^{(s-1)\ell})}{1-2^{s-1}}+\frac{2^{(\frac{n}{2}-s)(\alpha-1)-\frac{1}{2}}(1-2^{(s-1)(\ell-\delta_{\alpha, 2\ell})})}{(1-2^{n-2s})(1-2^{s-1})}.
\end{align*}
By rewriting $\frac{1-q^{\ell}}{1-q}=\frac{1-q^{\ell-\delta_{\alpha,2\ell}}}{1-q}+\delta_{\alpha,2\ell}q^{\ell-1}$ for $q=2^{s-1},2^{n-1-s}$ and using the identity $1-2^{n-2s}-2^{\frac{n-1}{2}-s}=(1-2^{\frac{n+1}{2}-s})(1+2^{\frac{n-1}{2}-s})$ we can further simplify the above expression to get
\begin{align*}
\cZ_n^{(2)}(s;\bm{\lambda})&=\tfrac{(1-2^{\frac{n+1}{2}-s})(1+2^{\frac{n-1}{2}-s})(1-2^{(n-1-s)(\ell-\delta_{\alpha,2\ell})})}{(1-2^{n-2s})(1-2^{n-1-s})}-\tfrac{2^{(\frac{n}{2}-s)\alpha}(1-2^{\frac{n+1}{2}-s})(1+2^{\frac{n-1}{2}-s})(1-2^{(s-1)(\ell-\delta_{\alpha,2\ell})})}{(1-2^{n-2s})(1-2^{s-1})}\\
&+\delta_{\alpha,2\ell}2^{(n-1-s)(\ell-1)}(1- 2^{\frac{n+1}{2}-s})(1+ 2^{\frac{n-1}{2}-s}),
\end{align*}
from which one easily sees that $\cZ_n^{(2)}(\tfrac{n+1}{2};\bm{\lambda})=0$. This finishes the proof.
\end{proof}

The rest of this subsection is devoted to proving \thmref{thm:funequbofo} which, as illustrated before, implies the functional equation \eqref{e:funeq}. For the remaining of this subsection, we assume $d=1$ and fix a nonzero $\bm{\lambda}\in \Lambda^*$. 
As before we abbreviate $\e_n(s;\chi,\bm{\lambda})$ by $\e_n(s;\bm{\lambda})$ when $\chi$ is the trivial character. Then the definition \eqref{def:epsilonfac} reads as
\begin{align}\label{equ:efadeq1case}
\e_n(s;\bm{\lambda})&=\prod_{p\mid 2\|2\bm{\lambda}\|^2}\e_n^{(p)}(s;\bm{\lambda}),
\end{align}
where for any odd prime $p\mid \|2\bm{\lambda}\|^2$,
\begin{align}\label{equ:relazep}
\e_n^{(p)}(s;\bm{\lambda}):=\left\{\begin{array}{ll}
				(1-\chi_{-4}(p)^{\frac{n}{2}}p^{\frac{n}{2}-1-s})^{-1}Z^{(p)}_n(s;\bm{\lambda}) & n\equiv 0\Mod{2},\\
			\frac{1-\chi_{D}(p)p^{\frac{n-1}{2}-s}}{1-p^{n-1-2s}}Z_n^{(p)}(s;\bm{\lambda})	 & n\equiv 1\Mod{2},
			\end{array}\right.
\end{align}
and for $p=2$,
\begin{align}\label{equ:relazep2}
\e_n^{(2)}(s;\bm{\lambda})&:=\left\{\begin{array}{ll}
				(1-2^{\frac{n}{2}-1-s})^{-1}\cZ_n^{(2)}(s;\bm{\lambda}) & n\equiv 0\Mod{4},\\
				\cZ_n^{(2)}(s;\bm{\lambda}) & n\equiv 2\Mod{4},\\
			\frac{1-\chi_{D}(2)2^{\frac{n-1}{2}-s}}{1-2^{n-1-2s}}\cZ_n^{(2)}(s;\bm{\lambda})	 & n\equiv 1\Mod{2}.
			\end{array}\right.
\end{align}
Here $Z_n^{(p)}(s;\bm{\lambda})$ and $\cZ_n^{(2)}(s;\bm{\lambda})$ are as in \propref{prop:lozetanonzeroramified} and \propref{prop:efactocomp2} respectively and when $n$ is odd $D=(-1)^{\frac{n-1}{2}}\|2\bm{\lambda}\|^2$ as before. 
In order to make \eqref{equ:relazep2} more explicit, we need to compute $\chi_{D}(2)$. 

\begin{Lem}\label{lem:relaprimchar}
	Assume $n$ is odd. 
	Then 
	\begin{align*}
		\chi_{D}(2)&=\left\{\begin{array}{ll}
			(-1)^{\frac{n-t_2}{4}}\mathfrak{s}_n & \text{$\alpha_2$ even and}\ t_2\equiv n\Mod{4},\\
			0 & \text{otherwise}.
		\end{array}\right.
	\end{align*}
			\end{Lem}
\begin{proof}
	 {We first give a more precise description of the primitive quadratic character $\chi_D$.} Let $D_0\in \Z$ be the {unique} odd, square-free integer such that $D_0$ is of the same sign of $D$ and $t_2/|D_0|$ is a square. {Then we have \begin{align}\label{equ:prchardes}
\chi_D(\cdot)&=\left\{\begin{array}{ll}   
				\left(\tfrac{8D_0}{\cdot}\right) & \text{$\alpha_2$ odd},\\
				\left(\tfrac{D_0}{\cdot}\right) &  \text{$\alpha_2$ even}\ \text{and}\ D_0\equiv 1\Mod{4},\\
				\left(\tfrac{4D_0}{\cdot}\right)  &  \text{$\alpha_2$ even}\ \text{and}\ D_0\equiv 3\Mod{4}.
			\end{array}\right.
\end{align}}
We also note that $(-1)^{\frac{n-1}{2}}t_2/D_0$ is an odd square, {thus $(-1)^{\frac{n-1}{2}}t_2/D_0\equiv 1\Mod{8}$, or equivalently, $D_0\equiv (-1)^{\frac{n-1}{2}}t_2\Mod{8}$}. 

We now prove the above formulas. We first prove the second case which splits into two cases that either $\alpha_2$ is odd or $\alpha_2$ is even and $t_2\not\equiv n\Mod{4}$. For the first case we have by \eqref{equ:prchardes}, $\chi_{D}(2)=(\tfrac{8D_0}{2})=0$. For the second case, note that since $n$ is odd, $(-1)^{\frac{n-1}{2}}\equiv n\Mod{4}$. Thus $D_0\equiv (-1)^{\frac{n-1}{2}}t_2\equiv nt_2\equiv 3\Mod{4}$, where for the last equality we used the assumptions that $t_2\not\equiv n\Mod{4}$ and both $t_2$ and $n$ are odd.  Thus by \eqref{equ:prchardes} we have $\chi_D(2)=(\frac{4D_0}{2})=0$. 

Finally assume $\alpha_2$ is even and $t_2\equiv n\Mod{4}$. Then $D_0\equiv 1\Mod{4}$ and thus 
$$
\chi_D(2)=\left(\frac{D_0}{2}\right)=\left(\tfrac{(-1)^{\frac{n-1}{2}}t_2}{2}\right)=\left(\frac{t_2}{2}\right)=(-1)^{\frac{t_2^2-1}{8}}
$$
from which the first case of this lemma follows immediately. 
	\end{proof}
	\begin{remark}\label{rmk:qformula}
	Keep the notation and assumptions as above. In view of the above lemma,  the last case in \eqref{equ:relazep2} (when $n$ is odd) reads as
	\begin{align*}
		\e_n^{(2)}(s;\bm{\lambda})&=\left\{\begin{array}{ll}
			(1+\mathfrak{s}_n 2^{\frac{n-1}{2}-s})^{-1}\cZ_n^{(2)}(s;\bm{\lambda}) & \text{$\alpha_2$ even and}\ t_2\equiv n\Mod{8},\\
			(1-\mathfrak{s}_n 2^{\frac{n-1}{2}-s})^{-1}\cZ_n^{(2)}(s;\bm{\lambda}) & \text{$\alpha_2$ even and}\ t_2\equiv n+4\Mod{8},\\
			(1-2^{n-1-2s})^{-1}\cZ_n^{(2)}(s;\bm{\lambda}) & \text{otherwise}.
		\end{array}\right.
	\end{align*}
On the other hand, let $D=(-1)^{\frac{n-1}{2}}\|2\bm{\lambda}\|^2$ and $D_0$ be as above, and let $q$ be the modulus of $\chi_{D}$. Note that the numerators in the right hand side of \eqref{equ:prchardes} are \textit{quadratic discriminants} (see \cite[p. 296]{MontgomeryVaughan2007} for the definition) and thus the modulus of $\chi_D$ equals the absolute value of these numerators \cite[Theorem 9.13]{MontgomeryVaughan2007}. This fact, together with the analysis in the above proof gives that  
	\begin{align}\label{equ:primcharmodu}
		q&=\left\{\begin{array}{ll}   
			8|D_0| & \text{$\alpha_2$ odd},\\
			4|D_0| &  \text{$\alpha_2$ even\ \text{and}\ $t_2\not\equiv n\Mod{4}$},\\
			|D_0| &  \text{$\alpha_2$ even\ \text{and}\ $t_2\equiv n\Mod{4}$}.
		\end{array}\right.
	\end{align}
\end{remark}

\begin{remark}\label{rmk:polyargu}
Let $D=(-1)^{\frac{n-1}{2}}\|2\bm{\lambda}\|^2$ be as before. For any odd prime $p\mid D$ 
{as discussed in the proof of \thmref{thm:fouexpan} we have} that $\e_n^{(p)}(s;\bm{\lambda})$ is  a polynomial in $p^{-s}$ (noting that $\chi_D(p)=0$ if $\alpha_p=\nu_p(|D|)$ is odd and $\chi_D(p)\in \{\pm 1\}$ if $\alpha_p$ is even). 
Similarly, for $p=2$ we see from the definition \eqref{equ:relazep2} that
\begin{align*}
\left\{\begin{array}{ll}   
			\e_n^{(2)}(s;\bm{\lambda})=\cZ_n^{(2)}(s;\bm{\lambda}) & n\equiv 2\Mod{4},\\
			(1-2^{\frac{n}{2}-1-s})\e_n^{(2)}(s;\bm{\lambda})=\cZ_n^{(2)}(s;\bm{\lambda}) &  n\equiv 0\Mod{4},\\
			(1-2^{n-1-2s})\e_n^{(2)}(s;\bm{\lambda})=(1-\chi_D(2)2^{\frac{n-1}{2}-s})\cZ_n^{(2)}(s;\bm{\lambda}) &  n\equiv 1\Mod{2},
		\end{array}\right.
\end{align*}
is a polynomial in $2^{-s}$; see also \eqref{def:z2} and \lemref{lem:localzeta}. 
\end{remark}

\begin{remark}\label{rmk:efacaritfun}
When $n$ is even our formula for the odd prime local factors $\e_n^{(p)}(s;\bm{\lambda})$ can be rewritten  in terms of twisted divisor functions as follows
\begin{align*}
\e_n^{(p)}(s;\bm{\lambda})
&= \sum_{j=0}^{\ell_p}p^{(n-1-s)j}\tau_{\frac{n}{2}-s}(p^{\alpha_p-2j}; \chi_{-4}^{n/2}),
\end{align*}
where $\tau_{s}(m; \chi)=\sum_{d\mid m}\chi(d)d^s$.
In particular, if we set $a=\gcd(2\bm{\lambda})/2^{\ell_2}$ and $b=\|2\bm{\lambda}\|^2/2^{\alpha_2}$, then 
\begin{equation}\label{e:epsilonasdivisor}
\e_n(s;\bm{\lambda})=\e_n^{(2)}(s;\bm{\lambda})\sum_{d\mid a}d^{n-1-s}\tau_{\frac{n}{2}-s}\left(\tfrac{b}{d^2};\chi_{-4}^{n/2}\right).
\end{equation}
For $n=2$ we also get a simple formula for $\e_n^{(2)}(s;\bm{\lambda})$ that $\e_2^{(2)}(s;\bm{\lambda})=\left\{\begin{array}{ll}   
			-1 & \alpha_2=0,\\
			\tau_{1-s}(2^{\alpha_2-1}) &  \alpha_2>0.
		\end{array}\right.$
For $n=1$, we have that $2\bm{\lambda}=2^{\ell_2} a$  with $a\in \Z$ odd and $\epsilon_1(s;\bm\lambda)=\epsilon_1^{(2)}(s;\bm\lambda)\tau_{1-2s}(|a|)$, and while the formula for $\epsilon_1^{(2)}(s;\bm\lambda)$ can be written explicitly it is no longer a simple divisor function.
For larger odd $n>1$, the local factors for odd primes $p$ co-prime to the conductor of $\chi_D$ also have a similar simple expression as divisor functions, but for the prime divisors of the conductor (as well as for $p=2$) the formula is again more complicated.
\end{remark}

We are now ready to describe the functional equation satisfied by these local $\e$-factors $\e_n^{(p)}(s;\bm{\lambda})$ which eventually leads to the functional equation for $\e_n(s;\bm{\lambda})$. 
\begin{Prop}\label{prop:lononzeconterop}
	For any odd prime $p\mid \|2\bm{\lambda}\|^2$, the $\e$-factor $\e_n^{(p)}(s;\bm{\lambda})$ satisfies the functional equation that for $n$ even,
\begin{align*}
\e_n^{(p)}(n-s;\bm{\lambda})&=\chi_{-4}(p)^{\frac{n\alpha_p}{2}}p^{(s-\frac{n}{2})\alpha_p}\e_n^{(p)}(s;\bm{\lambda}),
\end{align*}
and for $n$ odd,
\begin{align*}
\e_n^{(p)}(n-s;\bm{\lambda})&=\left\{\begin{array}{ll}                               
							p^{(s-\frac{n}{2})(\alpha_p-1)}\e_n^{(p)}(s;\bm{\lambda})
&  \text{$\alpha_p$ odd},\\
			p^{(s-\frac{n}{2})\alpha_p}\e_n^{(p)}(s;\bm{\lambda})	  & \text{$\alpha_p$ even}.
			\end{array}\right.
\end{align*}
For $p=2$, the $\e$-factor $\e_n^{(2)}(s;\bm{\lambda})$ satisfies the functional equation that for $n$ even,
\begin{align*}
	\e_n^{(2)}(n-s;\bm{\lambda})&=
	\left\{\begin{array}{ll}
		\frac{2^{(s-\frac{n}{2})(\alpha_2-2)}(1-2^{\frac{n}{2}-1-s})}{1-2^{s-\frac{n}{2}-1}}\e_n^{(2)}(s;\bm{\lambda}) & n\equiv 4\Mod{8},\\
		(-1)^{\frac{n-2t_2}{4}}2^{(s-\frac{n}{2})(\alpha_2-1)}\e_n^{(2)}(s;\bm{\lambda}) & n\equiv 2\Mod{4},
	\end{array}\right.
\end{align*}
while for $n$ odd, if $\alpha_2$ is odd, then
\begin{align*}
	\e_n^{(2)}(n-s;\bm{\lambda})&=\frac{\mathfrak{s}_n2^{(s-\frac{n}{2})(\alpha_2-4)+\frac12}(1+\mathfrak{s}_n2^{\frac{n-1}{2}-s})}{1+\mathfrak{s}_n 2^{\frac{n+1}{2}-s}}\e_n^{(2)}(s;\bm{\lambda}),
\end{align*}
and if $\alpha_2$ is even, then
\begin{align*}
	\e_n^{(2)}(n-s;\bm{\lambda})&=
	\left\{\begin{array}{ll}
		\frac{\mathfrak{s}_n2^{(s-\frac{n}{2})(\alpha_2-1)+\frac12}(1+\mathfrak{s}_n2^{\frac{n-1}{2}-s})}{1+\mathfrak{s}_n2^{\frac{n+1}{2}-s}}\e^{(2)}_n(s;\bm{\lambda}) & t_2\equiv n\Mod{4},\\
		\frac{\mathfrak{s}_n2^{(s-\frac{n}{2})(\alpha_2-3)+\frac12}(1+\mathfrak{s}_n2^{\frac{n-1}{2}-s})}{1+\mathfrak{s}_n 2^{\frac{n+1}{2}-s}}\e_n^{(2)}(s;\bm{\lambda}) & t_2\not\equiv n\Mod{4}.
	\end{array}\right.
\end{align*}
\end{Prop}

\begin{proof}
	First we treat the case when $p$ is odd.
Combining \propref{prop:lozetanonzeroramified} and the relation \eqref{equ:relazep} we immediately get that for $n$ even,
	\begin{align*}
\e_n^{(p)}(s;\bm{\lambda})&=\frac{1}{1-\chi_{-4}(p)^{\frac{n}{2}}p^{\frac{n}{2}-s}}\left(\frac{1-p^{(n-1-s)(\ell_p+1)}}{1-p^{n-1-s}}-\frac{\chi_{-4}(p)^{\frac{n(\alpha_p+1)}{2}}p^{(\frac{n}{2}-s)(\alpha_p+1)}(1-p^{(s-1)(\ell_p+1)})}{1-p^{s-1}}\right),	
\end{align*}
	and for $n$ and $\alpha_p$ both odd,
	\begin{align*}
\e_n^{(p)}(s;\bm{\lambda})&=\frac{1-p^{(n-1-s)(\ell_p+1)}}{(1-p^{n-2s})(1-p^{n-1-s})}-\frac{p^{(\frac{n}{2}-s)(\alpha_p+1)}(1-p^{(s-1)(\ell_p+1)})}{(1-p^{n-2s})(1-p^{s-1})},
\end{align*}
	while for $n$ odd and $\alpha_p$ even,
\begin{align*}
\e_n^{(p)}(s;\bm{\lambda})&=\tfrac{(1-\chi_{D}(p)p^{\frac{n-1}{2}-s})(1-p^{(n-1-s)(\ell_p+1)})}{(1-p^{n-2s})(1-p^{n-1-s})}+\tfrac{\chi_{D}(p)p^{(\frac{n}{2}-s)(\alpha_p+1)-\frac12}(1-\chi_{D}(p)p^{\frac{n+1}{2}-s})(1-p^{(s-1)(\ell_p+1)})}{(1-p^{n-2s})(1-p^{s-1})}.
\end{align*}
From these formulas one easily verifies the desired functional equations for $\e_n^{(p)}(s;\bm{\lambda})$.

Next, we treat the case when $p=2$ and we assume $n\not\equiv 0\Mod{8}$. The case when $\ell_2=0$ can be checked directly using the definition \eqref{equ:relazep2}, \rmkref{rmk:qformula} and the facts that in this case $\cZ_n^{(2)}(s;\bm{\lambda})=-1$ (see \propref{prop:efactocomp2}) and $\alpha_2=\nu_2(n)$. 

For the remaining we thus assume $\ell_2\geq 1$.
We first prove the case when $n\equiv 2\Mod{4}$. In this case by \propref{prop:efactocomp2} and \eqref{equ:relazep2} we have
\begin{align*}
	\e_n^{(2)}(s;\bm{\lambda})=\frac{1-2^{(n-1-s)\ell_2}}{1-2^{n-1-s}}+\frac{(-1)^{\frac{n-2t_2}{4}}2^{(\frac{n}{2}-s)(\alpha_2-1)}(1-2^{(s-1)\ell_2})}{1-2^{s-1}}+\delta_{\bm{\lambda}}2^{(n-1-s)\ell_2}.
\end{align*}
Let $A_1(s)$ be the sum of the above first two terms and let $A_2(s):=\delta_{\bm{\lambda}}2^{(n-1-s)\ell_2}$ so that $\e_n^{(2)}(s;\bm{\lambda})=A_1(s)+A_2(s)$. It thus suffices to show that $A_1(s)$ and $A_2(s)$ both satisfy the functional equation
\begin{align*}
	A_i(n-s)=(-1)^{\frac{n-2t_2}{4}}2^{(s-\frac{n}{2})(\alpha_2-1)}A_i(s),\quad i=1,2.
\end{align*}
For $A_1(s)$ we can directly check that
\begin{align*}
	A_1(n-s)&=\frac{1-2^{(s-1)\ell_2}}{1-2^{s-1}}+\frac{(-1)^{\frac{n-2t_2}{4}}2^{(s-\frac{n}{2})(\alpha_2-1)}(1-2^{(n-1-s)\ell_2})}{1-2^{n-1-s}}=(-1)^{\frac{n-2t_2}{4}}2^{(s-\frac{n}{2})(\alpha_2-1)}A_1(s),
\end{align*}
satisfies the desired functional equation. For $A_2(s)$, if $\delta_{\bm{\lambda}}=0$ then this functional equation holds trivially. If $\delta_{\bm\lambda}=1$, then $(2\bm{\lambda})/2^{\ell_2}$ has only odd entries, implying that $\|2\bm{\lambda}\|^2/2^{2\ell_2}\equiv n\Mod{8}$. Since $n\equiv 2\Mod{4}$, this implies that $\alpha_2=2\ell_2+1$ and thus $2t_2=\|2\bm{\lambda}\|^2/2^{2\ell_2}\equiv n\Mod{8}$ so that $(-1)^{\frac{n-2t_2}{4}}=1$. One then easily checks that in this case $A_2(s)=2^{(n-1-s)\ell_2}$ satisfies that 
$$(-1)^{\frac{n-2t_2}{4}}2^{(s-\frac{n}{2})(\alpha_2-1)}A_2(s)=2^{(2s-n)\ell_2+(n-1-s)\ell_2}=2^{(s-1)\ell_2}=A_n(n-s),
$$
which is as desired.

Next, assume $n\equiv 4\Mod{8}$. In this case in view of the first relation in \eqref{equ:relazep2}, it suffices to show that $\cZ_n^{(2)}(s;\bm{\lambda})$ satisfies the functional equation 
\begin{align*}
	\cZ_n^{(2)}(n-s;\bm{\lambda})&=2^{(s-\frac{n}{2})(\alpha_2-2)}\cZ_n^{(2)}(s;\bm{\lambda}).
\end{align*}
Now note that in this case by \propref{prop:efactocomp2}
\begin{align*}
	\cZ_n^{(2)}(s;\bm{\lambda})&=\frac{1-2^{(n-1-s)\ell_2}}{1-2^{n-1-s}}+\delta_{\bm\lambda} 2^{(n-1-s)\ell_2}-\frac{2^{\frac{n}{2}-s}(1-2^{(n-1-s)(\ell_2-\delta_{\alpha_2,2\ell_2})})}{(1-2^{\frac{n}{2}-s})(1-2^{n-1-s})}\\
	&+\frac{2^{(\frac{n}{2}-s)(\alpha_2-2)+1}(1-2^{\frac{n}{2}-1-s})(1-2^{(s-1)(\ell_2-\delta_{\alpha_2,2\ell_2})})}{(1-2^{\frac{n}{2}-s})(1-2^{s-1})}.
\end{align*}
If $\alpha_2>2\ell_2$ then $\delta_{\alpha_2,2\ell_2}=0$ and
\begin{align*}
	\cZ_n^{(2)}(s;\bm{\lambda})&=\tfrac{(1-2^{\frac{n}{2}+1-s})(1-2^{(n-1-s)\ell_2})}{(1-2^{\frac{n}{2}-s})(1-2^{n-1-s})}
	+\tfrac{2^{(\frac{n}{2}-s)(\alpha_2-2)+1}(1-2^{\frac{n}{2}-1-s})(1-2^{(s-1)\ell_2})}{(1-2^{\frac{n}{2}-s})(1-2^{s-1})}+\delta_{\bm\lambda} 2^{(n-1-s)\ell_2}.
\end{align*}
Then doing similar arguments as in the previous case we can see that in this case $\cZ_n^{(2)}(s;\bm{\lambda})$ satisfies the above desired functional equation. Here when $\delta_{\bm{\lambda}}=1$ we need to use the fact that $\alpha_2=2\ell_2+2$. When $\alpha_2=2\ell_2$ (in particular, $\delta_{\bm{\lambda}}=0$) by rewriting $\frac{1-2^{(n-1-s)\ell_2}}{1-2^{n-1-s}}=\frac{1-2^{(n-1-s)(\ell_2-1)}}{1-2^{n-1-s}}+2^{(n-1-s)(\ell_2-1)}$ we similarly see that in this case
\begin{align*}
	\cZ_n^{(2)}(s;\bm{\lambda})&=\tfrac{(1-2^{\frac{n}{2}+1-s})(1-2^{(n-1-s)(\ell_2-1)})}{(1-2^{\frac{n}{2}-s})(1-2^{n-1-s})}
	+\tfrac{2^{(\frac{n}{2}-s)(\alpha_2-2)+1}(1-2^{\frac{n}{2}-1-s})(1-2^{(s-1)(\ell_2-1)})}{(1-2^{\frac{n}{2}-s})(1-2^{s-1})}+ 2^{(n-1-s)(\ell_2-1)},
\end{align*}
from which one can similarly check the desired functional equation for $\cZ_n^{(2)}(s;\bm{\lambda})$.

Finally, we treat the case when $n$ is odd. If $\alpha_2$ is odd, in view of the relation $\e_n^{(2)}(s;\bm{\lambda})=(1-2^{n-1-2s})^{-1}\cZ_n^{(2)}(s;\bm{\lambda})$ (cf. \rmkref{rmk:qformula}) and \propref{prop:efactocomp2} we see that in this case
\begin{align*}
	\e_n^{(2)}(s;\bm{\lambda})&=\frac{(1+\mathfrak{s}_n2^{\frac{n+1}{2}-s})(1-2^{(n-1-s)\ell_2})}{(1+\mathfrak{s}_n2^{\frac{n-1}{2}-s})(1-2^{n-2s})(1-2^{n-1-s})}-\frac{\mathfrak{s}_n2^{(\frac{n}{2}-s)(\alpha_2-2)+\frac12}(1-2^{(s-1)\ell_2})}{(1-2^{n-2s})(1-2^{s-1})},
\end{align*}
from which one directly checks the desired functional equation. Next, assume $\alpha_2$ is even and $t_2\not\equiv n\Mod{4}$ so that $b_{n-t_2}=0$ and  $(-1)^{\frac{n+t_2}{2}}=1$ (since both $n$ and $t_2$ are odd). Note also that in this case $\delta_{\bm{\lambda}}=0$ since otherwise we would have $\alpha_2=2\ell_2$ and $t_2=\|2\bm{\lambda}\|^2/2^{2\ell_2}\equiv n\Mod{8}$, contradicting the assumption that $t_2\not\equiv n\Mod{4}$. Thus applying \propref{prop:efactocomp2} we have
\begin{align*}
	\cZ_n^{(2)}(s;\bm{\lambda})&=\frac{1-2^{(n-1-s)\ell_2}}{1-2^{n-1-s}}+\frac{\mathfrak{s}_n2^{\frac{n-1}{2}-s}(1-2^{(n-1-s)(\ell_2-\delta_{\alpha_2, 2\ell_2})})}{(1-2^{n-2s})(1-2^{n-1-s})}\\
	&-\frac{\mathfrak{s}_n2^{(\frac{n}{2}-s)(\alpha_2-1)-\frac{1}{2}}(1-2^{(s-1)\ell_2})}{1-2^{s-1}}-\frac{\mathfrak{s}_n2^{(\frac{n}{2}-s)(\alpha_2-1)-\frac{1}{2}}(1-2^{(s-1)(\ell_2-\delta_{\alpha_2, 2\ell_2})})}{(1-2^{n-2s})(1-2^{s-1})}.
\end{align*}
Now applying the identity 
$\frac{1-q^{\ell}}{1-q}=\frac{1-q^{\ell-\delta_{\alpha,2\ell}}}{1-q}+\delta_{\alpha,2\ell}q^{\ell-1}$ for $q=2^{n-1-s}, 2^{s-1}$, and doing some further simplifications we get
\begin{align*}
	\cZ_n^{(2)}(s;\bm{\lambda})&=\tfrac{(1+\mathfrak{s}_n2^{\frac{n+1}{2}-s})(1-\mathfrak{s}_n2^{\frac{n-1}{2}-s})(1-2^{(n-1-s)(\ell_2-\delta_{\alpha_2,2\ell_2})})}{(1-2^{n-2s})(1-2^{n-1-s})}-\tfrac{\mathfrak{s}_n2^{(\frac{n}{2}-s)(\alpha_2-1)+\frac{1}{2}}(1-2^{n-1-2s})(1-2^{(s-1)(\ell_2-\delta_{\alpha_2,2\ell_2})})}{(1-2^{n-2s})(1-2^{s-1})}\\
	&+\delta_{\alpha_2,2\ell_2}2^{(n-1-s)(\ell_2-1)}(1-\mathfrak{s}_n 2^{\frac{n-1}{2}-s}).
\end{align*}
This, together with the relation $\e_n^{(2)}(s;\bm{\lambda})=(1-2^{n-1-2s})^{-1}\cZ_n^{(2)}(s;\bm{\lambda})$ (cf. \rmkref{rmk:qformula}), implies that
\begin{align*}
	\e_n^{(2)}(s;\bm{\lambda})&=\tfrac{(1+\mathfrak{s}_n2^{\frac{n+1}{2}-s})(1-2^{(n-1-s)(\ell_2-\delta_{\alpha_2,2\ell_2})})}{(1+\mathfrak{s}_n2^{\frac{n-1}{2}-s})(1-2^{n-2s})(1-2^{n-1-s})}-\tfrac{\mathfrak{s}_n 2^{(\frac{n}{2}-s)(\alpha_2-1)+\tfrac12}(1-2^{(s-1)(\ell_2-\delta_{\alpha_2,2\ell_2})})}{(1-2^{n-2s})(1-2^{s-1})}+\tfrac{\delta_{\alpha_2,2\ell_2}2^{(n-1-s)(\ell_2-1)}}{1+\mathfrak{s}_n 2^{\frac{n-1}{2}-s}}, 
\end{align*} 
from which one checks the desired functional equation. It now remains to treat the case when $\alpha_2$ is even and $t_2\equiv n\Mod{4}$ (so that $b_{n-t_2}=(-1)^{\frac{n-t_2}{4}}4$ and $(-1)^{\frac{n+t_2}{2}}=-1$). We further split the discussion into two cases. First we assume $t_2\equiv n+4\Mod{8}$ so that $b_{n-t_2}=-4$ and $\e_n^{(2)}(s;\bm{\lambda})=(1-\mathfrak{s}_n 2^{\frac{n-1}{2}-s})^{-1}\cZ_n^{(2)}(s;\bm{\lambda})$ (cf. \rmkref{rmk:qformula}). From the above discussion we also know that in this case $\delta_{\bm{\lambda}}=0$. Then similarly we can deduce that in this case
\begin{align*}
	\e_n^{(2)}(s;\bm{\lambda})&=\frac{(1+\mathfrak{s}_n2^{\frac{n+1}{2}-s})(1-2^{(n-1-s)(\ell_2-\delta_{\alpha_2,2\ell_2})})}{(1-2^{n-2s})(1-2^{n-1-s})}-\frac{2^{(\frac{n}{2}-s)\alpha_2}(1+\mathfrak{s}_n2^{\frac{n+1}{2}-s})(1-2^{(s-1)(\ell_2-\delta_{\alpha_2,2\ell_2})})}{(1-2^{n-2s})(1-2^{s-1})}\\
	&+\delta_{\alpha_2,2\ell_2}2^{(n-1-s)(\ell_2-1)}(1+\mathfrak{s}_n 2^{\frac{n+1}{2}-s}),
\end{align*}
from which one can also check the desired functional equation. Finally we assume $t_2\equiv n\Mod{8}$ so that $b_{n-t_2}=4$ and again by \rmkref{rmk:qformula}, $\e_n^{(2)}(s;\bm{\lambda})=(1+\mathfrak{s}_n 2^{\frac{n-1}{2}-s})^{-1}\cZ_n^{(2)}(s;\bm{\lambda})$. Doing similar calculations we get in this case
\begin{align*}
	\e_n^{(2)}(s;\bm{\lambda})&=\tfrac{(1+\mathfrak{s}_n2^{\frac{n+1}{2}-s})(1-\mathfrak{s}_n2^{\frac{n-1}{2}-s})(1-2^{(n-1-s)(\ell_2-\delta_{\alpha_2,2\ell_2})})}{(1+\mathfrak{s}_n2^{\frac{n-1}{2}-s})(1-2^{n-2s})(1-2^{n-1-s})} +\tfrac{2^{(\frac{n}{2}-s)\alpha_2}(1-\mathfrak{s}_n2^{\frac{n+1}{2}-s})(1-2^{(s-1)(\ell_2-\delta_{\alpha_2,2\ell_2})})}{(1-2^{n-2s})(1-2^{s-1})}\\ 
	&+\frac{\delta_{\alpha_2, 2\ell_2}2^{(n-1-s)(\ell_2-1)}(1+\delta_{\bm{\lambda}}2^{n-1-s}+2^{n-2s}+\mathfrak{s}_n 2^{\frac{n-1}{2}-s})}{1+\mathfrak{s}_n 2^{\frac{n-1}{2}-s}},
\end{align*}
from which one checks the desired functional equation. This finishes the proof.
\end{proof}

In order to prove \eqref{equ:crilinebd}, we collect the following estimates on the size of these local $\e$-factors on the line $\Re(s)=\frac{n}{2}$.	
\begin{Prop}\label{prop:estsizeepfac}
Assume $n$ is even. For any prime $p\mid 2\|2\bm{\lambda}\|^2$ and for any $t\in \R$ we have
\begin{align*}
|\e_n^{(p)}(\tfrac{n}{2}+it;\bm{\lambda})|\leq (1+\delta_{p,2})(\alpha_p+1)(\ell_p+1)p^{(\frac{n}{2}-1)\ell_p}.
\end{align*}
\end{Prop}
\begin{proof}
First assume $p$ is odd and we see from \rmkref{rmk:efacaritfun} that for $s=\frac{n}{2}+it$
\begin{align*}
\left|\e_n^{(p)}(s;\bm{\lambda})\right|&\leq \sum_{j=0}^{\ell_p}p^{(\frac{n}{2}-1)j}\sum_{l=0}^{\alpha_p-2j}\left|p^{-itl}\right|\leq (\alpha_p+1)(\ell_p+1)p^{(\frac{n}{2}-1)\ell_p},
\end{align*}
as claimed.
Next, we treat the case when $p=2$. In view of the first two equations in \eqref{equ:relazep2} and the estimate $\frac12\leq |1-2^{\frac{n}{2}-1-s}|\leq \frac32$ for $\Re(s)=\frac{n}{2}$, it suffices to show that 
\begin{align*}
|\cZ_n^{(2)}(\tfrac{n}{2}+it;\bm{\lambda})|\leq (\alpha_2+1)(\ell_2+1)2^{(\frac{n}{2}-1)\ell_2}.
\end{align*}
If $\ell_2=0$, then \propref{prop:efactocomp2} implies $\cZ_n^{(2)}(s;\bm{\lambda})=-1$ and the above bound holds trivially. Now we assume $\ell_2\geq 1$. 
When $n\equiv 2\Mod{4}$, we have by \propref{prop:efactocomp2} that for $s=\tfrac{n}{2}+it$,
\begin{align*}
|\cZ_n^{(2)}(s;\bm{\lambda})|&\leq \left|\sum_{j=0}^{\ell_2-1}2^{(n-1-s)j}\right|+2^{(\frac{n}{2}-1)\ell_2}+\left|\sum_{i=0}^{\ell_2-1}2^{(s-1)j}\right|\\
&\leq (2\ell_2+1)2^{(\frac{n}{2}-1)\ell_2}< (\alpha_2+1)(\ell_2+1)2^{(\frac{n}{2}-1)\ell_2}.
\end{align*}
Next, when $n\equiv 0\Mod{4}$ we have again by \propref{prop:efactocomp2} that for $s=\frac{n}{2}+it$ (denoting by $\tilde{\ell}_2:=\ell_2-\delta_{\alpha,2\ell{_2}}$),
\begin{align*}
|\cZ_n^{(2)}(s;\bm{\lambda})|&\leq (\ell_2+1) 2^{(\frac{n}{2}-1)\ell_2}+\left|\tfrac{1-2^{(n-1-s)\tilde{\ell}_2}}{(1-2^{\frac{n}{2}-s})(1-2^{n-1-s})}
-\tfrac{2^{(\frac{n}{2}-s)(\alpha_2-3)+1}(1-2^{\frac{n}{2}-1-s})(1-2^{(s-1)\tilde{\ell}_2})}{(1-2^{\frac{n}{2}-s})(1-2^{s-1})}\right|.
\end{align*}
If $\tilde{\ell}_2=0$, that is, $\alpha_2=2\ell_2=2$, then the terms in the above absolute value sign vanish and we have
\begin{align*}
|\cZ_n^{(2)}(\tfrac{n}{2}+it;\bm{\lambda})|&\leq (\ell_2+1)2^{(\frac{n}{2}-1)\ell_2}<(\alpha_2+1)(\ell_2+1)2^{(\frac{n}{2}-1)\ell_2}.
\end{align*}
Otherwise we have
\begin{align*}
&\frac{1-2^{(n-1-s)\tilde{\ell}_2}}{(1-2^{\frac{n}{2}-s})(1-2^{n-1-s})}
-\frac{2^{(\frac{n}{2}-s)(\alpha_2-3)+1}(1-2^{\frac{n}{2}-1-s})(1-2^{(s-1)\tilde{\ell}_2})}{(1-2^{\frac{n}{2}-s})(1-2^{s-1})}\\
&=\sum_{j=0}^{\tilde{\ell}_2-1}2^{(n-1-s)j}\frac{1-2^{(\frac{n}{2}-s)(\alpha_2-3-2j)+1}(1-2^{\frac{n}{2}-1-s})}{1-2^{\frac{n}{2}-s}}\\
&=\sum_{j=0}^{\tilde{\ell}_2-1}2^{(n-1-s)j}\left(\sum_{l=0}^{\alpha_2-4-2j}2^{(\frac{n}{2}-s)j}-2^{(\frac{n}{2}-s)(\alpha_2-3-2j)}\right).
\end{align*}
From this we have that when $\tilde{\ell}_2\geq 1$, 
\begin{align*}
|\cZ_n^{(2)}(\tfrac{n}{2}+it;\bm{\lambda})|&\leq (\ell_2+1) 2^{(\frac{n}{2}-1)\ell_2}+(\alpha_2-2)\ell_2 2^{(\frac{n}{2}-1)\ell_2}<(\alpha_2+1)(\ell_2+1)2^{(\frac{n}{2}-1)\ell_2},
\end{align*}
finishing the proof.
\end{proof}

We now give the 
\begin{proof}[Proof of \thmref{thm:funequbofo}]
We first treat the case when $n$ is even. In view of the relation \eqref{equ:efadeq1case} and using the functional equations satisfied by $\e_n^{(p)}(s;\bm{\lambda})$ for primes $p$ dividing $2\|2\bm{\lambda}\|^2$ (see \propref{prop:lononzeconterop}) we get 
\begin{align*}
\e_n(n-s;\bm{\lambda})&=\e_n(s;\bm{\lambda})
\left\{\begin{array}{ll}
				 \frac{2^{(s-\frac{n}{2})(\alpha_2-2)}(1-2^{\frac{n}{2}-1-s})}{1-2^{s-\frac{n}{2}-1}}\prod'_{p\mid \|2\bm{\lambda}\|^2}p^{(s-\frac{n}{2})\alpha_p} & n\equiv 4\Mod{8},\\
				(-1)^{\frac{n-2t_2}{4}}2^{(\alpha_2-1)(s-\frac{n}{2})}\prod'_{p\mid \|2\bm{\lambda}\|^2}\chi_{-4}(p)^{\alpha_p}p^{(s-\frac{n}{2})\alpha_p} & n\equiv 2\Mod{4}.
				\end{array}\right.
\end{align*} 
When $n\equiv 4\Mod{8}$ the above functional equation further reads as (using that $\|2\bm{\lambda}\|^2=\prod_{p\mid 2\|2\bm{\lambda}\|^2}p^{\alpha_p}$)
\begin{align*}
\e_n(n-s;\bm{\lambda})=\frac{2^{n-2s}(1-2^{\frac{n}{2}-1-s})}{1-2^{s-\frac{n}{2}-1}}\|2\bm{\lambda}\|^{2s-n}\e_n(s;\bm{\lambda})
\end{align*}
which is as desired. When $n\equiv 2\Mod{4}$ the above functional equation further reads as (noting that $\chi_{-4}(p)=\left(\frac{-1}{p}\right)$ for any odd prime $p$)
\begin{align*}
\e_n(n-s;\bm{\lambda})&=(-1)^{\frac{n-2t_2}{4}}\text{$\prod$}'_{p\mid \|2\bm{\lambda}\|^2}\left(\tfrac{-1}{p}\right)^{\alpha_p}2^{\frac{n}{2}-s}\|2\bm{\lambda}\|^{2s-n}\e_n(s;\bm{\lambda})\\
&=(-1)^{\frac{n-2}{4}}(-1)^{\frac{1-t_2}{2}}\text{$\prod$}'_{ p\mid \|2\bm{\lambda}\|^2}\left(\tfrac{-1}{p}\right)^{\alpha_p}2^{\frac{n}{2}-s}\|2\bm{\lambda}\|^{2s-n}\e_n(s;\bm{\lambda}).
\end{align*}
But $t_2=\|2\bm{\lambda}\|^2/2^{\alpha_2}=\prod'_{p\mid \|2\bm{\lambda}\|^2}p^{\alpha_p}$ and thus
\begin{align*}
\text{$\prod$}'_{p\mid \|2\bm{\lambda}\|^2}\left(\tfrac{-1}{p}\right)^{\alpha_p}=\left(\tfrac{-1}{\prod'_{p\mid \|2\bm{\lambda}\|^2}p^{\alpha_p}}\right)=\left(\tfrac{-1}{t_2}\right)=(-1)^{\frac{1-t_2}{2}},
\end{align*}
implying that 
\begin{align*}
\e_n(n-s;\bm{\lambda})&=\mathfrak{s}_n2^{\frac{n}{2}-s}\|2\bm{\lambda}\|^{2s-n}\e_n(s;\bm{\lambda}),
\end{align*}
as desired. Finally, we treat the case when $n$ is odd and as before $D=(-1)^{\frac{n-1}{2}}\|2\bm{\lambda}\|^2$. Similarly, from the functional equations proved in \propref{prop:lononzeconterop} we get 
\begin{align*}
\e_n(n-s;\bm{\lambda})&=\e_n(s;\bm{\lambda})\left\{\begin{array}{ll}
\frac{\mathfrak{s}_n2^{(s-\frac{n}{2})(\alpha_2-4)+\frac12}(1+\mathfrak{s}_n2^{\frac{n-1}{2}-s})}{1+\mathfrak{s}_n 2^{\frac{n+1}{2}-s}}\prod'_{p\mid \|2\bm{\lambda}\|^2}p^{(s-\frac{n}{2})\tilde{\alpha}_p}  & \text{$\alpha_2$ odd},\\
				 \frac{\mathfrak{s}_n2^{(s-\frac{n}{2})(\alpha_2-1)+\frac12}(1+\mathfrak{s}_n2^{\frac{n-1}{2}-s})}{1+\mathfrak{s}_n2^{\frac{n+1}{2}-s}}\prod'_{p\mid \|2\bm{\lambda}\|^2}p^{(s-\frac{n}{2})\tilde{\alpha}_p} & \text{$\alpha_2$ even}\ \text{and}\ t_2\equiv n\Mod{4},\\
				\frac{\mathfrak{s}_n2^{(s-\frac{n}{2})(\alpha_2-3)+\frac12}(1+\mathfrak{s}_n2^{\frac{n-1}{2}-s})}{1+\mathfrak{s}_n 2^{\frac{n+1}{2}-s}}\prod'_{p\mid \|2\bm{\lambda}\|^2}p^{(s-\frac{n}{2})\tilde{\alpha}_p}  & \text{$\alpha_2$ even}\ \text{and}\ t_2\not\equiv n\Mod{4}.
				\end{array}\right.
\end{align*}
Here $\tilde{\alpha}_p=\alpha_p-1$ if $\alpha_p$ is odd and $\tilde{\alpha}_p=\alpha_p$ if $\alpha_p$ is even. From the above formulas, together with the modulus formula of $\chi_{D}$ (see \eqref{equ:primcharmodu}) one can then check the desired functional equation for  $\e_n(s;\bm{\lambda})$. 
This finishes the proof of the functional equation. 

For \eqref{equ:crilinebd} write $\bm{\lambda}=a\bm{\lambda}'$ with $\bm{\lambda}'\in \Lambda_{\rm pr}^{*}$ and $a\in \N$. Note that either $\bm{\lambda}'\in \Z^n_{\rm pr}$ or $2\bm{\lambda}'\in \Z^n_{\rm pr}$. From this we see that for any odd prime $p\mid \|2\bm{\lambda}\|^2$, 
$$
\ell_p=\nu_p(\gcd(2\bm{\lambda}))=\nu_p(a)+\nu_p(\gcd(2\bm{\lambda}'))=\nu_p(a).
$$
Similarly, we have $0\leq \ell_2-\nu_2(a)\leq 1$, implying that $
a\asymp \text{$\prod$}_{p\mid 2\|2\bm{\lambda}\|^2} p^{\ell_p}$.
Now, by \propref{prop:estsizeepfac} and the relation \eqref{equ:efadeq1case} we can bound
\begin{align*}
	|\e_n(\tfrac{n}{2}+it;\bm{\lambda})|&\ll \prod_{p\mid 2\|2\bm\lambda\|^2}(\alpha_p+1)(\ell_p+1)p^{(\frac{n}{2}-1)\ell_p}\ll_{\e}a^{\frac{n}{2}-1}\prod_{p\mid 2\|2\bm{\lambda}\|^2}p^{\frac{\alpha_p\e}{2}}\asymp a^{\frac{n}{2}-1+\e}\|\bm{\lambda}'\|^{\e},
\end{align*}
proving \eqref{equ:crilinebd}. This then finishes the proof.
\end{proof}

\bibliographystyle{alpha}
\bibliography{DKbibliog}

\end{document}